\documentclass{article}

\usepackage{amsmath,amssymb,amsthm}
\usepackage[mathscr]{eucal}
\usepackage{a4wide}
\usepackage{graphicx}
\usepackage{psfrag}
\usepackage[all]{xy}

\title{From Fractal Groups to Fractal Sets}

\author{Laurent Bartholdi\footnote{Supported by the ``Swiss National
    Foundation for Scientific Research'' grant 83R-064282.}
  \and Rostislav Grigorchuk\footnote{Supported by the program
    ``Leading Scientific Schools of the Russian Federation'' project N
    00-15-96107.}
  \and Volodymyr Nekrashevych\footnote{Supported by the
    ``A.~von~Humboldt Foundation''.}}

\newtheorem{theorem}{Theorem}[section]
\newtheorem{proposition}[theorem]{Proposition}
\newtheorem{corollary}[theorem]{Corollary}
\newtheorem{lemma}[theorem]{Lemma}

\theoremstyle{definition}
\newtheorem{defi}{Definition}[section]

\newtheorem*{examp}{Example}

\newenvironment{reformulate}[2]{\begingroup
  \expandafter\def\csname the#1\endcsname{#2}\def\envname{#1}\begin{\envname}}%
  {\end{\envname}\addtocounter{\envname}{-1}\endgroup}

\newcommand{\mapdown}[1]%
{\Big\downarrow\rlap{$\vcenter{\hbox{$\scriptstyle#1$}}$}}

\newcommand{\be}{\alpha}
\newcommand{\C}{{\mathbb{C}}}

\newcommand{\dom}{\mathop{\mathrm{Dom}}}
\newcommand{\en}{\omega}
\newcommand{\emp}{\emptyset}
\newcommand{\f}{\mathcal{F}}
\newcommand{\G}{\Gamma}
\newcommand{\IMG}{\mathop{IMG}}
\newcommand{\lims}[1]{\mathscr{J}_{#1}}
\newcommand{\nuke}{\mathcal{N}}
\newcommand{\On}{\mathcal{O}_d}
\newcommand{\OPhi}{\mathcal{O}_\Phi}
\newcommand{\prr}{\dashrightarrow}
\newcommand{\ran}{\mathop{\mathrm{Ran}}}
\newcommand{\redu}{\mathcal{A}_\Phi}
\newcommand{\scalar}[2]{\left\langle\left. #1\right| #2\right\rangle}
\newcommand{\shift}{\mathsf{s}}
\newcommand{\solen}[1]{\mathscr{S}_{#1}}
\newcommand{\Ss}{\Sigma}
\newcommand{\si}{\mathfrak{s}}
\newcommand{\til}{\mathscr{T}}
\newcommand{\UPhi}{\mathcal{F}(\Phi)}
\newcommand{\xo}{X^\omega}
\newcommand{\xmo}{X^{-\omega}}
\newcommand{\xs}{X^*}
\newcommand{\xz}{X^\Z}

\newcommand{\img}[1]{\mathop{\mathrm{IMG}}\left(#1\right)}

\newcommand\X{{\mathfrak X}}
\newcommand\lang{{\mathcal L}}
\newcommand\autom{{\mathcal A}}
\newcommand\aut{{\mathop{\textsf{Aut}}}}
\newcommand\stab{{\mathop{\textsf{St}}}}
\newcommand\rist{{\mathop{\textsf{Rist}}}}
\newcommand\N{{\mathbb N}}
\newcommand\Q{\mathbb{Q}}
\newcommand\Z{{\mathbb Z}}
\newcommand\R{{\mathbb R}}
\newcommand{\spec}{\mathop{\mathrm{spec}}}
\newcommand\hilb{{\mathcal H}}

\begin{document}
\maketitle
\tableofcontents

\section{Introduction}
The idea of self-similarity is one of the most fundamental in the
modern mathematics. The notion of ``renormalization group'', which
plays an essential role in quantum field theory, statistical physics
and dynamical systems, is related to it. The notions of fractal and
multi-fractal, playing an important role in singular geometry, measure
theory and holomorphic dynamics, are also related.  Self-similarity
also appears in the theory of $C^*$-algebras (for example in the
representation theory of the Cuntz algebras) and in many other branches of
mathematics. Starting from 1980 the idea of self-similarity entered
algebra and began to exert great influence on asymptotic and geometric
group theory.

The aim of this paper is to present a survey of ideas, notions and
results that are connected to self-similarity of groups, semigroups
and their actions; and to relate them to the above-mentioned classical
objects. Besides that, our aim is to exhibit new connections of groups
and semigroups with fractals objects, in particular with Julia sets.

Let us review shortly some historical aspects of our research and list
its main subjects.

\subsection{Burnside groups}
The second named author has constructed in 1980~\cite{grigorchuk:80_en} two
Burnside groups which (especially the first one) played a decisive
role in the development of the idea of a self-similar group.
Originally the Grigorchuk group (let us denote it by $\mathbf{G}$) was defined
as a transformation group of the segment $[0,1]$ without dyadically
rational points. The generators were defined in a simple way as
permutations of subintervals. One of the main properties of this
action is the fact that if we restrict to an arbitrary dyadic
subinterval $I=[(k-1)/2^n,k/2^n]$ the action of its stabilizer
$\mathbf{G}_I=\mathsf{st}_\mathbf{G}(I)$, then the restriction will coincide with the
action of the group $\mathbf{G}$ on $[0, 1]$ (after the natural identification of $I$ with the whole
interval $[0, 1]$).

Another fundamental property is the fact that the action of the group
$\mathbf{G}$ is \emph{contracting}, i.e., that the canonical homomorphism
$\phi_I:\mathbf{G}_I\to \mathbf{G}$ contracts the length of the group elements by a
constant $\lambda>1$. Finally, the third fundamental property is the
\emph{branching} nature of the action.  This means that, up to
finite-index inclusions, the stabilizers of the partition into the
dyadic subintervals of the $n$th level are direct products of $2^n$
isomorphic groups and that the lattice of subnormal groups has a
tree-like structure.

Between 1983 and 1985 the second named author~\cite{grig:85a,grigorchuk:growth_en}, and
independently N.~Gupta and S.~Sidki~\cite{gupta-s:ext},
constructed other examples of similar groups and established
other their important properties. It became clear that these examples
are related to some big classes of groups.  Namely, they belong to the
class of \emph{finitely automatic groups} (this was noted for the first time
by Ju.~I.~Merzlyakov~\cite{merzlyakov:grigorchuk}), and to the class
of \emph{branch groups}. The first class was defined in the early
1960's~\cite{hor:fin_en}, while the definition of the class of branch groups
was given in 1997 by R.~Grigorchuk in his talk at the St-Andrews
conference in Bath (see~\cite{grigorchuk:branch}). The methods used
in~\cite{gupta-sidki_group,grigorchuk:80_en} were new to the theory of
automata groups, and heralded the application of these groups to many
new problems in group theory --- see for
instance~\cite{bartholdi-g:lie}. 

One of the distinctive features of the branch groups is their actions
on spherically homogeneous rooted trees. Such trees appear naturally
in the study of unimodal transformations~\cite{bass}, in particular in
problems of renormalization and in the study of dynamical systems in a
neighborhood of Lyapunov-stable attractors~\cite{buescu}. A regular
rooted tree is an example of a geometric object most closely related
with the notion of self-similarity (every rooted subtree is naturally
isomorphic to the original tree). The actions of the groups
from~\cite{grigorchuk:80_en,grigorchuk:growth_en,grig85_en,gupta-sidki_group,gupsid2,bsv:jns}
on the respective regular trees have self-similarity properties
similar to the self-similarity property of the action of the group $\mathbf{G}$
on the interval $[0,1]$.  These examples can be formalized in the
general notion of a self-similar set and a self-similar action.

In Section~\ref{s:sss} we discuss in detail the notion of a
self-similar set and then define the notion of a self-similar action
of a group and of an inverse semigroup. (The fact that self-similar
inverse semigroups appear naturally in connection with self-similar
sets, in particular, with Penrose tilings, was noted by
V.~Nekrashevych~\cite{nek:ssinv}.) One of the main sources of examples
of self-similar group actions are the actions generated by finite
automata. Groups of finite automata and their actions are discussed in
Section~\ref{s:automata}.

\subsection{Growth}
The second named author noted in 1983 that the growth of the group $\mathbf{G}$,
(which was initially introduced as a Burnside group, i.e., an infinite
finitely generated torsion group), is intermediate between polynomial
and exponential (and thus provides an answer to Milnor's
question~\cite{milnor:5603}). This observation was developed in
different directions~\cite{grigorchuk:growth_en}. Besides, it was recently proved
that many such groups, and in particular the Gupta-Sidki
groups, also have intermediate growth~\cite{bartholdi:interm}.
\emph{In fine}, intermediate growth of these groups follows from
strong contraction properties and branch structure.  It became clear
that the groups of such a type can be related with problems of fractal
geometry, in particular with the problems of computation of the Hausdorff
dimension. Later this conjecture was confirmed in different ways.

Most problems of computation of Hausdorff dimension (and dimensions of
other type) reduce to the problem of finding the degree of polynomial
growth or the base of exponential growth of some formal language over a finite alphabet.
Section~\ref{s:growth} is devoted to questions of growth of formal
languages, groups, semigroups, graphs and finite automata. In the same
section the notion of amenability, which also plays a role in the
theory of fractals, is considered.

Besides playing an important role in the study of growth of finitely
generated groups, self-similar groups also appear in the study of the
Hausdorff dimension of profinite groups. For instance, the profinite
completion of the group $\mathbf{G}$ has Hausdorff dimension
$5/8$~\cite{0964.20018,grigorchuk:branch}.

Another relation of self-similar groups to fractals was found
accidentally while studying the spectra of the
non-commutative dynamical systems generated by the actions of
self-similar groups (like the group $\mathbf{G}$ mentioned above, or the
Gupta-Sidki group)~\cite{bgr:spec,bartholdi-g:parabolic}.
First, it turned out that one has to use multi-dimensional rational
mappings for the solution of the spectral problem and to study
their invariant subsets (similar to the Julia sets). Secondly, the
spectra turn out to be the Julia sets of polynomial mappings of the
interval. Finally, amenability of the objects of polynomial growth (in
particular groups and graphs) imply coincidence of the spectra of the
above-mentioned dynamical systems with the spectra of the discrete
Laplace operator on the Schreier graphs of the self-similar groups
(where the Schreier graphs are defined with respect to the stabilizer
of an end of the tree).

\subsection{Schreier graphs}
The Schreier graphs themselves are also interesting objects of
investigation. They have polynomial growth in the case when the group action is
contracting (though often the degree of the growth is non-integral).
It was discovered that the Schreier graphs (defined
with respect to the stabilizers of different points of the tree
boundary) behave similarly to quasi-periodic tilings (like the Penrose
tilings). Often there exist uncountably many non-isomorphic Schreier
graphs in the bundle, while they all are locally isomorphic. Many
other analogies can be found, in particular, the inflation and
adjacency rules of the Penrose tilings have their counterparts in the
Schreier graphs of self-similar groups.  The finite Schreier graphs
(defined with respect to the stabilizers of the tree vertices) of
self-similar actions are often substitution graphs (here again an
analogy with substitution dynamical systems and L-systems appears).

Their limit spaces (for example, in the sense of M.~Gromov) often have a
fractal nature. Moreover, recently V.~Nekrashevych has introduced a
notion of an iterated monodromy group (i.m.g.) of a branched covering
and proved that the i.m.g.\ of a postcritically finite rational
mapping is contracting (and thus is generated by a finite automaton)
with the limit space homeomorphic to the Julia set of the mapping.
Probably in the future the iterated monodromy groups and their
Schreier graphs will play an important role in the holomorphic
dynamics and the methods of the asymptotic and geometric group theory
will be actively used in this part of mathematics. The basic
definitions, facts and some examples of iterated monodromy groups are
described in Section~\ref{s:img}.

Schreier graphs of self-similar groups have very interesting spectral
properties,  as was discovered in~\cite{bgr:spec} and~\cite{gr_zu:lamp}.  For instance,
in~\cite{bgr:spec}  the first examples of regular graphs with a Cantor spectrum are constructed,
while in~\cite{gr_zu:lamp}  the first example of a group with discrete spectral measure is given
(which solves a question of Atyiah, see~\cite{glsz:atiyah}).

 Some new examples of computations of the sepctra and an example from~\cite{bgr:spec}
 concerning the Gupta-Fabrikovsky group is considered in Section~\ref{s:spectrum}.

As we have noted above, the Schreier graphs of contracting
self-similar groups converge to nice fractal topological spaces. This
was formalized by V.~Nekrashevych~\cite{nek:lim} in the notion of the
limit space of a self-similar action. The limit space has a rich
self-similarity structure with which self-replicating tiling systems
and limit solenoids are related. The limit solenoid can be defined
also for self-similar inverse semigroups. The obtained constructions
agree with the well known notions of a self-affine tiling of the
Euclidean plane (which correspond in our situation to self-similar
actions of commutative groups). In the case of the iterated monodromy
group of a rational mapping of the complex sphere the limit space is
homeomorphic to the Julia set of the mapping.

\subsection{Virtual endomorphisms and $L$-presentations}
An important role in the study of self-similar groups is played by
\emph{virtual endomorphisms}, i.e., endomorphisms defined on a
subgroup of finite index.  Every self-similar action of a group
defines an associated virtual endomorphism of the group. Together with
some simple additional data the virtual endomorphism determines uniquely
the action.  In this way the self-similar actions can be interpreted
as abstract numeration systems on the groups with the virtual
endomorphism playing the role of the base. Such numeration systems are
natural generalizations of the usual numeration systems on the group
$\Z$. Self-similar actions are also associated with the well-known
numeration systems on the free abelian groups $\Z^n$
(see~\cite{vince:digtile}).  The respective tilings of Euclidean space
can also be interpreted in the terms of self-similar actions and
generalized to non-commutative groups.

I.~G.~Lysenok has obtained in 1985 a finite \emph{L-presentation} of
the Grigorchuk group $\mathbf{G}$, i.e.\ a presentation in which the defining
relations are obtained from a finite set of relations using iterated
application of a substitution $f$ over the alphabet of generators
(equivalently, $f$ can be viewed as an endomorphism of the respective
free group). Such presentations were obtained independently by
S.~Sidki for the Gupta-Sidki group~\cite{sidki:pres}, and later
L.~Bartholdi gave a universal method of constructing L-presentations
of branch groups using the virtual
endomorphisms~\cite{bartholdi:lpres}. The L-presentations are
convenient to construct embeddings of a group into a finitely
presented group (using only one HNN-extension).  If the L-presentation
is defined using a usual endomorphism, then this embedding preserves
amenability of the group.  This was used in~\cite{grigorchuk:notEG} to
construct a finitely presented amenable but not elementary amenable
group. On one hand, L-presentations are similar to L-systems, which
are well known in formal language theory.  On the other hand, they
have analogies with substitution dynamical systems~\cite{substdyn}.  The
substitution dynamical systems have a direct relation with fractals.

\subsection{Boundaries}
One of the sources of fractal sets are the various boundaries of
groups. There exist a great variety of different notions of boundary
of a group connected with different compactifications: Baily-Borel
compactification, Freudenthal boundary
(the space of the ends)~\cite{freudenthal:ueberenden,freudenthal:endenth}
boundary (see~\cite{woess:rw}), Poisson-Furstenberg boundary~\cite{furst:pois},
Gromov boundary~\cite{gro:hyperb},
Higson-Roe corona~\cite{dranishnikov-f:hrcorona},  Floyd
compactification~\cite{karlsson1,karlsson2,floyd:fursten}, etc.
There exists a rather general method to
construct a boundary of a finitely generated group based on the use of
the metrics (or uniform structures) satisfying the condition
$d(gx,gy)\to 0$ when $g$ tends to infinity (see, for
instance,~\cite{floyd:kleinian} where a partial case of such metrics
is considered).  As it was noted probably for the first time by
A.~S.~Mishchenko, the respective boundaries play an important role in
the topics related to such famous conjectures as the Novikov
conjecture or the Baum-Connes
conjecture~\cite{dranishnikov-f:hrcorona}.

An important class of metric spaces and finitely generated groups are
the Gromov-hyperbolic groups~\cite{gro:hyperb}. They possess a natural
boundary (the Gromov boundary) which is one of the most well studied boundaries.
It is known, that the boundary of a hyperbolic group can
be homeomorphic to the Sierpinski carpet, the Menger curve and to other
sets of a fractal nature. The action of a hyperbolic space on its
boundary is an example of a finitely presented dynamical system (the
notion belongs also to M.~Gromov, see~\cite{gro:hyperb}
and~\cite{curnpap:symb}). The boundary is a semi-Markovian space, and 
has a rich self-similarity structure.  Many Kleinian groups are word hyperbolic
and their limit sets are often homeomorphic to their Gromov
boundaries.

A natural generalization of the boundary of hyperbolic groups are the Dynkin sets of $\mathcal{D}$-groups,
defined by H.~Furstenberg~\cite{furst:dyn}. The hyperbolic groups and
their finitely generated subgroups, which are not virtually cyclic,
belong to the class of $\mathcal{D}$-groups.

As the third named author noted recently, a naturally-defined hyperbolic
graph can be associated with every contracting self-similar group
action; and its boundary is homeomorphic to its limit space.  A short
survey of notions and facts about hyperbolic spaces, groups and their
boundaries is presented in Sections~\ref{s:hyp},~\ref{s:fpds}.

The classical notion of self-similarity (of topological spaces) and
the notion of a self-similar action can be interpreted from a common
point of view using the notion of a Hilbert bimodule over a
$C^*$-algebra. Some $C^*$-algebras (for instance the Cuntz-Pimsner
algebras) associated with such self-similarity bimodules were studied
in~\cite{nek:bim}; however the study of some other algebras is only at
its beginning~\cite{GZ*}.  This is the topic of
Section~\ref{s:cstar}.

Finally, a word on ``fractals''. We use the term ``fractal group''
only in the title, and not in the text. The reason is that, just as
there is no unique generally accepted definition of a fractal set,
there is no definition of a fractal group. However, in some papers,
including the works of the authors, different variants of a definition
corresponding to the notion of recurrent action (see
Definition~\ref{defi:recur}) where proposed. Roughly speaking, fractal
groups are the groups acting self-similarly on a self-similar set and
such that their geometry, analysis and dynamics are related in some
way with fractal objects. We hope that numerous examples given in this
paper will convince the reader that we are on the right path toward
the notion of a fractal group.

\paragraph{Acknowledgements.} The authors are grateful to K.~Falconer,
H.~Furstenberg, T.~Giordano, P.~de~la~Harpe, A.~Henriques and V.~Jones
for useful discussions and
 interest.

The authors are also grateful to the Swiss National Science Foundation
and especially to P.~de~la~Harpe for invitations to Geneva University,
where a major part of the work on the article was made.

The third named author thanks R.~Pink for an invitation to ETH Z\"urich and
interesting discussions which eventually led to the definition
of the notion of an iterated monodromy group.

\section{Preliminary definitions}
\subsection{Spaces of words}
For a finite set $X$ (an \emph{alphabet}) we denote by
$\xs=\{x_1x_2\ldots x_n| x_i\in X, n\ge 0\}$ the set of all finite
words over the alphabet $X$, including the empty word $\emp$.  We have
$\xs=\cup_{n\geq 0} X^n$ (we put $X^0=\{\emp\}$).  If $v=x_1x_2\ldots
x_n\in X^n$ then $n$ is the length of the word $v$ and is written
$|v|$.  The product (concatenation) $v_1v_2$ of two words $v_1, v_2\in
\xs$ is naturally defined.

By $\xo$ we denote the set of all infinite unilateral sequences
(words) of the form $x_1x_2\ldots$, $x_i\in X$. If $v\in\xs$ and
$w\in\xo$, then the product $vw\in\xo$ is also naturally defined.

The set $\xo$ is equipped with the topology of direct product of the
discrete finite sets $X$. The basis of open sets in this topology is
the collection of all \emph{cylindrical sets}
\[a_1a_2\ldots a_n\xo=\{x_1x_2\ldots\in \xo | x_i=a_i, 1\le i\le n\}\]
where $a_1a_2\ldots a_n$ runs through $\xs$. The space $\xo$ is
totally disconnected and homeomorphic to the Cantor set.

In a similar way we can introduce a topology on the set $\xo\cup\xs$
taking a basis of open sets $\left\{v\xs\cup v\xo:v\in\xs\right\}$,
where $v\xs\cup v\xo$ is the set of all words (finite and infinite)
beginning with $v$.

The topological space $\xo\cup\xs$ is compact, the set $\xo$ is closed
in it and the set $\xs$ is a dense subset of isolated points.

The \emph{shift} on the space $\xo$ is the map $\shift:\xo\to\xo$,
that deletes the first letter of the word:
\[\shift(x_1x_2\ldots)=x_2x_3\ldots.\]

The space $\xo$ is also called the \emph{full one-sided shift
  (space)}.

\begin{defi}
  A subset $\mathcal{F}\subseteq\xo$ is called a \emph{subshift (space)}
  if it is closed and invariant under the shift $\shift$, i.e., if
  $\shift(\mathcal{F})\subseteq\mathcal{F}$.
\end{defi}

\begin{defi}
  A subset $\mathcal{F}\subseteq\xo$ is a \emph{subshift (space) of  finite type} if
there exists a number $m\in\N$ and a subset $A\subset X^m$ of \emph{admissible words},
such that a word $w\in\xo$ belongs to $\mathcal{F}$ if and
every its subword of length $m$ belongs to $A$.
 If $m=2$ then the subshift $\mathcal{F}$ is also called a \emph{topological Markov chain}.
\end{defi}

It is easy to prove that every subshift of finite type is a subshift
space.

\subsection{Graphs}
We will use two versions of the notion of a graph. The first is the most general one (directed graphs with
loops and multiple edges). It will be used to construct Moore diagrams of automata
(Subsection~\ref{ss:automata}), structural graphs of iterated function systems (Definition~\ref{defi:gdifs}),
Schreier graphs of groups (Definition~\ref{ss:schdef}). All the other graphs, which will appear in our
paper, are defined using a more classical notion of a (simplicial) graph, i.e., a non-directed graph without loops
and multiple edges.

\begin{defi}
A \emph{graph} $\G$ is defined by a set of \emph{vertices} $V(\G)$, a set of \emph{edges}
(\emph{arrows}) $E(\G)$ and maps $\be, \en:E(\G)\to V(\G)$. Here $\be(e)$
is the \emph{beginning} (or \emph{source}) of the edge $e$ and $\en(e)$
is its \emph{end} (or \emph{range}).
\end{defi}

Two vertices $v_1, v_2$ are \emph{adjacent} if there exists an edge
$e$ such that $v_1=\be(e)$ and $v_2=\en(e)$ or $v_2=\be(e)$ and
$v_1=\en(e)$.  Then we say that the edge $e$ \emph{connects} the
vertices $v_1$ and $v_2$.

The \emph{(edge-)labeled graph} is a graph together with a map
$l:E(\G)\to S$, which assigns a \emph{label} $l(e)\in S$ to every edge
of the graph. Here $S$ is a given \emph{label set}.

A \emph{morphism} of graphs $f:\G_1\to\G_2$, is a pair of maps
$f_v:V(\G_1)\to V(\G_2), f_e:E(\G_1)\to E(\G_2)$ such that
\begin{eqnarray*}
  \be\left(f_e\left(e\right)\right) &=& f_v\left(\be\left(e\right)\right)\\
  \en\left(f_e\left(e\right)\right) &=& f_v\left(\en\left(e\right)\right)
\end{eqnarray*}
for all $e\in E(\G_1)$. A \emph{morphism of labeled graphs} is a
morphism of graphs preserving the labels of the edges.

A \emph{path} in a graph $\G$ is a sequence of edges $e_1 e_2\ldots
e_n$, with $\en(e_i)=\be(e_{i+1})$ for every $1\leq i\leq n-1$. The
vertex $\be(e_1)$ is called the \emph{beginning of the path,} and the
vertex $\en(e_n)$ is its \emph{end}. The number $n$ is called the
\emph{length} of the path. In the similar way define \emph{infinite to
  the right} paths $e_1e_2\ldots$, \emph{infinite to the left} paths
$\ldots e_2e_1$ and the \emph{bi-infinite} paths $\ldots
e_{-1}e_0e_1e_2\ldots$.

\begin{defi}
A \emph{simplicial graph} $\G$ is defined by its set of vertices $V(\G)$ and a
set of edges $E(\G)$, where every edge is a set $\{v_1, v_2\}$ of two
different vertices $v_1, v_2\in V(\G)$.
\end{defi}

Thus, a simplicial graph is not directed and does not have loops or multiple edges.

If $\{v_1, v_2\}\in E$ then we say that the vertices $v_1$ and $v_2$
are \emph{adjacent}, and that the edge $e=\{v_1, v_2\}$
\emph{connects} the vertices.

A \emph{morphism} of simplicial graphs $f:\G_1\to\G_2$ is a map $V(\G_1)\to V(\G_2)$ which
preserves the adjacency of the vertices.

The degree $\deg v$ of a vertex $v$ of a simplicial graph is the number
of edges to which it belongs.

If $\G$ is a graph, then its \emph{associated simplicial graph} is the
simplicial graph with the same set of vertices, which contains an edge
$e=\{v_1, v_2\}$ if and only if the vertices $v_1$ and $v_2$ were
adjacent in the original graph and $v_1\neq v_2$.

A \emph{path} in a simplicial graph $\G$ is a sequence of vertices $v_1 v_2\ldots
v_{n+1}$, with $\{v_i, v_{i+1}\}\in E(\G)$ for every $1\leq i\leq n$. The
vertex $v_1$ is called the \emph{beginning of the path,} and the
vertex $v_{n+1}$ is its \emph{end}. The number $n$ is called the
\emph{length} of the path.

A \emph{geodesic path}, connecting the vertices $v_1$ and $v_2$ is a
path of minimal length, whose beginning and end are $v_1$ and $v_2$
respectively.

The length of a geodesic path connecting the edges $v_1$ and $v_2$ is
called their \emph{distance} and is written $d(v_1, v_2)$.  We define
$d(v, v)=0$. If the graph $\G$ is connected, then the distance $d(v_1, v_2)$ is defined
for every pair of vertices $v_1, v_2\in V(\G_1)$ and is called the \emph{natural} (or
\emph{combinatorial}) metric on the graph.

The distance between two vertices of a non-simplicial graph is defined
as the distance between the vertices in the associated simplicial graph.

For a graph $\G$, a vertex $v\in V(\G)$, and $r\in\N$, the
\emph{ball} $B(v, r)$ of \emph{radius} $r$ with \emph{center} at the
point $v$ is the set of the vertices $\{u\in V: d(v, u)\leq r\}$.

A graph $\G$ is \emph{locally finite} if for every vertex $v$ the ball
$B(v, 1)$ is finite. If the graph is locally finite, then every ball
$B(v, r)$ is finite.

\section{Self-similar sets and (semi)group actions}
\label{s:sss}
\subsection{Self-similar sets}
\label{ss:selfsims}
\begin{defi}
  \label{defi:gdifs}
  A \emph{graph-directed iterated function system} with
  \emph{structural graph} $\Gamma=\langle V, E, \be, \en\rangle$ is a
  finite collection of sets $\{F_v\}_{v\in V}$ together with a
  collection of injective maps $\{\phi_e:F_{\be(e)}\to
  F_{\en(e)}\}_{e\in E}$, such that for every $v\in V$
  \begin{equation}
    \label{eq:defifs}
    F_v=\bigcup_{\en(e)=v} \phi_e(F_{\be(e)}).
  \end{equation}
  
  If the structural graph contains only one vertex $v$, then we say
that we have an \emph{iterated function system} on the set $F=F_v$.
  
  If the sets $F_v$ are subsets of a common set $F$ and $F=\cup_{v\in
 V} F_v$, then the graph-directed iterated function system will
also be called a \emph{self-similarity structure} on $F$. A set with a self-similarity structure on it is called
\emph{self-similar}.
  
  A set with a fixed self-similarity structure is called
  \emph{self-similar}.
\end{defi}

Note that the structural graph of an iterated function system
on a single set is the graph with a single vertex and $|E|$ loops.

\begin{defi}
  A \emph{topological} graph-directed iterated function system is a
system $\left(\{F_v\}, \{\phi_e\}\right)$ for
which the sets $F_v$ are compact Hausdorff topological spaces and
the injections $\phi_e$ are continuous.
\end{defi}

In the case of the classical notion of a graph-directed iterated
function system, the sets $F_v$ are usually subsets of the Euclidean
space $\R^n$ and the maps $\phi_e$ are contractions.

For more on the notion of (graph-directed) iterated function systems
and self-similarity structures see the papers~\cite{hutchinson,hata,kigami:harm,bandt:ss1,bandt:ss3} and the
book~\cite{falconer:tech}.

Graph-directed iterated function systems can be also viewed as
particular examples of \emph{topological graphs}
(see~\cite{deac:solen}).

\paragraph{The full one-sided shift $\xo$} is one of the most
important examples of the spaces with a standard self-similarity
structure. The respective iterated function system on $\xo$ has the structural graph with one vertex and
$|X|$ loops and is equal to the collection of the maps $\{\mathsf{T}_x\}_{x\in X}$, where
\[\mathsf{T}_x(w)=xw.\]

The image of $\xo$ under the map $\mathsf{T}_x$ is the cylindrical set
$x\xo$, so that $\xo=\cup_{x\in X} \mathsf{T}_x(\xo)$.

\paragraph{The shifts of finite type.} Let $\f\subseteq\xo$ be a
topological Markov chain. Let $\f_x=\f\cap x\xo$ be the set of all the
sequences $x_1x_2\ldots\in\f$ such that $x_1=x$. For every admissible word $xy\in
X^2$ we define $\mathsf{T}_{xy}:\f_y\to\f_x$ by the formula
\[\mathsf{T}_{xy}(w)=xw.\]

It is easy to check that we get in this way a graph-directed
iterated function system. The set of vertices of its structural graph is identified with $X$;
two vertices $x, y\in X$ are connected with an arrow starting in $x$ and ending in $y$ if
and only if the word $xy$ is admissible. 

The above two examples are in some sense universal, since they are
used to encode all the other self-similar sets.

\begin{defi}
  Let $\left(\{F_v\}_{v\in V}, \{\phi_e\}_{e\in E}\right)$ be a
  graph-directed iterated function system. Let $p\in F_v$ be an
  arbitrary point. We define a \emph{code} of the point $p$ as an
  infinite sequence $e_1e_2\ldots \in E^\omega$ such that for every
  $k\in\N$ we have
  \[p\in\phi_{e_1}\left(\phi_{e_2}\left(\ldots\phi_{e_{k-1}}\left(\phi_{e_k}\left(F_{\be(e_k)}\right)\right)\ldots\right)\right).
  \]
\end{defi}

In general, a point can have different codes, since the sets $F_i$ can
overlap. Also, different points can have the same codes.

It follows directly from the definition that for every code
$e_1e_2\ldots$ of a point we have $\be(e_k)=\en(e_{k+1})$ for all
$k\geq 1$, i.e., the sequence $\ldots e_2e_1$ is a (left-)infinite
path in the structural graph.

On the other hand, in the case of a topological iterated function
system, for every infinite path $\ldots e_2e_1$ in the associated
graph we have a decreasing sequence of compact sets
\[
\phi_{e_1}\left(F_{\be(e_1)}\right)\supseteq\phi_{e_1}\left(\phi_{e_2}\left(F_{\be(e_2)}\right)\right)\supseteq
\phi_{e_1}\left(\phi_{e_2}\left(\phi_{e_3}\left(F_{\be(e_3)}\right)\right)\right)\supseteq\ldots
\]
and every point in the intersection of these sets will have the
code $e_1e_2\ldots$; therefore the set of all possible codes is a
topological Markov chain in $E^\omega$ (which is defined by the set of
admissible words $\{xy: \be(x)=\en(y)\}$.

\paragraph{Markov partitions.} If the sets $F_v$ form a covering of a
set $F$, and if all the maps $\phi^{-1}$ (which are defined on subsets
of the sets $F_v$) are restrictions of a single map $f:F\to F$, then
the collection $\{\phi_e(F_{\be(e)})\}_{e\in E}$ is called the
\emph{Markov partition} of the dynamical system $(F, f)$
(though, usually more restrictions are imposed on the sets $F_v$).
In this case the described encoding of the
points of $F$ is the classical tool of the Symbolic Dynamics.
For the notion of a Markov partition see the survey~\cite{adler:markov}, the
book~\cite{kitchens} and the bibliography in them.

Note that the maps $T_x$ defining the iterated function system on
$\xo$ are the inverses of the shift $\shift$, so the sets $x\xo$ form
a Markov partition of the dynamical system $(\xo, \shift)$.

\begin{figure}[h]
  \begin{center}
    \includegraphics{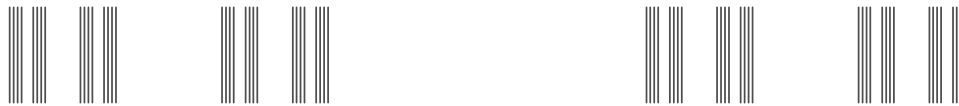}
  \end{center}
  \caption{The Cantor set}
  \label{fig:cantor}
\end{figure}

\paragraph{Other examples of self-similar sets.\\
  1. The Cantor middle-third set.}\label{ex:cantor} The full shift
$\xo$ for $|X|=2$ can be naturally interpreted as the classical Cantor
set. This is the set $C$ obtained from the segment $[0, 1]$ by
successive removing from the segments their middle thirds (see Figure~\ref{fig:cantor}). 

It follows that $C$ is the set of the numbers $x\in[0, 1]$ which
have only the digits $0$ and $2$ in their triadic expansion. In other
words
\[C=\left\{\sum_{n=1}^\infty\frac{d_n}{3^n} :d_1d_2d_3\ldots\in\{0,2\}^\omega\right\}.
\]

Moreover, the map $\Phi:d_1d_2d_3\ldots\mapsto\sum_{n=1}^\infty
d_n3^{-n}$ is a homeomorphism $\{0,2\}^\omega\to C$.

The standard iterated function system on $\{0,2\}^\omega$ can be
identified via the map $\Phi$ with the following natural iterated
function system on $C$: this system is the collection of two maps
$\phi_0:x\mapsto x/3$ and $\phi_1:x\mapsto x/3+2/3$.

It is easy to see now that the triadic numeration system gives exactly
the standard encoding of the Cantor set with respect to the described
iterated function system.

\paragraph{2. The segment $[0, 1]$.} Let $\phi_0(x)=x/2$ and $\phi_1(x)=x/2+1/2$
be two maps of the interval $[0, 1]$ to itself. In this way we obtain an iterated function system
on the segment $[0, 1]$.  The code of a point $x\in [0, 1]$ in this
case will be the sequence of digits in its dyadic expansion. Note that
here the dyadic expansion is not uniquely defined. For instance, the
point $1/2$ has two different codes: $100000\ldots$ and
$011111\ldots$.

\begin{figure}[ht]
  \begin{center}
    \psfrag{a1}{$A_1$} \psfrag{a2}{$A_2$} \psfrag{a3}{$A_3$}
    \psfrag{a4}{$A_4$} \psfrag{a}{(a)} \psfrag{b}{(b)} \psfrag{c}{(c)}
    \includegraphics{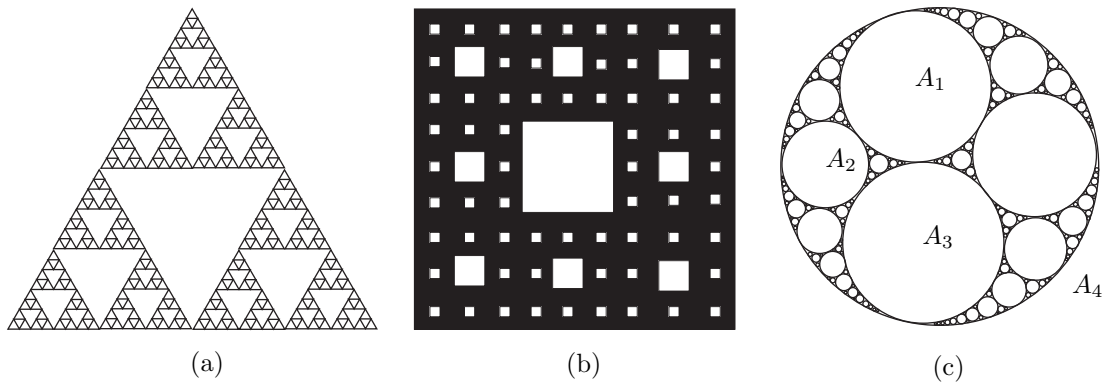}
  \end{center}
  \caption{Self-similar sets}
  \label{fig:apolo}
\end{figure}

\paragraph{3. The Sierpi\'nski gasket and Sierpi\'nski carpet.} The
Sierpi\'nski gasket (Figure~\ref{fig:apolo} (a)) is constructed from the
triangle with vertices $(0, 0)$, $(1, 0)$ and $(1/2,
\sqrt{3}/2)$ by successive deletion of the central triangles
(see~\cite{falcon:geom}).

It is a self-similar set with the iterated function system consisting
of three affine transformations $\phi_1\left(\vec x\right)=\vec x/2,
\phi_2\left(\vec x\right)=\vec x/2+(1/2, 0), \phi_3\left(\vec
  x\right)=\vec x/2+(1/4, \sqrt{3}/4)$.

The Sierpi\'nski carpet (Figure~\ref{fig:apolo} (b)) is constructed in a
similar way starting from the square $[0, 1]\times [0, 1]$ by deletion
of the central squares. It is also a self-similar set, for the
iterated function system
\begin{eqnarray*}
  \phi_1\left(\vec x\right)=\vec x/3,          &\quad & \phi_2\left(\vec x\right)=\vec x/3+(1/3, 0),\\
  \phi_3\left(\vec x\right)=\vec x/3+(2/3, 0), &\quad & \phi_4\left(\vec x\right)=\vec x/3+(0, 1/3),\\
  \phi_5\left(\vec x\right)=\vec x/3+(0, 2/3), &\quad & \phi_6\left(\vec x\right)=\vec x/3+(1/3, 2/3),\\
  \phi_7\left(\vec x\right)=\vec x/3+(2/3, 2/3), &\quad &
  \phi_6\left(\vec x\right)=\vec x/3+(2/3, 1/3).
\end{eqnarray*}

It is easy to see that the Sierpi\'nski gasket (resp.\ the Sierpi\'nski carpet) is uniquely determined
by the condition~\eqref{eq:defifs} for the corresponding set of affine transformations.
See the paper~\cite{hutchinson}, where properties of contracting iterated
function systems are investigated in a general setup.

\paragraph{4. The Apollonian net.}\label{ex:apollonian} The Apollonian
net (see~\cite{mandelbrot}) is a subset of the Riemann sphere, constructed in the following
way. Take four pairwise tangent circles $A_1, A_2, A_3, A_4$ (see Figure~\ref{fig:apolo} (c)).
We will get four curvilinear triangles with the vertices in the tangency points of the circles and
the sides coinciding with the respective arcs of the circles. Let us
denote the obtained triangles $T_1, T_2, T_3, T_4$.

Let us remove from the sphere the open disks bounded by the circles
$A_i$. Denote the obtained closed set by $\mathcal{P}_1$.  It is the
union of the triangles $T_i$, $i=1, \ldots, 4$. At the next stage we
inscribe into every triangle $T_i$ a maximal circle (this can be done
in a unique way), and remove the open disks bounded by these circles.
The obtained set $\mathcal{P}_2$ will be a union of $12$ triangles.
Then again, we inscribe a circle into each of these $12$ triangles and
remove from $\mathcal{P}_2$ the open disks bounded by these circles.
We denote the obtained union of $36$ triangles by $\mathcal{P}_3$. We
continue in a similar fashion. The \emph{Apollonian net} is the
intersection $\mathcal{P}=\cap_{i=1}^\infty\mathcal{P}_i$. See
Figure~\ref{fig:apolo} (c) for a picture of this fractal set.

Let us denote now by $\mathcal{T}_i$ the part of the Apollonian net
bounded by the triangle $T_i$. The sets $\mathcal{T}_i$ are called
\emph{Apollonian gaskets}. It is easy to prove that the Apollonian
gasket is homeomorphic to the Sierpinski gasket, but as metric spaces
they are very different.  For instance, the Hausdorff dimension of the
Sierpinski gasket is equal to $\log 3/\log 2$, while the exact value
of the Hausdorff dimension $h$ of the Apollonian gasket is not known.
At the moment very precise estimates of $h$ exists, and it can also be 
computed with arbitrary precision (see~\cite{boyda,boydb} 
and~\cite{MR2000d:37055}).  C.~T.~McMullen has shown that $h\approx 
1.305688$. 

The Apollonian gasket can be obtained by the same procedure of
removing inscribed disks but starting from a curvilinear triangle
formed by three tangent circles. It is a self-similar set with
self-similarity structure defined by the rational transformations
$f_i, i=1, 2, 3$, which map the triangle onto its three subtriangles.
For instance, if we take the triangle with the vertices $1, \exp(2\pi i/3), \exp(-2\pi i/3)$ then the
iterated function system will be $f_1(z)=\frac{(\sqrt{3}-1)z+1}{-z+\sqrt{3}+1}$, $f_2(z)=\exp(2\pi i/3)f_1(z)$
and $f_3(z)=\exp(-2\pi i/3)f_1(z)$ (see~\cite{mauldinurb}).

\subsection{Direct limits and self-replicating tilings}
\label{ss:dirlim}
Let $\left(\{F_v\}_{v\in V}, \{\phi_e\}_{e\in E}\right)$ be a
graph-directed iterated function system.

Let $e_1e_2\ldots$ be a path in the structural graph of the iterated
function system.  The \emph{leaf} defined by the path $e_1\ldots$ is
the direct limit of the sequence
\begin{equation}
  \label{eq:tillim}
  F_{\be(e_1)}\stackrel{\phi_{e_1}}{\longrightarrow}F_{\be(e_2)}\stackrel{\phi_{e_2}}{\longrightarrow}\cdots.
\end{equation}

We identify in a natural way the leaf defined by the path
$e_1e_2e_3\ldots$ with the leaf defined by the path $f_1f_2f_3\ldots$
if the paths are \emph{cofinal}, i.e., if $e_i=f_i$ for all $i$ big
enough. For a fixed path $e_1e_2e_3\ldots$, the \emph{central tile}
defined by the path is the image of $F_{\be(e_1)}$ in the direct
limit~\eqref{eq:tillim}. The central tiles defined by other paths,
cofinal with $e_1e_2e_3\ldots$ are called \emph{tiles} of the leaf.

It follows from Definition~\ref{defi:gdifs} that a leaf is the union of
its tiles. A \emph{tiled leaf} is a leaf together with its
decomposition into the union of its tiles.

\paragraph{Fibonacci tilings.} Let us construct an iterated function
system with two sets $F_1$ and $F_2$, where $F_1$ is the segment $[0,
1]$, and $F_2$ is the segment $[0, \tau]$, where $\tau=(1+\sqrt{5})/2$
is the \emph{golden mean}.

The functions are $\phi_{21}:F_2\to F_1$, $\phi_{12}:F_1\to F_2$ and
$\phi_{22}:F_2\to F_2$, where $\phi_{21}(x)=x/\tau$,
$\phi_{12}(x)=x/\tau$ and $\phi_{22}(x)=(x+1)/\tau$. So that
$\phi_{21}(F_2)=F_1$ and $\phi_{12}(F_1)=[0, 1/\tau]$,
$\phi_{22}(F_2)=[1/\tau, \tau]$, since $1+\tau=\tau^2$. In this way we
obtain the \emph{Fibonacci iterated function system}.

Its leaves are homeomorphic to the real lines tiled by segments of
length $1$ and $\tau$. For every such leaf we get a bi-infinite
sequence of symbols $1$ and $\tau$, denoting the lengths of the
respective tiles. The set of all such sequences is a closed subset of
the space of all bi-infinite sequences over the alphabet $\{1,\tau\}$
and is an example of a \emph{substitution dynamical system}. It can also
be viewed as a one-dimensional analogue of the Penrose tiling. See
for example the paper~\cite{bruijn:one}, where such tilings are
studied and are constructed using the \emph{projection method}. (The
Fibonacci tiling appears at the very end of the article.)

\paragraph{Penrose tilings.} The Penrose tiling is the most famous
example of an aperiodic tiling of the plane.  A tiling is
\emph{aperiodic} if it is not preserved by a non-trivial translation.

See the interesting description of the way in which it was invented by
R.~Penrose in~\cite{penrose}, the first article on these tilings by
M.~Gardner in~\cite{mgard} and its detailed analysis in the
book~\cite{tilings}.

\begin{figure}[ht]
  \psfrag{a}[Bc]{(a)} \psfrag{b}[Bc]{(b)}
  \begin{center}
    \includegraphics{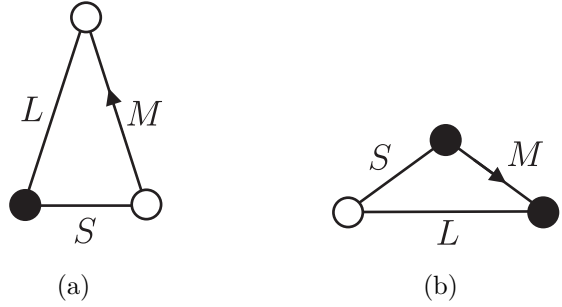}
  \end{center}
  \caption{Tiles and matching rules} \label{fig:match}
\end{figure}

The tiles of the Penrose tiling are the triangles shown on
Figure~\ref{fig:match}. The angles of the triangle~(a) are equal to
$\frac{\pi}{5}, \frac{2\pi}{5}, \frac{2\pi}{5}$ and the angles of the
triangle~(b) are equal to $\frac{3\pi}{5}, \frac{\pi}{5},
\frac{\pi}{5}$. These triangles are cut from a regular pentagon by the
diagonals issued from a common vertex, so that the ratio between the
lengths of the shorter and the longer sides of the triangles is equal
to $\tau=\frac{1+\sqrt{5}}{2}$.

A tiling of the plane by the triangles is a Penrose tiling if it
satisfies the \emph{matching rules}, which require that the common
vertices of two tiles are marked by the same color
(black or white) and that the arrows on the adjacent sides of the
triangles point at the same direction. The colors of the vertices
and the arrows on the sides are shown on Figure~\ref{fig:match}.

It follows from the matching rules that if we mark the sides of the
tiles by the letters $S, L$ and $M$ in the way it is done on
Figure~\ref{fig:match}, then in every Penrose tiling the common sides of
any two adjacent tiles will be marked by the same letter.

\begin{figure}[ht]
  \begin{center}
    \includegraphics{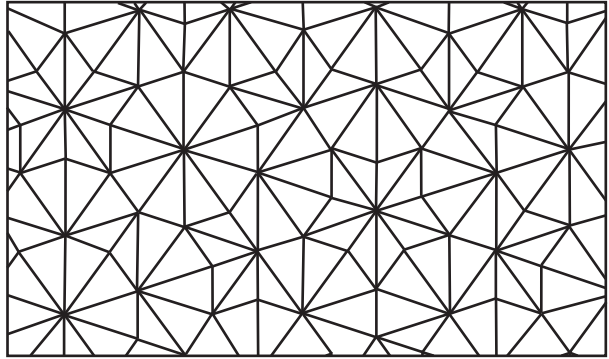}
  \end{center}
  \caption{Penrose tiling} \label{fig:tiling}
\end{figure}

An example of a patch of a Penrose tiling is shown on
Figure~\ref{fig:tiling}.

\begin{figure}[ht]
  \psfrag{1}{$a$} \psfrag{2}{$b$} \psfrag{3}{$c$} \psfrag{4}{$a$}
  \psfrag{5}{$c$}
  \begin{center}
    \includegraphics{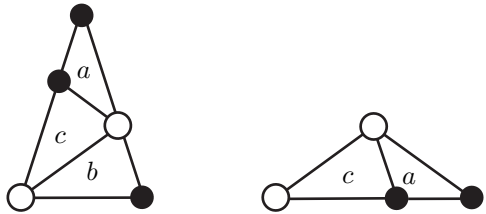}
  \end{center}
  \caption{Inflation} \label{fig:inflation}
\end{figure}

It follows from the matching rules that in any Penrose tiling one can
group the tiles into blocks of two or three tiles, shown on
Figure~\ref{fig:inflation} (their mirror images are also allowed).

The grouping into blocks is unique~\cite{tilings} and it is easy to
see that the blocks are similar to the original tiles, with similarity
coefficient $\tau$. Moreover, the grouping agrees with the matching
rules, so that the blocks also form a Penrose tiling. This tiling is
called the \emph{inflation} of the original one.

It follows from existence and uniqueness of the inflation that every Penrose tiling
is aperiodic. It also follows from inflation that any two Penrose
tilings are locally isomorphic. See the book~\cite{tilings} for the proofs
of these and other properties of the Penrose tilings.

\begin{figure}[ht]
  \begin{center}
    \includegraphics{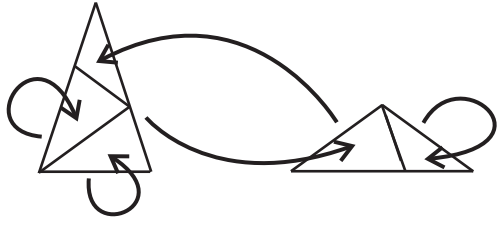}
  \end{center}
  \caption{Penrose iterated function system}
  \label{fig:infl}
\end{figure}

The inflation rule defines obviously a graph-directed iterated
function system with two sets $F_1$ and $F_2$, equal to the triangles
and five similarities, as shown on Figure~\ref{fig:infl}. This
iterated function system is called the \emph{Penrose iterated function
  system}.

The symmetry group of a Penrose tiling is either trivial, or is of
order two (then the only non-trivial symmetry of the tiling is a
symmetry with respect to a line), or is the dihedral group $D_5$ (then
the only symmetries of the tiling are rotations around a point on the
angles $k\cdot\frac{2\pi}{5}$ and symmetries with respect to five
lines passing through the rotation point). More detailed analysis of
the symmetric situations will be made later.

The tiled leaves of the direct limit of the Penrose iterated function
system will be either a Penrose tiling of the plane with trivial
symmetry group, or the fundamental domain (a half-plane, or an
$36$-degree angle) of a Penrose tiling with a non-trivial symmetry
group.

\subsection{Self-similar groups}
\label{ss:ssgr}

The aim of this section is to define \emph{self-similar group actions},
which are analogous to self-similar spaces. The analogy will become more clear in
Section~\ref{s:cstar}, where both notions will be treated within a common
algebraic approach.

\begin{defi}
  \label{defi:ssact}
  Suppose a collection $\{\phi_e\}_{e\in E}$ of functions $F\to F$
  defines an iterated function system on the space $F$.  An action of
  a group $G$ on $F$ is then said to be \emph{self-similar} if for
  every $e\in E$ and for every element $g\in G$ there exists $f\in E$
  and $h\in G$ such that
  \[
\left(\phi_e\left(p\right)\right)^g=\phi_f\left(p^h\right)\quad\text{
  for all }p\in F.
  \]
\end{defi}

The most favorable case is when the images of the set $F$ under the
maps $\phi_e$ do not intersect. 
We also assume that different points of $F$ have different codes.  The space $F$ is then
homeomorphic to $\xo$ for the alphabet $X=E$, and the functions
$\phi_e$ become equal to the maps $\mathsf{T}_x:w\mapsto xw$ on $\xo$.

We therefore reformulate Definition~\ref{defi:ssact} for the space $F=\xo$
in the following way:
\begin{reformulate}{defi}{\ref{defi:ssact}'}\label{defi:ssact2}
  An action of a group $G$ on the space $\xo$ is \emph{self-similar}
  if for every $g\in G$ and every $x\in X$ there exist $h\in G$ and
  $y\in X$ such that
  \begin{equation}
    (xw)^g=y(w^h)
    \label{eq:ssim}
  \end{equation}
  for every $w\in\xo$.
\end{reformulate}

Applying Equation~\eqref{eq:ssim} several times we see that for every
finite word $v\in\xs$ and every $g\in G$ there exist $h\in G$ and a
word $u\in\xs$ such that $|u|=|v|$ and
\[(vw)^g=u(w^h)\]
for all $w\in\xo$. Hence we get a naturally-defined associated
action of $G$ on $\xs$, for which the word $u$ is the image of $v$
under the action of $g$. If $G$ acts faithfully then the group element
$h$ is defined uniquely.  It is called then the \emph{restriction} of
$g$ in the word $v$ and is written $h=g|_v$.

It follows directly from the definitions that for any $v, v_1,
v_2\in\xs$ and $g, g_1, g_2\in G$ we have
\begin{equation}\label{eq:restr}
g|_{v_1v_2}=\left(g|_{v_1}\right)|_{v_2}\qquad
(g_1g_2)|_v=\left(g_1|_v\right)\left(g_2|_{v^{g_1}}\right).
\end{equation}

\begin{defi}
  A self-similar action is \emph{level-transitive} for every $n\ge 0$
  the associated action on $\xs$ is transitive on the set $X^n$ of the
  words of length $n$.
\end{defi}

All self-similar actions in this paper are assumed to be
level-transitive.

\subsection{Examples of self-similar group actions}
\label{ss:exss}
\paragraph{The adding machine.} Let $a$ be the transformation of the
space $\{0,1\}^\omega$ defined by the following recursion formul\ae:
\begin{eqnarray*}
  (0w)^a & = & 1w\\
  (1w)^a & = & 0w^a,
\end{eqnarray*}
where $w$ is an arbitrary infinite word over the alphabet $\{0,1\}$.

These formul\ae\ can be interpreted naturally as the rule of adding $1$ to
a dyadic integer. More precisely, identifying $\{0,1\}^\omega$ with
the dyadic integers $\Z_2$ via $\Phi:x_1x_2\dots\mapsto
\sum_{i\ge1}x_i2^{i-1}$, we have
\[\Phi(w^a)=\Phi(w)+1.\]
The term ``adding machine'' originates from
this interpretation of the transformation $a$.

The transformation $a$ generates an infinite cyclic group of
transformations of the space $\{0,1\}^\omega$.  Thus we get an action
of the group $\Z$, which will be also called \emph{adding machine
  action}. It is easy to see that the adding machine action is
self-similar.

\paragraph{The dihedral group.} Let $a$ and $b$ be the transformations
of the space $\xo=\{0,1\}^\omega$, defined by the rules
\[
\begin{array}{lcl}
  (0w)^a=1w &\quad & (0w)^b=0w^a \\
  (1w)^a=0w &\quad & (1w)^b=1w^b,
\end{array}
\]
were  $w\in\xo$ is arbitrary.

The group generated by the transformations $a$ and $b$ is
isomorphic to the infinite dihedral group $\mathbb{D}_\infty$, and
thus we get a self-similar action of this group on $\xo$.

\paragraph{The Grigorchuk group.} The Grigorchuk group
(see~\cite{grigorchuk:80_en} for the original definition) is the
transformation group of the space $\{0,1\}^\omega$ generated by the
transformations $a, b, c, d$, which are defined by the rules:
\[
\begin{array}{lll}
  (0w)^a=1w & \quad & (1w)^a=0w\\
  (0w)^b=0w^a &\quad & (1w)^b=1w^c\\
  (0w)^c=0w^a & \quad & (1w)^c=1w^d\\
  (0w)^d=0w &\quad & (1w)^d=1w^b.
\end{array}
\]

The Grigorchuk group is an infinite finitely generated torsion group. This fact relates it to the
General  Burnside Problem. It is also
the first example of a group of intermediate
growth~\cite{grigorchuk:milnor_en} (about growth see Section~\ref{s:growth}).  For
more on the Grigorchuk group see~\cite{grigorchuk:branch} and the last chapter of~\cite{harpe}.

\paragraph{Other examples.} The ``\emph{lamplighter group}''
$(\Z/2\Z)^\Z\rtimes\Z$, where $\Z$ acts on $(\Z/2\Z)^\Z$ by the shift
is isomorphic to the transformation group of the space
$\{0,1\}^\omega$ generated by the transformations
\begin{eqnarray*}
(0w)^a=1w^b & \quad & (0w)^b=0w^b\\
(1w)^a=0w^a & \quad &(1w)^b=1w^a.
\end{eqnarray*}
To see this, identify $\{0,1\}^\omega$ with $(\Z/2\Z)[[t]]$ via
$\Phi: x_1x_2\dots\mapsto\sum_{i\ge1}x_it^{i-1}$. We then have (see~\cite{grineksu_en})
\[\Phi(w^{b^{-1}a})=\Phi(w)+1,\quad\Phi(w^b)=(1+t)\Phi(w).\]
This action was used in~\cite{gr_zu:lamp} to compute the spectrum
of the lamplighter group.

The self-similar representations of \emph{abelian} and \emph{affine groups}
were studied in the papers~\cite{sibr:gln} and~\cite{neksid}.

\subsection{Self-similar inverse semigroups}
A \emph{partial permutation} of a set $M$ is a bijection
$\pi:\dom\pi\to\ran\pi$ between a subset $\dom\pi\subseteq M$ (the
\emph{domain} of the permutation) and a subset $\ran\pi$ (the
\emph{range} of the permutation).  We admit also the \emph{empty permutation}
$0$ with empty domain and range.

Product of two partial permutations $\pi_1, \pi_2$ is the partial permutation
$\pi_1\pi_2$ with the domain 
$$\dom\pi_1\pi_2=\left\{x\in\dom\pi_1: x^{\pi_1}\in\dom\pi_2\right\}$$
defined on it by the condition $x^{\pi_1\pi_2}=\left(x^{\pi_1}\right)^{\pi_2}$.
In particular, $\pi\cdot 0=0\cdot \pi=0$ for every partial permutation.

Every partial permutation has its \emph{inverse} $\pi^{-1}:\ran\pi\to\dom\pi$. Note that the
product $\pi^{-1}\pi$ is not the identity but an \emph{idempotent}. It is
the identical permutation with the domain coinciding with the range of
the permutation $\pi$.

An \emph{inverse semigroup} acting on a set $M$ is a set of partial
permutations of the set $M$ which is closed respectively to
the composition and the inversion.

An example of an inverse semigroup is the \emph{symmetric inverse
  semigroup} $IS(M)$ which is the semigroup of all partial
permutations of the set $M$.

If $\pi$ and $\rho$ are two partial permutations such that
$\dom\pi\cap\dom\rho=\emptyset$ and $\ran\pi\cap\ran\rho=\emptyset$
then the \emph{sum} of the permutations $\pi$ and $\rho$ is
permutation $\pi+\rho$ with domain $\dom\pi\cup\dom\rho$, defined by
the condition
\[
x^{\pi+\rho}=\left\{
  \begin{array}{lcl}
    x^\pi   & \quad & \text{if $x\in\dom\pi$}\\
    x^\rho & \quad & \text{if $x\in\dom\rho$.}
  \end{array}
\right.
\]

Just as in the case of group actions, we consider only actions of
inverse semigroups on topological Markov chains (with standard
iterated function systems on them).

\begin{defi}
  \label{defi:ssinv}
  Let $\f\subseteq\xo$ be a topological Markov chain.
  
  An inverse semigroup $G$ acting on $\f$ is \emph{self-similar} if
  for every $g\in G$ and $x\in X$ there exist $y_1, y_2, \ldots,
  y_k\in X$ and $h_1, h_2, \ldots, h_k\in G$ such that the sets $\dom
  h_i$ are disjoint, $\cup_{i=1}^k y_i\dom h_i=x\xo\cap\dom g$, and
  for every $w\in\xo$ we have
  \[(xw)^g=y_i (w)^{h_i},\]
  where $i$ is such that $w\in\dom h_i$.
\end{defi}
It follows from the definition that the sets $y_i\dom h_i$ are also
disjoint.

Let us denote by $\mathsf{T}_x$ the partial permutation of the space
$\f$ which acts by the rule $\mathsf{T}_x(w)=xw$. The domain of
$\mathsf{T}_x$ is the set of those elements $w\in\f$, for which $xw$
belongs to $\f$. The range of $\mathsf{T}_x$ is then the set
$x\xo\cap\f$.

The previous definition can be reformulated in the following way:

\begin{reformulate}{defi}{\ref{defi:ssinv}'}
  An inverse semigroup $G$ acting on the space $\xo$ is
  \emph{self-similar} if for every $g\in G$ and $x\in X$ there exist
  $y_1, y_2, \ldots, y_k\in X$ and $h_1, h_2, \ldots, h_k\in G$ such
  that
  \begin{equation}
    \label{eq:ssinv}
    \mathsf{T}_xg=h_1\mathsf{T}_{y_1}+h_2\mathsf{T}_{y_2}+\cdots+h_k\mathsf{T}_{y_k}.
  \end{equation}
\end{reformulate}

\subsection{Examples of self-similar inverse semigroups.}
\paragraph{Fibonacci transformations.} Let $\f$ be the set of all
infinite words over the alphabet $X=\{0,1\}$ which do not contain a
subword $11$. The space $\f$ is a shift of finite type called the
\emph{Fibonacci shift} (see~\cite{marcus,kitchens}).

Let us define two partial homeomorphisms $a, b$ of the space $\f$
having domains $0\xo\cap\f$ and $1\xo\cap\f$
respectively, by the inductive formula
\[
\left\{
  \begin{array}{lcl}
    (00w)^a & = & 10w\\
    (01w)^a & = & 0(1w)^b\\[1ex]
    (1w)^b  & = & 0(w^a).
  \end{array}
\right.
\]

It follows from the inductive definition that the transformation $a$
acts by the rule
\[a:(01)^n 00w\mapsto (00)^n10w,\]
where $n\geq 0$ and $w\in\f$ is an arbitrary word. Additionally, (in
some sense, for $n=\infty$) we have $(010101\ldots)^a=000000\ldots$.

The transformation $b$ acts by the rule
\[b:1(01)^n00w\mapsto 0(00)^n10w,\]
where $n\geq 0$ and $w\in\f$. In the limit case we have
$(101010\ldots)^b=000000\ldots$.

It follows from the formul\ae\ above that $\ran a$ is equal to the set
of all the words from $\f$ with even or infinite number of leading
zeros and $\ran b$ is the set of the words with odd or infinite number
of leading zeros. These two sets intersect in the point $000\ldots$.
Thus the sum $a+b$ (which is defined as the map $\f\to\f$ equal to $a$
on $0\xo\cap\f$ and equal to $b$ on $1\xo\cap\f$) is not invertible.

The transformations $a$ and $b$ generate a self-similar inverse semigroup.

The partial homeomorphisms $a$ and $b$ are associated to a special
numeration system (called the \emph{Fibonacci system}) on the
integers, in a similar way like the binary adding machine is
associated to the usual dyadic numeration system.

Let $u_1=1, u_2=2, u_3=3, u_4=5, u_5=8, \ldots$ be the Fibonacci
sequence defined by the recursion $u_n=u_{n-1}+u_{n-2}$. Then (see for
instance~\cite{knuth}) every positive integer $m$ can be uniquely
presented as a finite sum
\[a_1u_1+a_2u_2+\cdots+a_nu_n,\]
where $a_i\in\{0,1\}$ and $a_i=1$ for $i<n$ implies that $a_{i+1}=0$.

If we put into correspondence to the word $a_1a_2\ldots\in\f$ a formal
infinite sum
\[a_1u_1+a_2u_2+\cdots,\]
then the map $a+b$ can be interpreted as adding of $1$ to this formal
expression. The recurrent definitions of $a$ and $b$ come from the
``carrying'' rules for addition in the Fibonacci system, which in turn
come from the relations $u_1+1=u_2$ and $u_i+u_{i+1}=u_{i+2}$.

A relation between the Fibonacci transformations and the Figonacci tiling, introduced before, will
be clarified in Subsection~\ref{ss:lsis}.

\paragraph{A group and an inverse semigroup associated to the Penrose
  tilings.} Let us label the tiles of the Penrose tiling by letters
$a, b, c$ according to their participation in the inflation process
as it is shown in Figure~\ref{fig:inflation}. Note that every
obtuse-angled triangle is labeled by $a$ and the acute-angled
triangles are labeled either by $b$ or by $c$ depending on their role
in the inflation.

Let $u$ be a tile of a Penrose tiling of the plane. Let $x_1\in\{a, b,
c\}$ be its label. After inflation the tile $u$ becomes a part of a
tile $\shift(u)$ in the new inflated Penrose tiling.  Let $x_2$ be the
label of the tile $\shift(u)$. In general, let $x_n$ be the label of
the tile to which $u$ belongs in the Penrose tiling obtained from the
original one by $n-1$ successive inflations. The sequence
$x_1x_2\ldots\in\{a, b, c\}^\omega$ is called the \emph{code} of the
tile $u$.

A sequence $x_1x_2\ldots\in\{a, b, c\}^\omega$ is the code of a tile
in a Penrose tiling if and only if it does not contain a subsequence
$x_ix_{i+1}$ such that $x_i=b, x_{i+1}=a$.  Let
$\mathcal{P}\subset\{a, b, c\}^\omega$ be the set of such sequences.

The shift $\shift:\mathcal{P}\to\mathcal{P}$ encodes the
inflation, more precisely, it maps the code of a tile $u$ to the
code of the tile $\shift(u)$, to which $u$ belongs after the
inflation.

Two sequences $u, v\in\mathcal{P}$ are codes of tiles belonging to a
common tiling if and only if they are \emph{cofinal}, i.e., if for
some $n\in\N$ the sequences $\shift^n(u)$ and $\shift^n(v)$ are equal.
Since every cofinality class is countable, there exist uncountably
many non-isomorphic Penrose tilings.

Two tiles $u$ and $v$ have the same codes if and only if they belong
to a common tiling and there exists a symmetry of the tiling which
carries $u$ to $v$.

Therefore, we shall identify a tile with its code, keeping in mind the
non-uniqueness described in the previous paragraph. Then $\mathcal{P}$
becomes a union of the sets of the tiles of all Penrose tilings and
the set of cofinality classes on $\mathcal{P}$ becomes the set of all
Penrose tilings.

Let $L, M, S:\mathcal{P}\to\mathcal{P}$ be the maps which carry a tile
to the neighbor, adjacent to the side labeled by the letter $L, M$ or
$S$, respectively. We follow the labeling rules shown on
Figure~\ref{fig:match}.

It follows from the matching rules that the maps $L, M$ and $S$ are
involutions, since common sides of tiles are labeled by the same
letters. Therefore, the maps are invertible, and they give us an
action of the free product $F=\Z/2\Z*\Z/2\Z*\Z/2\Z$ of three cyclic
groups of order two on the set $\mathcal{P}$. Then the orbits of the
action of $F$ on $\mathcal{P}$ are in one-to-one correspondence with
the Penrose tilings.

The fact that the action of $F$ on $\mathcal{P}$ is continuous can be
easily deduced from the properties of the Penrose tilings, but this
also directly follows from the formul\ae\ defining $L, M$ and $S$,
written below:

\begin{theorem}
  \label{th:penr}
  The transformations $L, M, S$ act on $\mathcal{P}$ according to
  the following rules:
  \[
  \begin{array}{rclcrcl}
    (aw)^S   & = & cw           &\quad & (aaw)^L  & = & b\cdot (aw)^S\\
    (bw)^S   & = & b\cdot (w)^M  &\quad & (abw)^L  & = & a\cdot (bw)^M\\
    (cw)^S   & = & aw           &\quad & (acw)^L  & = & a\cdot (cw)^M\\
    &   &              &\quad & (bbw)^L  & = & b\cdot (bw)^S\\
    (aw)^M   & = & a\cdot (w)^L  &\quad & (bcw)^L  & = & a\cdot (cw)^S\\
    (bw)^M   & = & cw           &\quad & (cw)^L   & = & c\cdot (w)^L\\
    (caw)^M  & = & c\cdot (aw)^M\\
    (cbw)^M  & = & bbw\\
    (ccw)^M  & = & bcw
  \end{array}
  \]
\end{theorem}

Let us now define partial homeomorphisms $S_x, M_x, L_x$, respectively,
equal to the restriction of $S, M$ or $L$ onto the set
$x\mathcal{P}\cap\mathcal{P}$ when $x$ traverses $\{a,b,c\}$. Let $G$
be the inverse semigroup generated by the obtained partial
homeomorphisms. Then the homeomorphisms $S, M$ and $L$ are equal to
the sums $S_a+S_b+S_c, M_a+M_b+M_c$ and $L_a+L_b+L_c$ respectively. By
$\mathsf{T}_x, x\in\{a, b, c\}$ we denote, as usual, the partially
defined transformations $\mathsf{T}_x:w\mapsto xw$. Let
$P_x=\mathsf{T}_x^{-1}\mathsf{T}_x$ be the projection onto the
cylindrical set $x\xo\cap\mathcal{P}$. We have $1=P_a+P_b+P_c$.

In this notation the formul\ae\ from Theorem~\ref{th:penr} can be
rewritten as:
\[
\begin{array}{rclrcl}
  \mathsf{T}_aS_a  & = & \mathsf{T}_c & \mathsf{T}_a M_a & = &  L\mathsf{T}_a\\
  \mathsf{T}_b S_b & = & (M_b+M_c)\mathsf{T}_b & \mathsf{T}_bM_b & = & \mathsf{T}_c\\
  \mathsf{T}_c S_c & = & \mathsf{T}_a  & \mathsf{T}_c M_c & = & M_a
  \mathsf{T}_c+(P_b+P_c)\mathsf{T}_b\\[1ex]
  \mathsf{T}_aL_a & = & S_a \mathsf{T}_b+(M_b+M_c)\mathsf{T}_a\\
  \mathsf{T}_bL_b & = & S_b \mathsf{T}_b+S_c \mathsf{T}_a\\
  \mathsf{T}_c L_c & = & L\mathsf{T}_c.
\end{array}
\]
These are the only non-zero products of a generator with a
transformation of the form $\mathsf{T}_x$.

\begin{corollary}
  The inverse semigroup $G=\langle L_x, M_x, S_x, P_x: x=a, b,
  c\rangle$ is self-similar.
\end{corollary}

Suppose that we have a Penrose tiling with non-trivial symmetry.  Let
$u$ be a tile adjacent to the symmetry axes. Then either $u^M=u$, or
$u^L=u$, or $u^S=u$. A careful analysis of the formul\ae\ shows that in
the first case $u$ has to be of the form
$v_0c^{n_1}v_1c^{n_2}v_2c^{n_3}v_3\ldots$, where $v_0$ is equal either
to $ca$ or to $a$, the $v_i$ for $i\ge 1$ belong to the set $\{aca,
bbca\}$, and $n_i\geq 1$ arbitrary integers. In the second case $u$ is
of the form $c^{n_1}v_1c^{n_2}v_2c^{n_3}v_3\ldots$ and in the third
case it has the form $bbcac^{n_1}v_1c^{n_2}v_2c^{n_3}v_3\ldots$, with
the same conditions on the words $v_i$.

Thus the set of the tiles in the Penrose tilings with a non-trivial
symmetry is uncountable, although it is nowhere dense.

All these tilings have a symmetry group of order two (i.e., they have
only one non-trivial symmetry), except for a unique tiling with
symmetry group $D_5$. It is the tiling containing the tile with the
code $cacaca\ldots$.  This code is fixed by $M$ and by $L$.  See the
monograph~\cite{tilings} for more details about this exceptional
tiling.

\subsection{Semigroups and groups of self-similarities}
\label{ss:gofss}

The self-similarities also generate interesting (semi)group actions.

If $\left\langle\{F_v\}_{v\in V}, \{\phi_e\}_{e\in E}\right\rangle$ is
a self-similarity structure on a set $F$, then the functions $\phi_e$
generate a semigroup of \emph{self-similarities}, or an inverse
semigroup of self-similarities, if we take also the inverses of the
maps $\phi_e$.

\paragraph{Example.} Recall that the middle-thirds Cantor set $C$ from
page~\pageref{ex:cantor} is self-similar with respect to the
self-similarity structure $\{\phi_0(x)=\frac{1}{3}x,
\phi_1(x)=\frac{1}{3}x+\frac{2}{3}\}$. The semigroup $G$ generated by
the transformations $\phi_0$ and $\phi_1$ acts on the Cantor set $C$.
It is easy to see that the semigroup $G$ is free.

More generally, if we have an iterated function system $\left\langle
  F, \{\phi_i\}\right\rangle$ such that the sets $\phi_i(F)\subset F$
have disjoint subsets, then the semigroup generated by the maps
$\phi_i$ will be free.
\medskip

As another class of examples, consider the actions which locally
coincide with the self-similarities. More precisely: let $H_0$ be the
inverse semigroup of self-similarities. Let $H$ be the semigroup
consisting of all sums of the elements of the semigroup $H_0$. We say
that a group $G$ acts on $F$ by \emph{piecewise self-similarities}, if
it is a subgroup of the pseudogroup $H$.

\paragraph{The Thompson groups.} The maximal group acting by piecewise
self-similarities on the full shift space $\xo$ is called the
\emph{Higman-Thompson group}.

The elements of the Higman-Thompson group $V_d$ are the
transformations of the set $\xo$ (for $|X|=d$) defined by
\emph{tables} of the form
\[
\left(
\begin{array}{cccc}
  v_1 & v_2 & \ldots & v_m\\
  u_1 & u_2 & \ldots & u_m
\end{array}
\right),
\]
where $m\in\N$ and $v_i, u_i\in \xs$ are finite words such that for
every infinite word $w\in\xo$ exactly one word $v_i$ and exactly one
word $u_j$ is a beginning of $w$. The transformation $\tau$ defined by
the above table acts on the infinite words by the rule
\[(v_iw)^\tau=u_iw,\]
where $w$ is an arbitrary infinite word over the alphabet $X$.

Denote by $V_d'=\left[V_d, V_d\right]$ the commutant of the group
$V_d$.  The following theorem holds
(see~\cite{hgthomp,intro_tomp,thompson}).

\begin{theorem}
  \label{th:thomps}
  For every $d>1$ the Higman-Thompson group $V_d$ is finitely
  presented. For even $d$ the groups $V_d$ and $V_d'$ coincide.  For
  odd $d$ the group $V_d'$ is of index $2$ in $V_d$.
  
  The group $V_d'$ is the only non-trivial normal subgroup of $V_d$
  and is simple.
\end{theorem}

The group $V_2$ and its analogs where constructed by R.~Thompson
in~1965.  These group were used in~\cite{mcthomp} to construct an
example of a finitely presented group with unsolvable word problem and
in~\cite{thompson} for embeddings of groups into finitely presented
simple groups.

The maximal group of piecewise self-similarities of the segment $[0,
1]$ (with respect to the iterated function system $\{x\mapsto x/2,
x\mapsto x/2+1/2\}$ is also called the Thompson group, and is denoted
$F$. It is the group of the piecewise linear continuous transformations
of the segment $[0,1]$, which are differentiable in every point,
except for a finite number of points of the form $m/2^n, (n\geq 1,
0\le m\le 2^n)$ and such that the derivative in all the points where
it exists is an integral power of 2.

The group $F$ is a finitely presented group. Its commutant is simple
and $F/F'$ is isomorphic to the group $\Z^2$.  The group $F$ does not
contain a free subgroup and has no group laws, i.e., no identities of the
form $w(x_1, x_2,\ldots, x_n)=1$, which are true for any substitution of the group elements
into the variables $x_1, x_2,\ldots, x_n$. Here $w(x_1, x_2, \ldots, x_n)$ is a nonempty freely
reduced word in the alphabet $x_1, x_1^{-1}, x_2, x_2^{-1}, \ldots, x_n, x_n^{-1}$.

See the survey~\cite{intro_tomp} for more properties of the Thompson
groups.

\paragraph{The group associated with the Apollonian net.} We use here
the notations from the definition of the Apollonian net in
Subsection~\ref{ss:selfsims}, page~\pageref{ex:apollonian}.

Let $C_i$, $i=1, \ldots, 4$ denote the circle circumscribed around the
triangle $T_i$. The circle $C_i$ is orthogonal to the circles forming
the triangle $T_i$ (in other words, it is orthogonal to the sides of
the triangle $T_i$).

Let $\gamma_i$ be the inversion of the Riemann sphere with respect to
the circle $C_i$. It is easy to deduce from the definition of the
Apollonian net that it is invariant under each of the transformations
$\gamma_i$, so it is invariant under the action of the group $\Gamma$
generated by the set $\{\gamma_1, \gamma_2, \gamma_3, \gamma_4\}$.

It is also easy to see that $\gamma_i$ maps the gasket $\mathcal{T}_i$
onto the union of the other three gaskets $\mathcal{T}_j$, $j\neq i$.
On the other hand, the union of the sets $\mathcal{T}_j$ for $j\neq i$
is mapped by $\gamma_i$ onto $\mathcal{T}_i$ and each of the sets
$\mathcal{T}_j$ is mapped onto one of the three smaller triangles
forming $\mathcal{T}_i$.  Thus the inversions $\gamma_i$ restricted
onto the respective triangles $\mathcal{T}_i$ form a self-similarity
structure on the Apollonian gasket $\mathcal{P}$.

Therefore, we can define the encoding of the points of the net using
this self-similarity structure. Let $X=\{1,2,3,4\}$ be our
alphabet. Recall that the encoding is defined in such a way that if $x$
is a point of the Apollonian net $\mathcal{P}$ then the first letter
of the code is the index $i_1$ such that $x\in\mathcal{T}_{i_1}$. The
next letters of the code are defined by the condition that if the word
$i_1v\in\xs$ is the code of the point $x$ then the word $v$ must be
the code of the point $\gamma_{i_1}(x)$. Then the set of the points
whose code starts with a given word of length $n$ is one of the
$4\cdot 3^{n-1}$ triangles constructed on the $n$th stage of the
definition of the Apollonian net. It then follows that there exists
not more than one point of the net having a given code $w\in\xo$. It
also follows from the construction that a word $w\in\xo$ is a code of
a point of $\mathcal{P}$ if and only if it has no equal consecutive
letters.

We immediately conclude that the action of the generators $\gamma_i$
of the group $\Gamma$ is the following:
\begin{eqnarray*}
  (iw)^{\gamma_i} &=& w;\\
  (w)^{\gamma_i} &=& iw\quad\text{if the first letter of $w$ is not $i$.}
\end{eqnarray*}

It then easily follows that the group $\Gamma$ is isomorphic to the
free product of four groups of order $2$; the set of finite-length
codes is a principal $\Gamma$-space through the natural identification
$i_1i_2\dots i_n\leftrightarrow
\gamma_{i_1}\gamma_{i_2}\dots\gamma_{i_n}$.

The above formul\ae\ also imply that the action of the group $\Gamma$ is
minimal, i.e.\ that the closure of every $\Gamma$ orbit is the whole
Apollonian net.

\section{Automata groups and actions on rooted trees}
\label{s:automata}
\subsection{Automata}
\label{ss:automata}
Equation~\eqref{eq:ssim} in Subsection~\ref{ss:ssgr} can be interpreted as a work of the machine,
which being in a state $g$ and receiving as input a letter $x$, goes
into state $h$ and outputs the letter $y$. Such a machine is
formalized by the following definition:

\begin{defi}
  Let $X$ be an alphabet. An \emph{automaton} (\emph{transducer}) $\autom$ over an alphabet $X$ is a triple
$\langle Q, \lambda, \pi\rangle$, where
  \begin{enumerate}
  \item $Q$ is a set (\emph{the set of the internal states} of the
    automaton $\autom$);
  \item $\lambda:Q\times X\to X$ is a map, called \emph{the output
      function} of the automaton;
  \item $\pi:Q\times X\to Q$ is a map, called \emph{the transition
      function} of the automaton.
  \end{enumerate}
  
  An automaton is \emph{finite} if the set $Q$ is finite.
  
  A subset $S\subseteq Q$ is a \emph{subautomaton} of $\autom$ if for
  all $s\in S, x\in X$ the state $\pi(s, x)$ belongs to $S$.
\end{defi}

For a general theory of automata see~\cite{eil}.

\begin{figure}[ht]
  \psfrag{1,0}{$(1, 0)$} \psfrag{0,1}{$(0,1)$} \psfrag{0,0}{$(0,
    0)$} \psfrag{1,1}{$(1, 1)$}
  \begin{center}
    \includegraphics{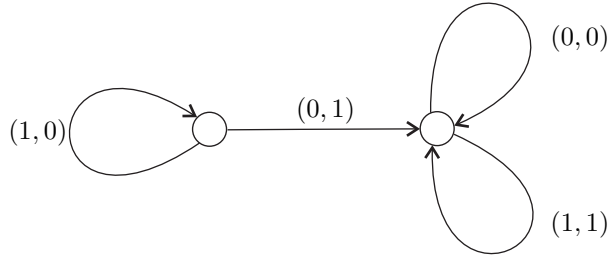}
  \end{center}
  \caption{A Moore diagram} \label{fig:admach}
\end{figure}

It is convenient to define automata by their diagram (their
\emph{Moore diagram}). Such a diagram is a directed labeled graph with
vertices identified with the states of the automaton. For every state
$q$ and letter $x\in X$ the diagram has an arrow from $q$ to $\pi(q,
x)$ labeled by the pair $(x, \lambda(q, x))$.  An example of Moore
diagram is given on Figure~\ref{fig:admach}.

We interpret the automaton as a machine, which being in a state $q$
and reading on the input tape a letter $x$, goes to the state $\pi(q,
x)$, types on the output tape the letter $\lambda(q ,x)$, then moves
both tapes to the next position and proceeds further.

In this way, we get a natural action of the automaton on the words.
Namely, we extend the functions $\lambda$ and $\pi$ to $Q\times\xs$
in the following way:
\begin{eqnarray}\label{recur1}
 \pi(q,\emp)=q    & \quad & \pi(q, xv)= \pi(\pi(q, x), v),\\
 \label{recur} \lambda(q, \emp)=\emp & \quad & \lambda(q,
 xv)=\lambda(q, x)\lambda(\pi(q, x), v).
\end{eqnarray}

Equations~\eqref{recur} also define uniquely a map
$\lambda:Q\times\xo\to\xo$.

Therefore every automaton with a fixed initial state $q$ (such
automata are called \emph{initial}) defines a transformation
$\lambda(q, \cdot)$ on the sets of finite and infinite words $\xs$ and
$\xo$. \emph{In fine}, we get a transformation of the set $\xs\cup\xo$
which will be written $\autom_q$.  The image of a word $x_1x_2\ldots$
under the map $\autom_q$ can be easily found using the Moore diagram
of the automaton. One must find a directed path starting in state $q$
with consecutive labels $(x_1, y_1), (x_2, y_2), \ldots$. Such a path
will be unique and the word $(x_1x_2\ldots)^{\autom_q}$ will be equal
to $y_1y_2\ldots$.

\begin{defi}
  \label{defi:finstate}
  A transformation $g:\xo\to\xo$ (or $g:\xs\to\xs$) is
  \emph{finite-state} if it is defined by a finite automaton.
\end{defi}

The product $\autom_q\cdot\mathcal{B}_p$ of two transformations
defined by automata $\autom=\langle Q_1, \lambda_1, \pi_1\rangle$
and $\mathcal{B}=\langle Q_1, \lambda_2, \pi_2\rangle$, with
respective initial states $q\in Q_1$ and $p\in Q_2$, is defined as the
\emph{composition} $\autom*\mathcal{B}$ of the automata with initial
state $(q, p)$. The composition $\autom*\mathcal{B}$ is defined as the
automaton $\langle Q_1\times Q_2, \lambda, \pi\rangle$, where
\begin{eqnarray*}
  \lambda((s, r), x) & = & \lambda_2(r, \lambda_1(s, x));\\
  \pi((s, r), x)     & = & \left(\pi_1(s, x), \pi_2(r, \lambda_1(s, x))\right).
\end{eqnarray*}

Thus the product of two finite-state transformations is again a
finite-state transformation.

An automaton is \emph{invertible} if each of its states defines an
invertible transformation of the set $\xo$ (or equivalently of the set
$\xs$). This condition holds if and only if for every state $q$ the
transformation $\lambda(q, \cdot):X\to X$ (i.e., the transformation
$\autom_q$ restricted onto the set $X^1\subset\xs$) is invertible
(see~\cite{eil,sush:avt_en,grineksu_en}).

If the automaton $\autom=\langle Q, \lambda, \pi\rangle$ is
invertible, then the inverse transformation to the transformation
$\autom_q$ is defined as the \emph{inverse automaton}
$\autom^*=\langle Q, \lambda^*, \pi^*\rangle$ with initial state
$q$, where
\begin{eqnarray*}
  \lambda^*(s, x) & = & y\\
  \pi^*(s, x) & = & \pi(s, y);
\end{eqnarray*}
here $y\in X$ is such that $\lambda(s, y)=x$ (such a $y$ exists and is
uniquely defined, since the automaton $\autom$ is invertible).

It also follows from these formul\ae\ that if an invertible
transformation is finite-state, then its inverse is also finite-state.

For more facts on automatic transformations
see~\cite{grineksu_en,sidki_monogr,sush:avt_en}.

\subsection{The complete automaton of a self-similar action}
Suppose that the group $G$ acts by a self-similar action on the space
$\xo$. Then this action defines an automaton over the alphabet $X$
with the set of states $G$ with the output and the transition
functions $\lambda$ and $\pi$ defined in such a way that
\[g(xw)=\lambda(g, x)w^{\pi(g, x)}\]
for all $w\in\xo$.

The obtained automaton $\autom$ has the property that the
transformation $\autom_g$ coincides with the action of the element
$g$. The automaton $\autom$ is called the \emph{(complete) automaton}
of the self-similar action.

Therefore the notion of a self-similar action of a group $G$ on the
space $\xo$ can be defined in terms of the automata theory in the
following way:
\begin{reformulate}{defi}{\ref{defi:ssact}''}
  An action of a group $G$ on the space $\xo$ is \emph{self-similar}
  if there exists an automaton $\autom=\langle G, \lambda,
  \pi\rangle$ such that
  \[w^{\autom_g}=w^g\]
  for all $w\in\xo$.
\end{reformulate}

We have the following obvious characterization of the automata which
are associated with the self-similar actions:
\begin{proposition}
  \label{pr:stcl2}
  Let $G$ be a group. An automaton $\autom=\langle G, \lambda,
  \pi\rangle$ is associated with a self-similar action of the group
  $G$ on the space $\xo$ if and only if $\autom_1$ is the identical
  transformation (here $1$ is the identity of the group $G$) and for
  all $g_1, g_2\in G$ we have
  \[\autom_{g_1}\autom_{g_2}=\autom_{g_1g_2}.\]
\end{proposition}

\subsection{Groups generated by automata}
The complete automaton of an action is infinite for infinite groups,
so it is not very convenient to define a group by its complete
automaton.  A better way is to define the group \emph{generated} by an
automaton in the following way.

Let $\autom=\langle Q, \lambda, \pi\rangle$ be an invertible
automaton. Then every transformation $\autom_q$ defined by an initial
state $q\in Q$ will be invertible. The group (or the group action) generated by the transformations
$\autom_q$ is called the \emph{group  generated by the automaton $\autom$} and is denoted
$G(\autom)$.

\begin{proposition}
  An action of a group on the set $\xo$ is self-similar if and only if
it is generated by an automaton.
\end{proposition}

The groups of the form $G(\autom)$, where $\autom$ is a finite
automaton, are the most interesting; we give here several
examples. These groups were introduced in Subsection~\ref{ss:exss}. More
examples will be given in Subsection~\ref{ss:exmono}.

\paragraph{The adding machine} action of the cyclic group is the action generated by the
automaton shown on Figure~\ref{fig:admach}. Namely, the state,
corresponding to the left vertex of the Moore diagram defines the adding
machine transformation $a$. The state corresponding to the right
vertex defines the trivial transformation. Thus the group defined by
the automaton is the adding machine action of the group $\Z$.

\begin{figure}[ht]
  \begin{center}
    \psfrag{1r}[rc]{$(1, 1)$} \psfrag{b}[cc]{$b$}
    \psfrag{0b}[cb]{$(0, 0)$} \psfrag{a}[cc]{$a$}
    \psfrag{10}[cb]{$(1, 0)$} \psfrag{01}[ct]{$(0, 1)$}
    \psfrag{1}[cc]{$1$} \psfrag{0l}[lb]{$(0, 0)$}
    \psfrag{1l}[lt]{$(1, 1)$}
    \includegraphics{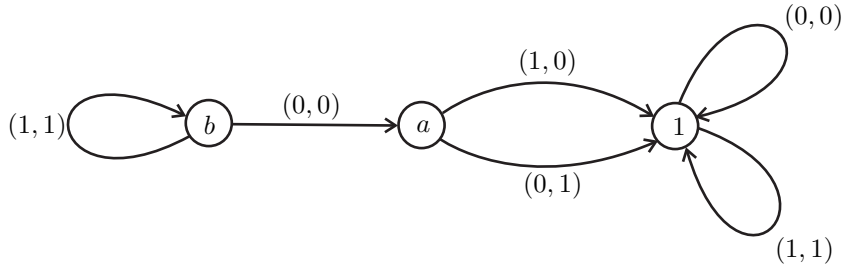}
  \end{center}
  \caption{An automaton generating $\mathbb{D}_\infty$}
  \label{fig:dihedr}
\end{figure}

\paragraph{The dihedral group} is generated by the automaton shown on
Figure~\ref{fig:dihedr}.

\begin{figure}[ht]
  \begin{center}
    \psfrag{a}[cc]{$a$} \psfrag{b}[cc]{$b$} \psfrag{c}[cc]{$c$}
    \psfrag{d}[cc]{$d$} \psfrag{e}[cc]{$1$} \psfrag{00}{$(0, 0)$}
    \psfrag{u}{$(1, 1)$} \psfrag{0}[lt]{$(0, 0)$} \psfrag{1}[rb]{$(1,
      1)$} \psfrag{01}[r]{$(0, 1)$} \psfrag{10}[l]{$(1, 0)$}
    \includegraphics{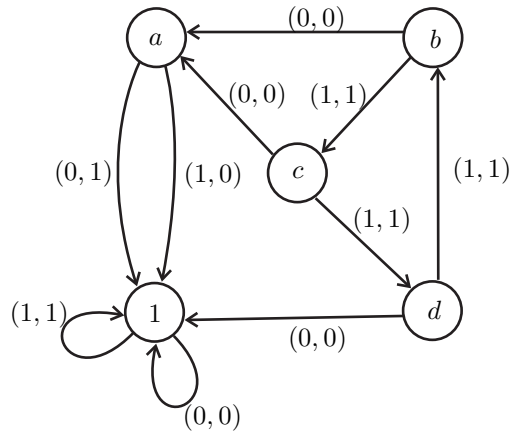}
  \end{center}
  \caption{The automaton generating the Grigorchuk group}
  \label{fig:grig}
\end{figure}

\paragraph{The Grigorchuk group} is defined by the automaton shown on
Figure~\ref{fig:grig}.

\subsection{Trees}
A \emph{tree} is a connected simplicial graph without cycles.

We consider only \emph{locally finite} trees, i.e., trees in which
every vertex has a finite degree.

A rooted tree is a tree with a fixed vertex called its \emph{root}. An
isomorphism of rooted trees $f:T_1\to T_2$ is an isomorphism of the
trees which maps the root of the tree $T_1$ to the root of the tree
$T_2$.

The vertices of a rooted tree are naturally partitioned into
\emph{levels}. If the distance between a vertex $v$ and the root is
equal to $n$ then we say that the vertex belongs to the $n$th level.
In particular, the $0$th level contains only the root.

A rooted tree is \emph{spherically homogeneous} (or \emph{isotropic}) if all the vertices
belonging to the same level have the same degree.

A spherically homogeneous tree is uniquely defined (up to isomorphism) by
its \emph{spherical index}. This is the sequence $(m_0,m_1,\ldots)$,
where $m_0$ is the degree of the root of the tree and $m_n+1$ is the
degree of the vertices of the $n$th level.

The regular $n$-ary tree is the spherically homogeneous tree with the
spherical index $(n, n, n, \ldots)$.

\paragraph{Example.} The set $\xs$ of finite words over the alphabet
$X$ is a vertex set of a naturally-defined rooted tree.  Namely, the
root of this tree is the empty word $\emp$, and two words are
connected by an edge if and only if they are of the form $v$ and $vx$
for some $x\in X$ and $v\in\xs$. We denote this rooted tree by $T(X)$.

The rooted tree $T(X)$ is spherically homogeneous with spherical index
$(d, d, d, \ldots)$, where $d=|X|$. Any regular $d$-ary rooted tree is
isomorphic to the tree $T(X)$.

The $n$th level of the tree $T(X)$ coincides with the set $X^n$.

An \emph{end} of a rooted tree $T$ is an infinite sequence of pairwise
different vertices (an infinite simple path) $v_0,v_1, v_2, \ldots$
such that $v_0$ is the root of the tree and for every $i$, the
vertices $v_i$ and $v_{i+1}$ are adjacent. The vertex $v_n$ will then
belong to the $n$th level of the tree.

The \emph{boundary} $\partial T$ of the tree $T$ is the set of all of its
ends. For every vertex $v$ denote by $\partial T_v$ the set of all the
ends passing through $v$. The sets $\partial T_v$ form a basis of
neighborhoods for the natural topology on $\partial T$. In this
topology the space $\partial T$ is totally discontinuous and compact.

Every automorphism of a rooted tree $T$ acts naturally on its boundary
and it directly follows from the definitions that it acts on
$\partial T$ by homeomorphisms.

In the case of the tree $T(X)$, every end has the form $(\emp, x_1,
x_1x_2, \ldots)$, and can be identified with the infinite word
$x_1x_2\ldots\in\xo$. Thus the boundary $\partial T(X)$ is naturally
identified with the space $\xo$. It is easy to see that this
identification agrees with the topologies on these sets, since
$\partial T(X)_v$ is identified with the cylindrical set $v\xo$.

\subsection{Action on rooted trees}
\label{ss:actrt}
It follows from the definition of the action of an initial automaton
$\autom_q$ on finite words that for any $v\in\xs, x\in X$ we have
$(vx)^{\autom_q}=v^{\autom_q}y$ for some $y\in X$. Thus the adjacent
vertices of the tree $T(X)$ are mapped onto adjacent vertices.
Therefore, if the transformation $\autom_q$ is invertible, then
$\autom_q$ defines an automorphism of the rooted tree $T(X)$.

More generally, the following holds:

\begin{proposition}
  \label{prop:automtree}
  A bijection $f:\xs\to\xs$ is defined by an automaton if and only if
  it induces an automorphism of the rooted tree $T(X)$.
  
  The action of an automorphism of the rooted tree $T(X)$ on the
  boundary $\partial T(X)=\xo$ coincides with the action of the
  respective initial automaton on the space $\xo$.
\end{proposition}

Let us denote by $\Sigma(X)$ the symmetric group of permutations of
the set $X$.  Every permutation $\alpha\in\Sigma(X)$ can be extended
to an automorphism of the whole tree $T(X)$ in a standard way:
\[(xw)^\alpha=x^\alpha w,\]
where $w\in\xs$ is arbitrary. This extension defines a canonical
embedding of the symmetric group $\Sigma(X)$ into the automorphism
group $\aut T(X)$ of the tree $T(X)$

Suppose $g$ is an automorphism of the tree $T(X)$.  It fixes the root
$\emp$ and it permutes the first level $X^1$ of the tree. Let
$\alpha_\emp\in\Sigma(X)$ be the permutation of the set $X^1$ induced
by $g$.  Then the automorphism $g\alpha_\emp^{-1}\in\aut T(X)$ fixes
the points of the first level (here $\alpha_\emp$ is identified with
its canonical extension).

\begin{defi}
  The set of all automorphisms which fix pointwise the first level is
  a subgroup called the \emph{first level stabilizer} and is written
  $\stab_1$. It is a normal subgroup of index $n!$ in the group $\aut
  T(X)$.
  
  In general, the \emph{$n$th level stabilizer} $\stab_n$ is the
  subgroup of those elements of the $\aut T(X)$, which fix all the
  elements of the $n$th level $X^n$ of the tree $T(X)$.  The $n$th
  level stabilizer is also a normal subgroup of finite index and the
  quotient $\aut T(X)/\stab_n$ is isomorphic to the automorphism group
  of the finite rooted subtree of $T(X)$ consisting of the first $n+1$
  levels.
  
  Let $G$ be a subgroup of $\aut T(X)$. We say $G$ is \emph{recurrent}
  if $G\cap\stab_n$ is a subgroup of $G^X$, where the latter group is
  a direct product of $|X|$ copies of $G$ acting independently on the
  level-$1$ subtrees of $T(X)$.
\end{defi}

The first level stabilizer is isomorphic to the direct product of
$d=|X|$ copies of the group $\aut T(X)$. Namely, every automorphism
$g\in \stab_1$ acts independently on the $d$ subtrees $T_x$, $x\in X$,
where $T_x$ is the subtree rooted at $x$, with the set of vertices $x\xs$.
Every subtree $T_x$ is isomorphic to the whole tree $T(X)$, with the
isomorphism given by restriction of the shift $\shift$ to the set
$x\xs$. Let us consider the restriction of the automorphism $g$ to the
subtree $T_x$ and conjugate it with the isomorphism $\shift:T_x\to
T(X)$ to get an automorphism of the whole tree $T(X)$. Let us denote
the obtained automorphism by $g|_x$.  It is easy to see that we have
\[(xw)^g=x(w)^{g|_x}\]
for all $w\in\xs$. Thus this notion of restriction agrees with the
one introduced before for self-similar groups.

\paragraph{Notation.} Set $X=\{x_1, x_2, \ldots, x_d\}$. Write
$g|_{x_i}=g_i$. We have a map $\Psi:g\mapsto (g_1, g_2, \ldots g_d):
\stab_1\to\aut T(X)^d$. It is easy to see that the map $\Psi$ is an
isomorphism.

Identifying the first level stabilizer with the direct product $\aut
T(X)^d$ we write $g=(g_1, g_2, \ldots g_d)$. In general, every element
of the group $\aut T(X)$ can be written as a product $(g_1, g_2,
\ldots g_d)\alpha_\emp$ of an element of the stabilizer and an element
of the symmetric group $\Sigma(X)$.

Therefore, the group $\aut T(X)$ is isomorphic to a semi-direct
product $\aut T(X)^X\rtimes\Sigma(X)$ with the natural action of $\Sigma(X)$ on the multiples
of the direct product $\aut T(X)^X$.

\begin{defi}
 Let $H$ be a group and let $G$ be a group acting on a set $M$. A \emph{(permutational) wreath product}
$H\wr G$ is the semi-direct product $H^M\rtimes G$, where $G$ acts on $H^M$ by the permutations
of the direct multiples coming from its action on $M$.
\end{defi}

Permutational wreath products are also called \emph{non-standard wreath products}. The
standard wreath product $H\wr G$ is the permutational wreath product respectively to the
regular action on $G$ on itself by right multiplication.

Hence, the automorphism group $\aut T(X)$ is isomorphic to the permutational wreath product $\aut
T(X)\wr\Sigma(X)$.

We can proceed further in this manner. Denote by $\alpha_x$ the
permutation of the points of the first level, defined by the
automorphism $g|_x$. Then $g|_x=(g|_{xx_1}, g|_{xx_2}, \ldots,
g|_{xx_d})\alpha_x$. In general, denote by $\alpha_v\in\Sigma(X)$ the
permutation of the first level of the tree, defined by the restriction
$g|_v$. In this way we get a labeling of the tree $T(X)$, where each
vertex $v$ is labeled by a permutation $\alpha_v$ from $\Sigma(X)$.
This labeled tree is called the \emph{portrait} of the automorphism
$g$.

The portrait defines the automorphism uniquely. The action of an
automorphism $g$ is computed via its portrait using the formula
\[
(a_1a_2\ldots)^g=(a_1)^{\alpha_\emp}(a_2)^{\alpha_{a_1}}(a_3)^{\alpha_{a_1a_2}}
(a_4)^{\alpha_{a_1a_2a_3}}\ldots.
\]

Let us recall here the definition of branch and weakly branch groups;
for more details see~\cite{grigorchuk:branch}.

\begin{defi}
  \label{defi:branch}
  Let $G$ be a level-transitive automorphism group of the rooted tree
  $T(X)$. The \emph{rigid stabilizer} $\rist_G(v)$ of a vertex
  $v\in\xs$ in the group $G$ is the set of all elements of the group
  $G$ which fix all the vertices of $T(X)$, except perhaps the
  vertices of the form $vu$, $u\in\xs$.  The \emph{$n$th level rigid
    stabilizer} $\rist_G(n)$ of the group $G$ is the subgroup
  generated by the rigid stabilizers $\rist_G(v)$ of all the
  vertices of the $n$th level (i.e., all the vertices $v\in X^n$).
  
  The group $G$ is \emph{weakly branch} if none of its rigid
  stabilizers $\rist_G(v)$ is trivial. It is \emph{branch} if
  $\rist_G(n)$ has finite index in $G$ for all $n\in\N$. The group is
  \emph{tough} if if it is not weakly branch.
  
  Furthermore, we say $G$ is \emph{regular (weakly) branch} if there
  is a subgroup $K<G$ such that $K^X<K$ as a geometric embedding
  induced by restriction to the level-$1$ subtrees in $T(X)$, and $K$
  has finite index in $G$ (respectively, is non-trivial).
\end{defi}

If the group is weakly branch, then all of its rigid stabilizers are
infinite. If it is tough, then for all $n$ big enough the $n$th level
rigid stabilizer is trivial.

\section{Iterated monodromy groups}
\label{s:img}
\subsection{Definitions and main properties}
Here we present a class of examples of self-similar group actions,
which are associated with topological dynamical systems. These
examples show the close relation between the classical self-similar
fractals, like the Julia sets of polynomials and the self-similar
group actions. We will show later how the Julia set of a rational
function can be reconstructed from its iterated monodromy group
(Subsection~\ref{ss:exlimit}).

Let $M_1$ and $M_2$ be arcwise connected and locally arcwise connected
topological spaces. A continuous map $f:M_1\to M_2$ is a \emph{local
  homeomorphism} at the point $p\in M_2$ if there exists a
neighborhood $U$ of the point $p$ such that $f^{-1}(U)$ is a disjoint
union of sets $U_i$ such that $f:U_i\to U$ is a homeomorphism.

We have the following classical result (see~\cite{massey:algtop}):

\begin{lemma}
  \label{l:paths}
  Let $f:M_1\to M_2$ be a local homeomorphism in every point of an arcwise
  connected open set $U\subset M_2$. Let $\gamma$ be a continuous path
  in $U$ starting in a point $p\in U$.  Then for every $x\in M_1$ such
  that $f(x)=p$ there exists a unique path $\gamma'$ in $M_1$ starting
  in $x$ such that $f(\gamma')=\gamma$.
\end{lemma}

Suppose $M$ is an arcwise connected and locally arcwise connected
topological space and let $f:M\to M$ be a branched $d$-fold
self-covering. In other words, there exists a subset $R\subset M$
(called the \emph{set of branching points}) such that for every
point $x\in M\setminus R$ the map $f$ is a local homeomorphism in $x$.

By $f^n$ we denote the $n$th iterate of the mapping $f$.  Let
$P=\overline{\cup_{n=0}^\infty f^n(R)}$ be the closure of the union of
the forward orbits of the branching points.  The set $P$ is called
the \emph{set of postcritical points}.  Then every preimage of a point
from $M\setminus P$ also belongs to $M\setminus P$, and $f$ is a local
homeomorphism at every point of $M\setminus P$.

An important example of a branched covering is the branched covering
of the Riemann sphere $M=\hat\C=\C\cup\{\infty\}$ defined by a rational
function $f$. The map $f$ is a $d$-fold branched covering, for $d$ equal to the
degree of the rational function, i.e., $d=\max(\deg p, \deg q)$, where $p, q\in\C[z]$
are such that $p/q=f$ is a reduced fraction. The set of branching points is in this case the set of
\emph{critical values} of the rational function $f$, i.e., the values of
the function at the critical points.

We impose throughout this section the condition that the $d$-fold branched covering $f:M\to M$ is such that
the space $M\setminus P$ is arcwise connected. If furthermore $P$ is finite, we
say $f$ is \emph{postcritically finite}. See the paper~\cite{DH:Thurston} and Appendix B
of the book~\cite{McMul:renorm} for an interesting Thurston's criterion for a
postcritically finite branched covering of a sphere to be defined by a rational function.

Let $t\in M\setminus P$ be an arbitrary point. Then it has $d^n$
preimages under $f^n$ for every $n\in\N$. Denote by $T$ the formal
disjoint union of the sets $f^{-n}(t)$ for $n\geq 0$, (where
$f^{-0}(t)=\{t\}$ and $f^{-n}(t)$ denotes the preimage of the
point $t$ under the map $f^n$). Formally the set $T$ can be defined as
\[T=\bigcup_{n=0}^\infty f^{-n}(t)\times \{n\},\]
so that its element $(z, n)$ is the point $z$ seen as an element of
the set $f^{-n}(t)$.

The set $T$ is a vertex set of a naturally-defined $d$-regular rooted
tree with the root $(t, 0)$ in which a vertex $(z, n)\in T$ is
connected with the vertex $(f(z), n-1)$. Let us call the tree $T$ the
\emph{preimages tree} of the point $t$.

Let now $\gamma$ be a loop in $M\setminus P$ based at $t$, i.e., a
path starting and ending at $t$.  Then, by Lemma~\ref{l:paths}, for
every element $v=(z, n)$ of the preimage tree $T$ there exists a
single path $\gamma_v$, which starts in $z$ and is such that
$f^n(\gamma_v)=\gamma$. Denote by $z^\gamma$ the end of the path
$\gamma_v$. Then we obviously have $f^n\left(z^\gamma\right)=t$, so
the element $v^\gamma=\left(z^\gamma, n\right)$ also belongs to the
preimage tree.

We have the following proposition:
\begin{proposition}
  \label{prop:img}
  The map $v\mapsto v^\gamma$ is an automorphism of the preimage tree,
  which depends only on the homotopy class of $\gamma$ in $M\setminus
  P$. The set of all such automorphisms is a group, which is a
  quotient of the fundamental group of the space $M\setminus P$. Up to
an  isomorphism, this group does not depend on the choice of the
  base point $t$.
\end{proposition}

\begin{defi}
  \label{defi:img}
  The group from Proposition~\ref{prop:img} is called the
  \emph{iterated monodromy group (i.m.g.)}  of the map $f$, and is
  written $\IMG(f)$.
\end{defi}

The term \emph{iterated monodromy group} comes from the fact that the
quotient of this group by the stabilizer of the $n$th level is the
monodromy group of the mapping $f^n$.

The preimage tree $T$ is $d$-regular, so that it can be identified with the
tree $T(X)$ for an alphabet $X$ of cardinality $d$.  There is no
canonical identification, but we will define a class of natural
identifications using paths in the space $M\setminus P$.

Let us take an alphabet $X$ with $d$ letters and a bijection $\Lambda:X\to f^{-1}(t)$.
For every point $x\in f^{-1}(t)$ choose a path $\ell_x$ in $M\setminus P$, starting at $t$ and ending
at $\Lambda(x)$.

Define a map $\Lambda:T(X)\to T$ inductively by the rules:
\begin{enumerate}
\item
$\Lambda(\emp)=(t, 0)$,
\item
for every $n\ge 1$, $v\in X^n$ and $x\in X$ the point $\Lambda(xv)$ is $(z, n+1)$,
where $z$ is the end of the path $\gamma$ which starts at $\Lambda(v)$ and is such that
$f^n(\gamma)=\ell_x$.
\end{enumerate}

\begin{proposition}
The constructed map $\Lambda: T(X)\to T$ is an isomorphism of the rooted trees.
\end{proposition}

\begin{defi}
  The \emph{standard action} of the i.m.g.\ on the tree $T(X)$ is
  obtained from its action on the preimage tree $T$ conjugating it by
  the isomorphism $\Lambda:T(X)\to T$.
\end{defi}

\begin{proposition}
  \label{prop:imgss}
 The standard action of an iterated monodromy group is self-similar. More precisely,
if $\gamma$ is a loop based at $t$ and $x\in X$ is a letter,  then, respectively to the standard action, for
every $v\in\xs$ we have
\begin{equation}
\label{eq:natur}
(xv)^\gamma=y\left(v^{\ell_x\gamma_x\ell_y^{-1}}\right),
\end{equation}
where $\gamma_x$ is the preimage of $\gamma$ starting at $\Lambda(x)$ and $y$ is such that
$\Lambda(y)$ is the end of $\gamma_x$ (i.e., $x^\gamma=y$).
\end{proposition}

\begin{figure}[h]
\begin{center}
\psfrag{g}[tr]{$\gamma$} \psfrag{lx}[br]{$\ell_x$}
\psfrag{ly}[tr]{$\ell_y$} \psfrag{gx}[bl]{$\gamma_x$}
\psfrag{x}[bl]{$x$} \psfrag{y}[tl]{$y$} \psfrag{t}[br]{$t$}
\includegraphics{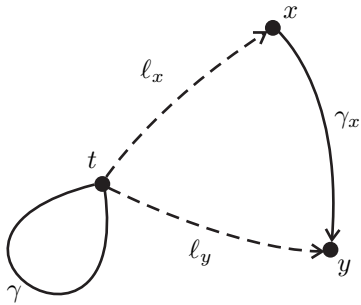}
\end{center}
\caption{The recurrent formula}\label{fig:recur}
\end{figure}

Note, that $\ell_x\gamma_x\ell_y^{-1}$ is a loop based at $t$ (see Figure~\ref{fig:recur}).

In many cases the action of the iterated monodromy group is generated
by a finite automaton.  This is the case, for instance, when the map $f$ is expanding with respect to some
Riemann metric on $M$, or for postcritically finite rational mappings of the complex sphere, as
follows from Theorem~\ref{th:limjul}.

\subsection{Examples of iterated monodromy groups}
\label{ss:exmono}
We give in this subsection some examples of iterated monodromy
groups. Most of these groups are described for the first time here,
and we will only mention their most elementary properties~--- a more
detailed study will be given in~\cite{bartholdi:imgexamples} and~\cite{nek:img}.

In this section we will use very often the notation from
Subsection~\ref{ss:actrt} coming from the decomposition of the group
$\aut T(X)$ into the semidirect product $\aut T(X)^d\rtimes\Sigma(X)$.

\paragraph{The adding machine.} Let $\mathbb{T}$ be the circle
$\{z\in\C: |z|=1\}$. The map $f:z\mapsto z^2$ is a two-fold
self-covering of the circle $\mathbb{T}$. Let us compute the iterated
monodromy group of the covering $f:\mathbb{T}\to\mathbb{T}$.

Let us chose the base-point $t$ equal to $1$. It has two preimages:
itself and $-1$. So we can take $\ell_1$ equal to the trivial path in the
point $1$ and $\ell_{-1}$ equal to the upper semicircle starting at $1$
and ending in $-1$.

The fundamental group of the circle is the cyclic group generated by
the loop $\tau$ which starts in $1$ and goes around the circle once in
the positive direction. So it is sufficient to compute the action of
the element $\tau$ on the tree $T(X)$. The path $\tau$ has two
preimages.  One is the upper semicircle $\tau_1$ starting in $1$ and
ending in $-1$ (i.e., equal to the path $\ell_{-1}$ defined above),
another is the lower semicircle $\tau_{-1}$ starting in $-1$ and
ending in $1$. Thus $\tau$ acts on the first level of the tree $T(X)$
by the transposition. So using Equation~\eqref{eq:natur} we get, in
notation of Subsection~\ref{ss:actrt},
\[
\tau=(\ell_1\tau_1\ell_{-1}^{-1}, \ell_{-1}\tau_{-1}\ell_1^{-1})\sigma=(1, \tau)\sigma,
\]
since the path $\ell_1\tau_1\ell_{-1}^{-1}$ is trivial and
$\ell_{-1}\tau_{-1}\ell_1^{-1}$ is equal to $\tau$. Therefore, $\tau$ acts
on the tree $T(X)$ as the adding machine.

\paragraph{The torus.} The above example can be obviously generalized
to tori in the following way:

Let $\mathbb{T}^n=\R^n/\Z^n$ be the $n$-dimensional torus. Let $A$ be
a $n\times n$-matrix with integral entries and with determinant equal
to $d>1$.  Then the linear map $A$ on the space $\R^n$ induces an
$d$-fold (non-branched) self-covering of the torus $\mathbb{T}^n$. Since the
fundamental group of the torus is the free abelian group $\Z^n$, the
i.m.g.\ of such coverings are abelian. More precisely, the following
proposition holds.

\begin{proposition}
\label{pr:torusimg}
  Let $A$ by the matrix defining the covering
  $f:\mathbb{T}^n\to\mathbb{T}^n$. Then the iterated monodromy group of
$f$ is the quotient  $\Z^n/H$, where 
$$
H=\bigcap_{n\ge 1} A^n\left(\Z^n\right).
$$
\end{proposition}

The group $H$ is trivial if and only if no eigenvalue of the matrix $A^{-1}$ is an algebraic integer
(see~\cite{neksid} Proposition~4.1 and~\cite{cuntz_rep} Proposition~10.1).

The corresponding actions of $\Z^n$ on rooted trees are studied in
Subsection~\ref{ss:abel}.

\paragraph{Chebyshev polynomials and example of S.~Latt\`es.}
Consider the Chebyshev polynomials
$T_d(z)=\cos(d\arccos z)$, satisfying the recursion
\[T_0(z)=1,\quad T_1(z)=z,\quad T_d(z)=2zT_{d-1}-T_{d-2}.\]
Then $T_d$ is an even or odd polynomial of degree $d$.

\begin{proposition}
\label{pr:cheb}
The group  $\IMG(T_d)$ is infinite dihedral for all $d\ge2$.
\end{proposition}

\begin{proof}
 We have the following commutative diagramm:
\begin{equation}
\label{eq:cheb}
\xymatrix{{\C} \ar[r]^{dz} \ar[d]_{\cos z} & {\C} \ar[d]^{\cos z}\\
{\C} \ar[r]_{T_d(z)} & {\C}
}
\end{equation}
Let $t\in (-1, 1)$ be a basepoint not belonging to the postcritical
set of $T_d$. From the diagram~\eqref{eq:cheb} follows that
$$
T_d^{-n}(t)=\left\{\cos\left(\pm\frac{\alpha+2\pi k}{d^n}\right): k\in\Z\right\}=
\left\{\cos\left(\frac{\alpha+2\pi k}{d^n}\right): k=0, 1, \ldots d^n-1\right\},
$$
where $\alpha=\arccos t$.

Let $\gamma$ be an arbitrary loop at $t$ not passing through the
branching points of the functions $T_d(z)$ and $\cos z$. The function
$\cos z$ is even and $2\pi$-periodic, thus the preimages of $\gamma$
under $\cos z$ are either paths starting at $\alpha+2\pi k$ and ending
at $\alpha+2\pi (k+l)$ and paths starting at $-\alpha-2\pi k$ and
ending at $-\alpha-2\pi (k+l)$, where $l$ is fixed and $k\in\Z$, or
paths starting at $\alpha+2\pi k$ and ending at $-\alpha+2\pi (k+l)$
and paths starting at $-\alpha-2\pi k$ and ending at $\alpha-2\pi
(k+l)$, where $l$ is fixed and $k\in\Z$.

This implies that the preimages of $\gamma$ under $T_d^n$ either start
at $\cos\left(\frac{\alpha+2\pi
    k}{d^n}\right)=\cos\left(\frac{-\alpha-2\pi k}{d^n}\right)$ and
end at $\cos\left(\frac{\alpha+2\pi
    (k+l)}{d^n}\right)=\cos\left(\frac{-\alpha-2\pi
    (k+l)}{d^n}\right)$, or they start at $\cos\left(\frac{\alpha+2\pi
    k}{d^n}\right)=\cos\left(\frac{-\alpha-2\pi k}{d^n}\right)$ and
end at $\cos\left(\frac{\alpha-2\pi
    (k+l)}{d^n}\right)=\cos\left(\frac{-\alpha+2\pi
    (k+l)}{d^n}\right)$, where $l$ is fixed and $k\in\Z$.

Therefore, the iterated monodromy group $\img(T_d)$ is isomorphic to
the group of affine functions $\{z+l, -z+l: l\in\Z\}$ under
composition. This group is isomorphic to the infinite dihedral group
$\mathbb{D}_\infty$.
\end{proof}

An explicit computation shows that the associated standard action of
$\mathbb{D}_\infty$ on the tree $T(X)$, for $d$ odd, is generated by
two involutions $a$ and $b$, where
$$
a=(a, 1, 1, \ldots, 1)\sigma_1, \quad b=(1, 1, \ldots, 1, b)\sigma_2,
$$
where $\sigma_1$ is the permutation $(2, 3)(4, 5)\ldots (d-1, d)$
and $\sigma_2=(1, 2)(3, 4)\ldots (d-2, d-1)$ and for $d$ even by
$$
a=(1, 1, \ldots, 1)\sigma_1,\quad b=(a, 1, \ldots, 1, b)\sigma_2,
$$
where $\sigma_1=(1, 2)(3, 4)\ldots (d-1, d)$ and $\sigma_2=(2,
3)(4, 5)\ldots (d-2, d-1)$.  Here the alphabet $X$ is $\{1, 2, \ldots,
d\}$.

Also related are the following examples of
S.~Latt\`es~\cite{lattes:iter}.  Let $\Lambda$ be a lattice in $\C$,
and let $\alpha$ be a multiplier of $\Lambda$, i.e., some
$\alpha\in\C$ such that $\alpha\Lambda\subset\Lambda$. Then
$\C/\Lambda$ is a torus, and the affine function $\alpha \cdot z$
induces an $|\alpha|^2$-fold self-covering of $\C/\Lambda$. Note that
the iterated monodromy group of this self-covering is isomorphic to
$\Lambda$, as follows from Proposition~\ref{pr:torusimg}.

The Weierstrass elliptic function
$$
\wp(z)=\frac{1}{z^2}+\sum_{\omega\in\Lambda\setminus\{0\}}\left[\frac{1}{(z+\omega)^2}-\frac{1}{\omega^2}\right]
$$
is $\Lambda$-periodic and even (see for
example~\cite{lang:ellipt}), so that it we get a well-defined $2$-fold
covering $\wp:\C/\Lambda\to\widehat\C$, branched at the four points
$\frac12\Lambda/\Lambda$.  Then the function
$f(z)=\wp(\alpha\wp^{-1}(z))$ is rational of degree $|\alpha|^2$. The
dynamics of such maps were first studied by
S.~Latt\`es~\cite{lattes:iter}.

Let $\Lambda$ be a lattice in $\C$, and let $\alpha$ be a multiplier
of $\Lambda$, i.e.\ some $\alpha\in\C$ such that
$\alpha\Lambda\subset\Lambda$. Then $\C/\Lambda$ is a torus, and
$\wp:\C/\Lambda\to\widehat\C$ is a degree-$2$ covering, branched at
the four points $\frac12\Lambda/\Lambda$.  Consider
$f(z)=\wp(\alpha\wp^{-1}(z))$. Its Julia set is the entire
$\widehat\C$.  Then $f$ is a rational map of degree $|\alpha|^2$, and
$\IMG(f)$ is a $4$-generated group.

 For instance, for $\alpha=2$ the function $f$ is
$$
f(z)=\frac{z^4+\frac{g_2}{2}z^2+2g_3z+\frac{g_2^2}{16}}{4z^3-g_2z-g_3},
$$
where $g_2=60s_4$ and $g_3=140 s_6$ for $s_m=\sum_{\omega\in\Lambda, \omega\ne 0} \omega^{-m}$,
see~\cite{beardon:iter} p.~74.

There exists a lattice $\Lambda$ with given values of $g_2$ and $g_3$ if and only if
$g_2^3-27g_3^2\ne 0$ (see~\cite{lang:ellipt}, p.~39).
In particular, there exists a lattice $\Lambda$ such that $g_3=0$ and $g_2=4$, so that
\begin{equation}
\label{eq:lattes}
f(z)=\frac{(z^2+1)^2}{4z(z^2-1)}.
\end{equation}

The following proposition is proved in the similar way as Proposition~\ref{pr:cheb}.

\begin{proposition}
Let $\Lambda$ be a lattice in $\C$ and let $\alpha\in\C$ be such that
$\alpha\Lambda\subset \Lambda$ and $|\alpha|\ne 1$. Let a rational function $f\in\C(z)$ be such
that $\wp(\alpha z)=f(\wp(z))$. Then the iterated monodromy group $\img{f}$ is isomorphic to the
semi-direct product $\Z^2\rtimes (\Z/2\Z)$, where $\Z/2\Z$ acts on $\Z^2$ by the automorphism
$(a, b)\mapsto (-a, -b)$.
\end{proposition}

\paragraph{$\mathbf{z^2-1}$.} The critical points of the polynomial $z^2-1$ are $\infty$ and $0$ and the
the postcritical set is $P=\{0, -1, \infty\}$.

\begin{figure}[h]
\begin{center}
\psfrag{a}[bl]{$a$}\psfrag{b}[bl]{$b$} \psfrag{-1}[ct]{$-1$}
\psfrag{0}[ct]{$0$} \psfrag{l}[br]{$\ell_1$}
\psfrag{fa}[bl]{$f^{-1}(a)$} \psfrag{fb}[bl]{$f^{-1}(b)$}
\psfrag{1}[ct]{$1$}
\includegraphics{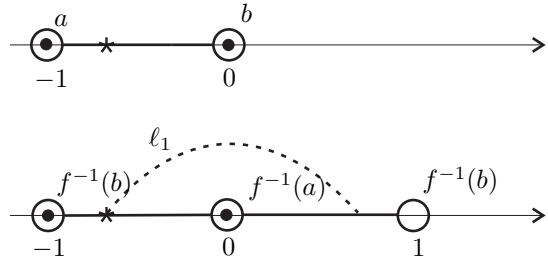}
\caption{Computation of the group $\img{z^2-1}$}\label{fig:min}
\end{center}
\end{figure}

 Choose as a basepoint the fixed point of the polynomial $t=\frac{1-\sqrt{5}}{2}$.
It has two preimages: itself, and $-t$. Choose the path $\ell_0$ to be trivial path at $t$ and
$\ell_1$ to be the path, connecting $t$ with $-t$ above the real axis, as on the lower part of
Figure~\ref{fig:min} (the point $t$ is marked there by a star). Let $a$ and $b$ be the generating elements of
$\img{z^2-1}$, defined by the small loops going in the positive direction around the points $-1$ and $0$,
respectively and connected to the basepoint by straight segments. The loops $a$ and $b$
are shown on the upper part of the figure.

\begin{figure}[h]
\begin{center}
\psfrag{a}[br]{$a$} \psfrag{b}[bl]{$b$} \psfrag{01}[bc]{$(0, 1)$}
\psfrag{00}[tc]{$(0, 0)$} \psfrag{10}[tr]{$(1, 0)$}
\psfrag{0}[tr]{$(0, 0)$} \psfrag{1}[tl]{$(1, 1)$}
\includegraphics{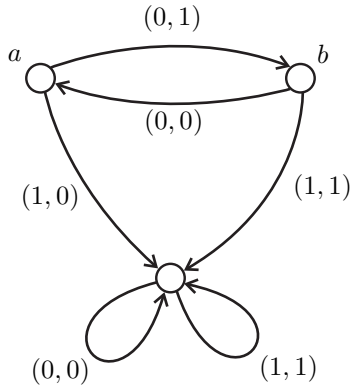}
\end{center}
\caption{The automaton generating the group $\img{z^2-1}$.}\label{fig:moor1}
\end{figure}

The preimages of the loops $a$ and $b$ are shown on the lower part of Figure~\ref{fig:min}. It follows that
$$
a=(b, 1)\sigma, \quad b=(a, 1),
$$
so that the group $\img{z^2-1}$ is generated by the automaton with the Moore diagram shown on
Figure~\ref{fig:moor1}.

The following properties of the group $\img{z^2-1}$ where proved by R.~Grigorchuk, A.~\.Zuk and L.~Bartholdi.

\begin{theorem}
\label{th:propimgmin}
 The group $\img{z^2-1}$
\begin{enumerate}
\item is torsion free;
\item has exponential growth (actually, the semigroup generated by $a$ and $b$ is free);
\item is just non-solvable, i.e., is not solvable, but all of its proper
  quotients are solvable;
\item has a presentation
$$
\img{z^2-1}=\left\langle a, b\left| \left[\left[a^{2^k}, b^{2^k}\right], b^{2^k}\right]=
\left[\left[b^{2^k}, a^{2^{k+1}}\right], a^{2^{k+1}}\right]=1, k=1, 2,\ldots\right.\right\rangle.
$$
\item
has solvable word and conjugacy problems;
\item
has no free non-abelian subgroups of rank 2.
\end{enumerate}
\end{theorem}

The properties of the group $\img{z^2-1}$ are similar to the
properties of the branch groups (see~\cite{grigorchuk:branch}) and of
the just non-solvable group of A.~Brunner, S.~Sidki and A.~Vieira
in~\cite{bsv:jns}.

\paragraph{More examples.} We present here the iterated monodromy groups
of quadratic rational maps with size of postcritical set at most $3$,
arranged in a table.  The algebraic structure of most of them is not
yet well understood.

\begin{tabular}{|p{10em}|p{8.5em}|p{17em}|}
\hline
$f(z)$ & standard action & comments \\
\hline \hline
$z^2$ & $\tau=(1, \tau)\sigma$ & The adding machine (see above)\\ \hline
$z^{-2}$ & $\mu=(1, \mu^{-1})\sigma$ & A conjugate of the adding machine.
The group, generated by the transformation $\mu$ and the adding machine $\tau$ was studied
by A.~Brunner, S.~Sidki and A.~C.~Vieira~\cite{bsv:jns}.\\ \hline
$z^2-1$ & $a=(b, 1)\sigma$,\newline  $b=(a, 1)$ & see
Theorem~\ref{th:propimgmin} above\\ \hline
$\frac{z^2-1}{z^2}$ & $a=(1,b)$,\newline  $b=(a^{-1},1)\sigma$ & \\ \hline
$z^2-2$ & $a=\sigma$, \newline $b=(a,b)$ & Isomorphic to $\mathbb{D}_\infty$, since $z^2-2$ is conjugate
to the Chebyshev polynomial $T_2$, (see above).\\ \hline
$z^2+c$,\newline with $c\in \R$ such that $c^3+2c^2+c+1=0$ &
$a=(1,b)\sigma$,\newline $b=(1,c)$,\newline $c=(a,1)$ & \\ \hline
$z^2+c$,\newline with $c^3+2c^2+c+1=0$, $c\notin\R$ & $a=(1,b)\sigma$, \newline
$b=(c,1)$,\newline  $c=(a,1)$ & The action is not conjugate to the previous one, though it
is not known if the groups are isomorphic. There closures in $\aut T(X)$ coincide. \\ \hline
$\frac{z^2-2}{z^2}$ & $a=(b,a)$,\newline $b=(b^{-1},a^{-1})\sigma$ &
The group is isomorphic to $\Z^2\rtimes(\Z/4)$, where $\Z/4\Z$ acts by
the matrix $(\begin{smallmatrix}0&-1\\1&0\end{smallmatrix})$.\\ \hline
$\frac{z^2-\phi^2}{z^2}$,\newline with $\phi=\frac{1+\sqrt5}2$ & 
$a=(b,1)$, \newline $b=(1,c)$,\newline $c=(a^{-1},b^{-1})\sigma$ & \\ \hline
$\frac{z^2-\phi^2}{z^2}$,\newline with $\phi=\frac{1-\sqrt5}2$ &
$a=(1,b)$,\newline $b=(1,c)$,\newline $c=(a^{-1},1)\sigma$ & \\ \hline
$\frac{z^2-1}{z^2+1}$ & $a=(1,b)\sigma$,\newline $b=(a,a^{-1})$ & \\ \hline
$\frac{z^2-1}{z^2-\omega}$,\newline with $\omega^3=1$, $\omega\ne 1$ & 
$a=(1,b)\sigma$,\newline $b=(c,1)$,\newline $c=(c^{-1}b^{-1},a^{-1})\sigma$ & \\ \hline
$z^2+i$ & $a=\sigma$,\newline $b=(a,c)$,\newline $c=(b,1)$ & This group has intermediate growth\\ \hline
\end{tabular}

\paragraph{Questions.}
\begin{enumerate}
\item When is the iterated monodromy group of a rational function torsion free?
\item When are these groups just non-solvable?
\item What do their presentations look like?
\item Do they all have solvable conjugacy problem?
\item Can any be non-amenable? or contain a free subgroup of rank $\ge
  2$?
\item Which rational functions have iterated monodromy groups of
  exponential growth?
\end{enumerate}

\paragraph{Galois groups.} The following construction is due to
R.~Pink (private communication).

Suppose the branched covering $f:M\to M$ is a polynomial mapping of
the complex sphere $M=\hat\C$.  Define the polynomials
$F_n(x)=f^n(x)-t$ over the field $\C(t)$. Let $\Omega_n$ be the
decomposition field of $F_n(x)$. Then $\Omega_n\subseteq\Omega_{n+1}$;
write $\Omega=\cup_{n\geq 0}\Omega_n$. Then the Galois group of the
extension $\C(t)\subset\Omega$ is isomorphic to the closure of the
group $\IMG(f)$ in $\aut T$, where $T$, the disjoint union of the
roots of $F_n(x)$, has the natural structure of a $\deg(f)$-regular
tree.

\section{Virtual endomorphisms}
\label{s:virt}
\subsection{Definitions}
\begin{defi}
  A \emph{virtual endomorphism} $\phi:G\prr G$ is a homomorphism from
  a subgroup of finite index $\dom\phi\leq G$ into $G$.
\end{defi}

The product $\phi=\phi_1\phi_2$ of two virtual endomorphisms is again a
virtual morphism with domain
\[
\dom\phi=\{g\in G:g\in\dom\phi_2, \phi_2(g)\in\dom\phi_1\}=\phi_2^{-1}
\left(\dom\phi_1\right).
\]

Let us fix a faithful self-similar action of a group $G$ on the space
$\xo$. For every $x\in X$ we denote by $G_x$ the stabilizer of the
one-letter word $x$ in the associated action of the group $G$ on
$\xs$. Define a map $\phi_x:G_x\to G$ by the formula
\[\phi_x(g)=g|_x.\]
It follows from Equation~\eqref{eq:restr} that $\phi_x$ is a
homomorphism from $G_x$ into $G$. The group $G_x$ is a subgroup of
finite index in $G$. The index is equal to the cardinality of $X$,
since we assume that the action is level-transitive.

Thus, given a self-similar action of a group $G$ for any $x\in X$ we
have a virtual endomorphism $\phi_x: G\prr G$. We call this
endomorphism the endomorphism \emph{associated with the self-similar
  action}.

For example, for the adding machine action, the associated virtual
endomorphism is the map $\Z\prr\Z:n\mapsto n/2$ with domain equal to
the set of even numbers.

Let again $\phi:G\prr G$ be any virtual endomorphism. Choose a right
transversal $T=\{g_0=1, g_1, g_2, \ldots, g_{d-1}\}$ to the subgroup
$\dom\phi$ (where $d=[G:\dom\phi]$), i.e., elements of $G$ such that
$G$ is a disjoint union of the sets $\dom\phi\cdot g_i$. Choose also a
sequence $C=(h_0=1, h_1, h_2, \ldots, h_{d-1})\in G^d$ of arbitrary
elements of the group. Define an automaton over the alphabet $X=\{x_0,
x_1, \ldots, x_{d-1}\}$, with the set of states equal to $G$, by the
equations
\begin{equation}\label{eq:sint}
\lambda(g, x_i)=x_j, \qquad \pi(g,
x_i)=h_i^{-1}\phi(g_igg_j^{-1})h_j,
\end{equation}
where $j$ is such that $g_igg_j^{-1}\in\dom\phi$.

The obtained automaton $\autom(\phi, T, C)$ will be called the
\emph{automaton defined by the virtual endomorphism $\phi$, the coset
  transversal $T$ and the sequence $C$}.

Using Proposition~\ref{pr:stcl2}, one can prove the following assertion:
\begin{proposition}
  \label{prop:generation}
  The automaton $\autom(\phi, T, C)$ defines a self-similar action of
  the group $G$ with the associated virtual endomorphism $\phi$.
\end{proposition}

A set $K\subseteq G$ is \emph{invariant} under a virtual endomorphism
$\phi$ if $K\subseteq\dom\phi$ and $\phi(K)\subseteq K$.

The kernel of a self-similar action can be determined using the
following description (see~\cite{nek:stab,neksid}):
\begin{proposition}
  \label{prop:kore}
 The kernel of a self-similar action of a group $G$ is the maximal
normal $\phi$-invariant subgroup of $G$, where $\phi$ is the virtual
endomorphism, associated with the action. It is equal to the group
$$
\bigcap_{n\ge 1}\bigcap_{g\in G} g^{-1}\cdot \dom\phi^n\cdot g.
$$
\end{proposition}

The above construction is universal, namely:
\begin{proposition}
  \label{prop:gen}
  Any self-similar action is defined by the associated virtual
  endomorphism $\phi_{x_0}$, a coset transversal $T=\{g_y: y\in X\}$
  and the sequence $C=\{h_y: y\in X\}$, where $x_0^{g_y}=y$ and
  $g_y|_{x_0}=h_y$.
\end{proposition}

\subsection{Recurrent actions and abstract numeration systems}

Let $G$ be a group with a self-similar action over the alphabet $X=\{x_0,x_1, \ldots, x_{d-1}\}$.

\begin{defi}
  \label{defi:recur}
  The self-similar action is \emph{recurrent} if the
associated virtual endomorphism $\phi_x$ is onto, i.e., if $\phi_x(\dom\phi_x)=G$.
\end{defi}

It follows from Proposition~\ref{prop:gen} and Equation~\eqref{eq:sint} that every recurrent action has
a \emph{digit set} in the sense of the following definition.

\begin{defi}
  \label{defi:digit}
A \emph{digit set} for the self-similar action is a set
$T=\{r_0=1, r_1, \ldots, r_{d-1}\}\subset G$ such that for every $g\in G$ we have
$x_0^{r_i}=x_i$ and $r_i|_{x_0}=1$.
  
  The self-similar action defined by a virtual endomorphism $\phi$ and
a digit set $T$ is the self-similar action defined by $\phi$, the
coset transversal $T$ and the sequence $C=(1, 1, \ldots, 1)$.
\end{defi}

Suppose now that the action of the group $G$ is recurrent
and let $T=\{r_0=1, r_1, \ldots, r_{d-1}\}$ be its digit
set. Then the action of the group $G$ can be interpreted in the
following way. Let $w=x_{i_1}x_{i_2}x_{i_3}\ldots\in\xo$ be an
infinite sequence. We put it in correspondence with a formal
expression
\[
w=\left.\left.\left.\cdots r_{i_4}\right)^{\phi^{-1}}r_{i_3}\right)^{\phi^{-1}}r_{i_2}\right)^{\phi^{-1}}r_{i_1},
\]
where $(\cdot)^{\phi^{-1}}$ is another notation for
$\phi^{-1}(\cdot)$.  In order to know the image of the sequence
$x_{i_1}x_{i_2}x_{i_3}\ldots$ under the action of an element $g\in G$,
we multiply the expression from the right by $g$ and then reduce it to
a similar form. There exists a unique index $j_1$ such that
$r_{i_1}g=\tilde g r_{j_1}$ for $\tilde g\in\dom\phi$, so it is
natural to write
\begin{eqnarray*}
  wg &=& \ldots r_{i_3})^{\phi^{-1}} r_{i_2})^{\phi^{-1}}r_{i_1}g\\
  &=& \ldots r_{i_3})^{\phi^{-1}} r_{i_2})^{\phi^{-1}}\tilde g r_{j_1}\\
  &=& \left.\ldots r_{i_3})^{\phi^{-1}} r_{i_2}\phi\left(\tilde g\right)\right)^{\phi^{-1}}r_{j_1}.
\end{eqnarray*}
Next, we determine the index $j_2$ such that $r_{i_2}\phi\left(\tilde
  g\right)=\tilde{\tilde g}r_{j_2}$ for $\tilde{\tilde g}\in\dom\phi$
and proceed further. It follows from formula~\eqref{eq:sint} that in
this way we will get correctly all the indices $j_1j_2\ldots$ of the
image $x_{j_1}x_{j_2}\ldots=(x_{i_1}x_{i_2}\ldots)^g$.  In the general
(non-recurrent) case the formul\ae\ are slightly more complicated.

This interpretation can be viewed as some sort of ``$\phi$-adic''
numeration system on the space $\xo$. In particular, in the case of
the binary adding machine we will get the usual dyadic numeration
system if we chose $r_0=0,r_1=1$ (recall that the virtual
endomorphism in this case is $\phi(n)=n/2$). So the elements of the
coset transversal $\{r_0,r_1, \ldots, r_{d-1}\}$ play the role of the
digits and $\phi^{-1}$ is the ``base'' of the numeration system.

\subsection{Contracting actions}
\begin{defi}
  \label{defi:contr}
  A self-similar action of a group $G$ is \emph{contracting} if there
  exists a finite set $\nuke\subset G$ such that for every $g\in G$
  there exists $k\in\N$ such that $g|_v\in\nuke$ for every word
  $v\in\xs$ of length $\geq k$.  The minimal set $\nuke$ with this
  property is called the \emph{nucleus} of the self-similar action.
\end{defi}

Obviously, every contracting action is finite state.

\paragraph{The adding machine.} Let $a$ be the adding machine
transformation of the space $\{0,1\}^\omega$ introduced
in Subsection~\ref{ss:exss}. Then for even $n$ we have
\[a^n|_0=a^{\frac{n}{2}}, \quad a^n|_1=a^{\frac{n}{2}},\]
and for odd $n$ we have
\[a^n|_0=a^{\frac{n-1}{2}},\quad a^n|_1=a^{\frac{n+1}{2}}.\]

It easily follows from this that the adding machine action of the
group $\Z$ is contracting with nucleus $\{-1, 0,1\}$.

\paragraph{The Grigorchuk group} is contracting with the nucleus
equal to the standard set of generators $\{1,a,b,c,d\}$.  In the
original paper~\cite{grigorchuk:80_en}, where the Grigorchuk group was
defined, the contraction of the group was used to prove that each of
its elements has finite order. Also many other properties of the
Grigorchuk group and its analogs are proved using the contraction
properties, since it allows to prove statements by induction on the
length of group elements.

\paragraph{The iterated monodromy groups.} A rational function is said
to be \emph{sub-hyperbolic} (see~\cite{milnor}) if there exists an orbifold metric on a
neighborhood of its Julia set such that the
rational function is expanding with respect to this metric. We have
the following theorem (for a proof see~\cite{milnor}).

\begin{theorem}
  A rational function is sub-hyperbolic if and only if every orbit of
  a critical point is either finite or converges to an attracting
  finite cycle.
\end{theorem}

Thus, in particular, every postcritically finite rational function is sub-hyperbolic.

If $f\in\C(x)$ is a sub-hyperbolic rational function, then the group
$\IMG(f),$ with respect to any natural action, is contracting (see
Theorem~\ref{th:limjul}).


\paragraph{Other examples} include the Gupta-Sidki
group~\cite{gupta-sidki_group}, the Fabrykowski-Gupta
group~\cite{gufabr}, the group from the paper~\cite{bsv:jns} and many
others~\cite{bartholdi-g:parabolic}.

\paragraph{A non-contracting action} is, for example, the action of the
lamplighter group described in Subsection~\ref{ss:exss}. See also the
paper~\cite{dahmani} for an example of a weakly branch non-contracting
group.
\medskip

\begin{figure}[ht]
\psfrag{0,0}[bc]{$(0, 0)$} \psfrag{0,1}[cl]{$(0, 1)$}
\psfrag{1,0}[cr]{$(1, 0)$} \psfrag{1,1}[tc]{$(1, 1)$}
\psfrag{0}[bc]{$(0, 1)$} \psfrag{1}[bc]{$(1, 0)$}
\psfrag{nP}[rt]{$1$} \psfrag{nZ}[lb]{$0$} \psfrag{nM}[cc]{$-1$}
\begin{center}
\includegraphics{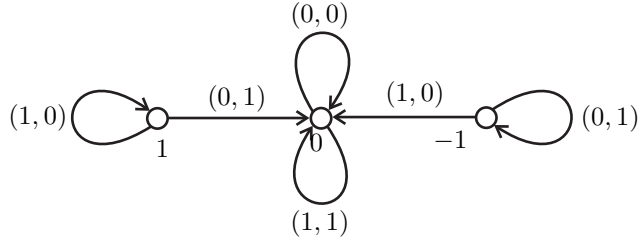}
\end{center}
\caption{The nucleus of the adding machine action} \label{fig:nuc}
\end{figure}

It follows from the definition that the restrictions of the elements
of the nucleus also belong to the nucleus. Thus, the nucleus is a
subautomaton of the complete automaton of the action. For instance,
the diagram of the nucleus of the adding machine is shown on
Figure~\ref{fig:nuc}.

For the Grigorchuk group the nucleus is the automaton defining the
generators (see Figure~\ref{fig:grig}, page~\pageref{fig:grig}).

An equivalent definition of contracting action uses the contraction of
the length of the group elements and is based on the following
proposition.

\begin{proposition}
  \label{pr:sprad}
  Let $G$ be a finitely generated group with a self-similar action.
  Let $|g|$ denote the word-length of $g\in G$ with respect to some
  fixed generating set of $G$. Then the number
  \[
  \rho=\limsup_{n\to\infty}\max_{v\in X^n}
  \sqrt[n]{\limsup_{|g|\to\infty}\frac{\left|g|_v\right|}{|g|}}
  \]
  is finite and does not depend on the choice of the generating system
  of the group.
  
  The number $\rho$ is less than $1$ if and only if the action of the
  group $G$ is contracting.
\end{proposition}

\subsection{Abelian self-similar groups}
\label{ss:abel}
Here we show, following~\cite{neksid}, some of the properties of
self-similar actions of the free abelian group $\Z^n$. We will use the
additive notation in this subsection.

Let $\phi=\phi_x:\Z^n\prr\Z^n$ be the virtual endomorphism associated
to a self-similar action over an alphabet of $d$ letters.

Let $G_1=\dom\phi$ be its domain. The map $\phi$ extends uniquely to a
linear map $\phi:\Q^n\to\Q^n$.  Since $\phi$ maps a subgroup of index
$d$ of $\Z^n$ to $\Z^n$, its matrix $A$ in a basis of $\Z^n$ has the
form
\[
\left(
  \begin{array}{cccc}
   \frac{a_{11}}{k_1} & \frac{a_{12}}{k_2} & \ldots & \frac{a_{1n}}{k_n}\\
   \frac{a_{21}}{k_1} & \frac{a_{22}}{k_2} & \ldots & \frac{a_{2n}}{k_n}\\
   \vdots             & \vdots             & \ddots & \vdots            \\
   \frac{a_{n1}}{k_1} & \frac{a_{n2}}{k_2} & \ldots & \frac{a_{nn}}{k_n}
  \end{array}
\right),
\]
where $a_{ij}$, $k_j$ are integers and $d=k_1k_2\cdots
k_n$. In particular, $d\cdot\det A\in\Z$ and the action is recurrent
if and only if $\det A=\pm d^{-1}$.

Every $\phi$-invariant subgroup must be trivial; this is equivalent to
the condition that the characteristic polynomial of the matrix $A$ is
not divisible by a monic polynomial with integral coefficients
(see~\cite{neksid} Proposition 4.1). In particular, the matrix $A$ is
non-degenerate.

If the action is recurrent then the matrix $A^{-1}$ has integral
entries. Any such action coincides with a natural action of $\Z^n$ as
the i.m.g.\ of the self-covering of the torus $\mathbb{T}^n$
defined by the matrix $A^{-1}$.

The following theorem is proved in~\cite{neksid} for the case
$d=2$. The proof in the general case is similar.

\begin{theorem}
\label{th:finab} A self-similar action of the group
$\Z^n$ is finite-state if and only if the associated
virtual endomorphism $\phi$ is a contraction, i.e., if its matrix
$A$ has spectral radius less than one.
\end{theorem}

Let us fix a recurrent self-similar action of the group $\Z^n$ over an
alphabet $X=\{0,1, \ldots, d-1\}$, with associated virtual
endomorphism $\phi$, and let $R=\{r_0=0,r_1, \ldots, r_{d-1}\}$ be a
digit set of the action.  Consider the closure $\hat\Z^n$ of the group
$\Z^n$ in $\aut T(X)$. Then $\hat\Z^n$ is an
abelian compact topological group. We have the following properties
(see~\cite{neksid}):
\begin{theorem}
  \label{th:closure}
  For every sequence $w=i_0i_1\ldots\in\xo$ the series
  \[\Phi(w)=r_{i_0}+\phi^{-1}(r_{i_1})+\phi^{-2}(r_{i_2})+\cdots+\phi^{-k}(r_{i_k})+\cdots\]
  converges in $\hat\Z^n$.
  
  The map $\Phi:\xo\to\hat\Z^n$ is a continuous bijection, which
  agrees with the action of the group $\Z^n$ on its closure, i.e.
  \[\Phi\left(w^g\right)=\Phi(w)+g\]
  for all $g\in\Z^n, w\in\xo$.
  
  The group $\hat\Z^n$ is isomorphic to the profinite completion of
  the group $\Z^n$ with respect to the series of subgroups of finite
  index
  \[\dom\phi\ge\dom\phi^2\ge\dom\phi^3\ge\ldots.\]
\end{theorem}

So, recurrent self-similar actions of abelian groups give us
``$\phi$-adic'' numeration systems on these groups and define
naturally their profinite $\phi$-adic completions $\hat\Z_{\phi}^n$.
For more on such numeration systems
see~\cite{vince:digtile,gilbert:three}.

\section{Schreier graphs}
\subsection{Definitions}
\label{ss:schdef}
Let $G$ be a group generated by a finite set $S$. We suppose
$1\notin S$ and $S=S^{-1}$.  Suppose that $G$ acts faithfully on a set
$M$.

Let us define the \emph{labeled Schreier graph} $\G(G, S, M)$ of the
group $G$ acting on $M$. It is a labeled graph with set of vertices
$M$ and set of edges $M\times S$. The label of every edge $(x, s)$ is
$s$. We set $\be(x, s)=x$ and $\en(x, s)=x^s$.

It is obvious that the labeled Schreier graph uniquely defines the
action of the generating elements on the set $M$, so it also
determines uniquely the group $G$.

Sometimes we will consider just \emph{Schreier graphs}, i.e., the graphs
defined in the same way but without the labeling.

A \emph{simplicial} Schreier graph is the simplicial graph associated
to a Schreier graph.

Obviously, the orbits of $G$ are exactly the vertex sets of the
connected components of the Schreier graph $\G(G, S, M)$.

If $x\in M$, then by $\G(G, S, x)$ we denote the Schreier graph of the
action of $G$ on the $G$-orbit of $x$. Such Schreier graphs are called
\emph{orbit Schreier graphs}.

If the group $G$ acts transitively on the set $M$, then the Schreier
graph $\G(G, S, M)$ can be interpreted in a more classical way:

\begin{defi}
  Let $G=\langle S\rangle$ be a group with distinguished generating
  set $S$, and let $H<G$ be a subgroup. The corresponding
  \emph{Schreier graph} is the graph whose vertices are the right
  cosets $H\backslash G=\{Hg: g\in G\}$ and whose set of edges is
  $(H\backslash G)\times S$, with maps $\be(Hg, s)=Hg$ and $\en(Hg,
  s)=Hgs$.
\end{defi}

If $G$ acts transitively on the set $M$, then the Schreier graph
$\G(G, S, M)$ is isomorphic to the Schreier graph corresponding to the
stabilizer $\stab(m)$, for any point $m\in M$.

In the special case of $G$ acting on itself by right-multiplication,
the Schreier graph $\G(G,S,G)$ is called the \emph{Cayley graph} of
$G$; it is the Schreier graph corresponding to the trivial subgroup.

\subsection{Schreier graphs of groups acting on rooted trees}
Suppose $G$ acts on the rooted tree $T(X)$ by automorphisms (this
holds, for instance, if $G$ is a self-similar group).

Then the levels $X^n$ are invariant under the action of $G$. Let us
denote by $\G_n(G, S)$ the Schreier graph of the action of $G$ on the
$n$th level. Then the Schreier graph $\G(G, S, \xs)$ of the action on
$\xs$ is the disjoint union of the graphs $\G_n(G, S)$.

For every $n\geq 0$, let $\pi_n:\G_{n+1}(G, S)\to\G_n(G, S)$ be the
map, defined on the vertex sets, given by $\pi_n(x_1\ldots
x_nx_{n+1})=x_1\ldots x_n$. Then, since $G$ acts by automorphisms of
the rooted tree, the maps $\pi_n$ induce surjective morphisms between
the labeled graphs. In this way we get an inverse spectrum of finite
labeled graphs
\begin{equation}
  \label{eq:spofgr}
  \G_0(G, S)\leftarrow\G_1(G, S)\leftarrow\G_2(G, S)\leftarrow\cdots
\end{equation}

We therefore get the following simple description of the graph $\G(G,
S, \xo)$:
\begin{proposition}
  \label{pr:profinitegr}
  The labeled Schreier graph $\G(G, S, \xo)$ is the inverse limit of
  the sequence~\eqref{eq:spofgr}.
\end{proposition}

The graphs which are isomorphic to an inverse limit of finite graphs
are called \emph{profinite graphs} (see~\cite{profinite_graphs}).

On the other side, it is possible to
interpret the orbit Schreier graphs on $\xo$ as limits of the finite
Schreier graphs $\G_n(G, S)$. The following proposition holds
(see~\cite{gr_zu:approx} for applications).

\begin{proposition}
  \label{pr:loclimit}
  Let $G$ be a finitely generated group acting on the tree $T(X)$ by
  automorphisms. Let $v=x_1x_2\ldots\in\xo$ be a point on the
  boundary. Then the pointed orbit Schreier graph $\left(\G(G, S, v),
    v\right)$ is isomorphic to the limit of the pointed Schreier
  graphs $(\G_n(G, S), x_1x_2\ldots x_n)$ with respect to the local
  topology on the space of pointed graphs.
\end{proposition}

Here the local topology on the space of pointed graphs is defined by the metric
$$d((\G_1, v_1), (\G_2, v_2))=2^{-R},$$ where $R$ is maximal among such that the ball
$B(v_1, R)$ in $\G_1$ and the ball $B(v_2, R)$ in $\G_2$ are isomorphic (respectively
to an isomorphism mapping $v_1$ to $v_2$).

\begin{defi}
  \label{defi:grcontr}
  A \emph{graph contraction} from a graph $(V, E)$ to a graph $(V',
  E')$ is a pair $f=(f_V,f_E)$ of maps $f_V:V\to V'$ and $f_E:E\to
  E'\cup\{\eth\}$ such that for edges $e\in E$ with $f_E(e)\neq\eth$,
  \[\be(f_E(e))=f_V(\be(e)),\quad\en(f_E(e))=f_V(\en(e));\]
  and $(f_E)_{|f_E^{-1}(E')}$ is a bijection from $f_E^{-1}(E')$ onto
  $E'$.
\end{defi}

In essence, some edges may be deleted (by sending them to $\eth$); all
other edges map bijectively onto the edges of $(V', E')$.

A graph $\Gamma=(V, E)$ is \emph{self-similar} if $V$ is infinite, and
there is a graph contraction $f:\Gamma\to\Gamma$ and a finite set
$N\subset V$ such that $\bigcup_{n\ge0}f^{-n}(N)=V$.

\begin{proposition}
  Let $G$ be a self-similar group generated by a finite set $S$ such
  that for every $g\in S$ and $x\in X$ the restriction $g|_x$ also
  belongs to $S$. Then for every $n\in\N$ the shift $\shift: xv\mapsto
  v:X^n\to X^{n-1}$ can be extended to a contraction of the Schreier
  graph $\G_n(G, S)$ onto the Schreier graph $\G_{n-1}(G, S)$.
\end{proposition}

\subsection{Examples of Schreier graphs of self-similar group actions}
\label{ss:exschreier}
\paragraph{The Grigorchuk group.} The Schreier graphs of the action of
the Grigorchuk group on the orbit of the point $111\ldots$, with its
self-similar nature, was described in~\cite{bgr:spec}.

\begin{figure}[ht]
  \begin{center}
    \psfrag{b}[br]{$b$} \psfrag{a}[tc]{$a$} \psfrag{c}[cr]{$c$}
    \psfrag{cl}[cl]{$c$} \psfrag{d}[tr]{$d$} \psfrag{db}[bc]{$d$}
    \psfrag{bc}[bc]{$b$} \psfrag{cc}[tc]{$c$}
    \includegraphics{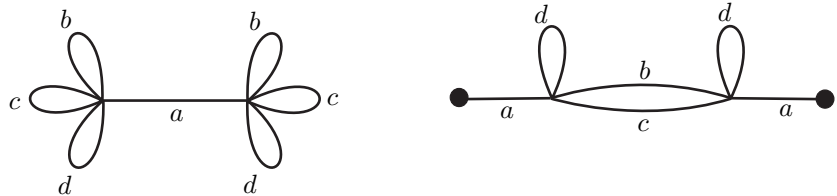}
  \end{center}
  \caption{Substitution rule for the Schreier graphs of the
    Grigorchuk group}
  \label{fig:grigsub}
\end{figure}

\begin{figure}[ht]
  \begin{center}
    \psfrag{b}[br]{$b$} \psfrag{a}[tc]{$a$} \psfrag{c}[cr]{$c$}
    \psfrag{cl}[cl]{$c$} \psfrag{d}[tr]{$d$} \psfrag{db}[bc]{$d$}
    \psfrag{bc}[bc]{$b$} \psfrag{cc}[tc]{$c$} \psfrag{bt}[tc]{$b$}
    \psfrag{cb}[bc]{$c$}
    \includegraphics{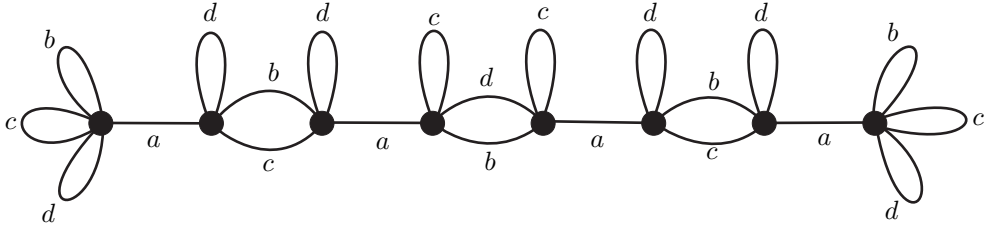}
  \end{center}
  \caption{A Schreier graph of the Grigorchuk group}
  \label{fig:gris3}
\end{figure}

The Schreier graph $\G_1(\mathbf{G}, S)$ of the action of the Grigorchuk
group is shown on left-hand side part of Figure~\ref{fig:grigsub}. In order to
obtain the Schreier graph $\G_n(\mathbf{G}, S)$, one has to replace in the
graph $\G_{n-1}(\mathbf{G}, S)$ simultaneously all the labels $b$ by the
labels $d$, the labels $c$ by $b$, the labels $d$ by $c$ and all the
edges labeled by $a$ by the graph depicted on the right part of
Figure~\ref{fig:grigsub} (the ends of the original edge correspond to
the marked vertices of the graph). See, for example, the Schreier
graph $\G_3$ in Figure~\ref{fig:gris3}. (We do not indicate the orientation of the arrows on the figures, since
the generators are involutions.)

\begin{figure}[ht]
  \begin{center}
    \psfrag{d}[cc]{$\cdots$} \psfrag{a}[cc]{(a)} \psfrag{b}[cc]{(b)}
    \includegraphics{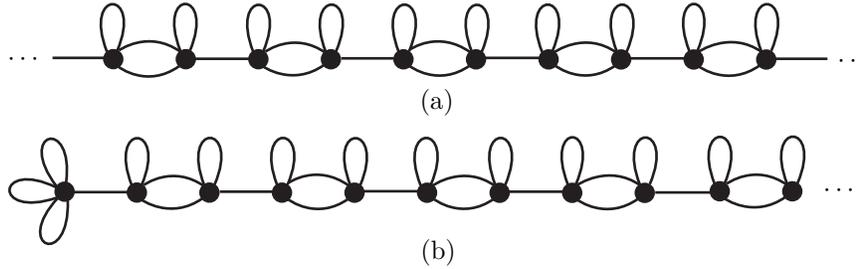}
  \end{center}
  \caption{The orbit Schreier graphs of the Grigorchuk group on the
    space $\xo$}
  \label{fig:grischrline}
\end{figure}

A component of the Schreier graph of the action of the Grigorchuk
group on the space $\xo$ is either the infinite line shown in
Figure~\ref{fig:grischrline}~(a), or the infinite ray shown in
Figure~\ref{fig:grischrline}~(b) (this last case occurs only for the
orbit of the point $1111\ldots$).

\paragraph{$\mathbf{\IMG(x^2-1)}$.} The iterated monodromy group of
the polynomial $x^2-1$ is generated by two generators
$a=(b, 1)\sigma$ and $b=(a, 1)$, where $\sigma$ is the
transposition.

\begin{figure}[ht]
  \begin{center}
    \unitlength=1pt\begin{picture}(0,0)
      \put(15,131){{\small $a$}}
      \put(35,131){{\small $b$}}
      \put(57,110){{\small $a$}}
      \put(57,142){{\small $a$}}
      \put(94,110){{\small $a$}}
      \put(94,142){{\small $a$}}
      \put(121,131){{\small $b$}}
      \put(136,131){{\small $a$}}
      \put(36,103){$\G_3$}
      \put(223,103){$\G_4$}
      \put(45,16){$\G_5$}
      \put(220,16){$\G_6$}
    \end{picture}
    \includegraphics{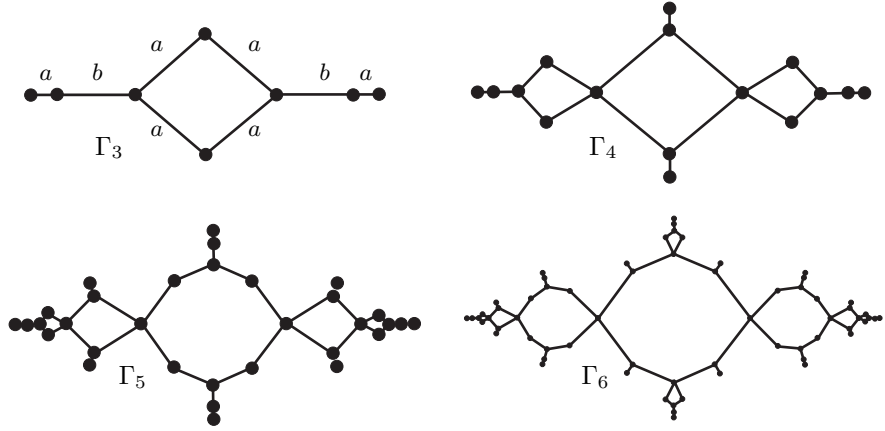}
\end{center}
  \caption{The Schreier graphs $\G_n(\img{z^2-1},\{a^{\pm1},b^{\pm1}\})$
    of the $n$th level action, for $3\le n\le 6$}
  \label{fig:g2shr}
\end{figure}

Some of the (simplicial) Schreier graphs of the group $\img{z^2-1}$ acting on
the finite levels are shown on Figure~\ref{fig:g2shr}. Compare with
the Julia set in Figure~\ref{fig:jul1}, page~\pageref{fig:jul1}.

The Schreier graphs of the groups $\img{z^2-1}$ are unions of
$2^n$-gons.

\begin{figure}[ht]
  \begin{center}
    \includegraphics[height=100mm]{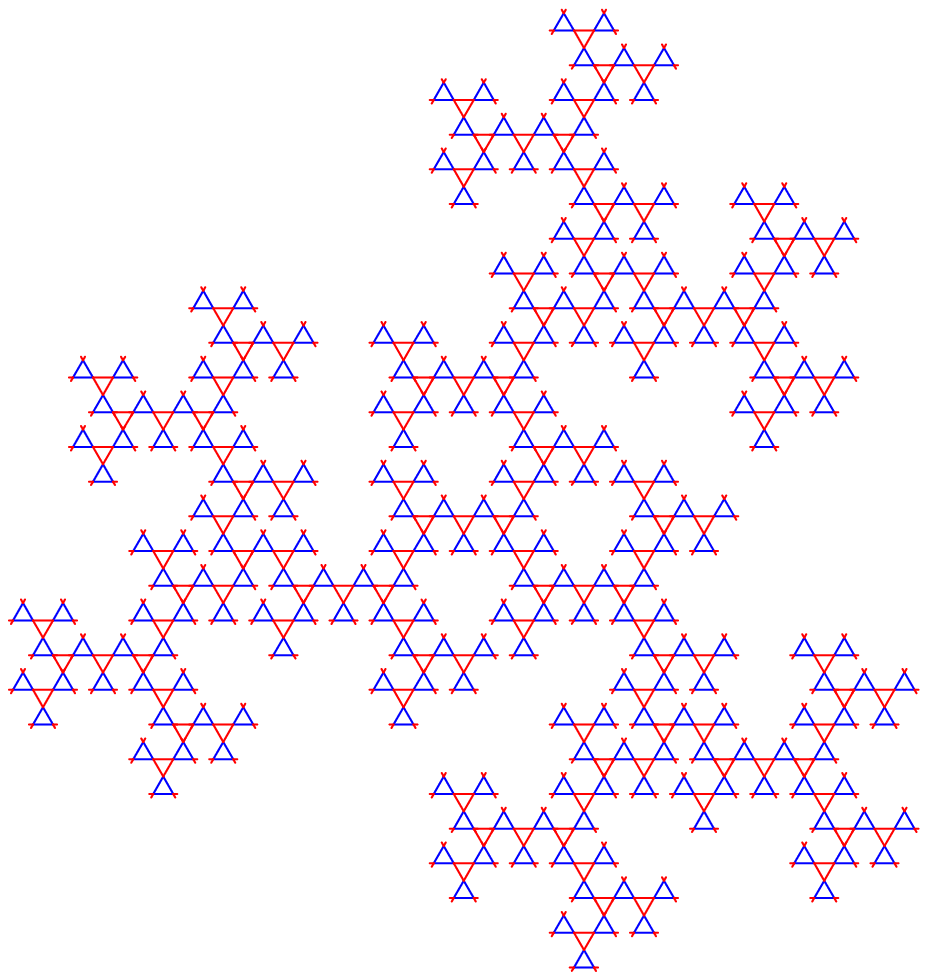}
    \caption{The Schreier graph $\G_6(G,S)$ of the Fabrykowski-Gupta group}
    \label{fig:fgshr}
  \end{center}
\end{figure}

\paragraph{The Fabrykowski-Gupta group.} The Schreier graph of this
group, introduced in Subsection~\ref{ss:exlimit}, is planar and is a
union of triangles.  The finite Schreier graph $\G_6(G,S)$ is given in
Figure~\ref{fig:fgshr}. As can be seen, the limit space and the
Schreier graph have a similar aspect.

\paragraph{Penrose tilings.} If we take the group $F$ generated by the
transformations $L$, $M$ and $S$ defined by the formul\ae\ from
Theorem~\ref{th:penr}, then it will act on the space $\mathcal{P}$,
with Schreier graphs isomorphic to the dual graphs of the Penrose
tilings (i.e., to the graphs whose vertices are tiles of the tiling,
with two vertices connected by an edge if and only if the respective
tiles have a common side), except for the Penrose tilings having non-trivial
symmetry. In that case the corresponding Schreier graph will be
isomorphic to the adjacency graph of the fundamental domain of the
symmetry group of the tiling, with loops at the vertices bounding the
domain.

\section{Growth and languages}
\label{s:growth}
In its most general form, the problem we deal with here is the
association to a geometric or combinatorial object of a numeric
invariant, the \emph{degree} or \emph{rate} of growth, or of a string
of numeric invariants, the \emph{growth power series}.  We sketch in
this section the main notions of growth, and present them in a unified
way.

The geometric objects described in this paper are of two natures: some
are compact ($X^\omega$, or the closure of $G$ in $\aut T(X)$), while
some are discrete ($G$, its Cayley graph, Schreier graphs, etc.)

Some other, more algebraic notions of growth or dimension may also be
integrated to this picture.  To name the main ones, growth of monoids
and automata (that are intimately connected to growth of groups);
cogrowth of groups (related to spectral properties of groups --- see
Section~\ref{s:spectrum}); subgroup growth~\cite{lubotzky:sg}; growth
of number of irreducible representations~\cite{passman-t:reps}; growth
of planar algebras~\cite{jones:planar}; growth of the lower central
series~\cite{grigorchuk:hilb_en,bartholdi-g:lie,petrogradsky:polynilpotent},
etc.

\subsection{Compact spaces}
\label{ss:grcompact}
Let $K$ be a compact metric space. Its \emph{Hausdorff dimension} (see~\cite{falconer:tech,MR2002e:28017})
is defined as follows: for $\beta>0$, the \emph{$\beta$-volume} of $K$
is
\[
H^\beta(K) = \lim_{\epsilon\searrow0}\inf_{\substack{\text{covers
      $\{\mathcal U_i\}$ of $K$ with}\\ \text{diameter at most
    $\epsilon$}}}\sum\operatorname{diam}(\mathcal U_i)^\beta.
\]
Clearly $H^\beta(K)$ is a decreasing function of $\beta$. The
\emph{Hausdorff dimension} $\dim_H(K)$ of $K$ is defined as the unique
value in $[0,\infty]$ such that $H^\beta(K)=\infty$ if
$0<\beta<\dim_H(K)$ and $H^\beta(K)=0$ if $\beta>\dim_H(K)$.

A connected, but easier-to-grasp notion, is that of \emph{box
  dimension}. It is defined, when it exists, as
\[\dim_\square(K) = -\lim_{\epsilon\searrow0}\frac{\ln(\text{number of
    $\epsilon$-balls needed to cover $K$})}{\ln\epsilon}.\] If
$\dim_\square(K)$ exists, then $\dim_H(K)$ exists too and takes the
same value.

For arbitrary topological spaces $F$, the following notion, which does
not refer to any metric, has been introduced: the \emph{topological
  dimension}, also called \emph{(Lebesgue) covering dimension}
$\dim_T(F)$ of $F$ is the minimal $n\in\N$ such that any open cover of
$F$ admits an open refinement of order $n+1$, i.e.\ such that no point
of $F$ is covered by more than $n+1$ open sets.

\subsection{Discrete spaces}
\label{ss:discretegrowth}
Let $\Gamma$ be a connected, locally finite graph, viewed as a
discrete metric space by assigning length $1$ to each edge. Choose a
base vertex $v\in V$. Then the growth of $\Gamma$ at $v$ is the
integer-valued function $\gamma_{\Gamma,v}:n\mapsto|B(v,n)|$ measuring
the volume growth of balls at $v$.

We introduce a preorder on positive-real-valued functions: say
$\gamma\precsim\delta$ if there is an $N\in\N$ such that
$\gamma(n)\le\delta(n+N)$ for all $n\in\N$; and say $\gamma\sim\delta$
if $\gamma\precsim\delta$ and $\delta\precsim\gamma$.

Clearly $\gamma_{\Gamma,v}(n)\le\gamma_{\Gamma,w}(n+d(v,w))$, so the
$\sim$-equivalence class of $\gamma_{\Gamma,v}$ does not depend on
$v$; we call it the \emph{growth of $\Gamma$}, written $\gamma_\Gamma$.

Note that if $\Gamma$ has degree bounded by a constant $D$, then
$\gamma_\Gamma\precsim D^n$.  The graph $\Gamma$ has \emph{polynomial
  growth} if $\gamma_\Gamma\precsim Kn^d$ for some $K,d\in\R$; the
infimal such $d$ is called the \emph{degree} of $\Gamma$. The graph
has \emph{exponential growth} if $\gamma_\Gamma\succsim b^n$ for some
$b>1$; the supremum of such $b$'s is called the \emph{growth rate} of
$\Gamma$.  In all other cases, $\Gamma$ has \emph{intermediate
  growth}.

The (polynomial) degree of growth is an exact analogue of the box
dimension defined above. Indeed, given a graph $\Gamma$ and a vertex
$v$, consider the metric spaces $K_n=\frac1n B(v,n)$, namely the balls
of radius $n$ with the metric scaled down by a factor of $n$. Then
each $K_n$ is compact (of diameter $1$). Assume $\Gamma$ has growth
degree $d$. Take the limit $K$ of a convergent subsequence (in the
Gromov-Hausdorff metric~\cite{gro:gr}) of $(K_n)_{n\ge1}$. Then
$\dim_\square(K)=d$.

Conversely, let $K$ be a compact space of box dimension $d$, with a
fixed point $*$. For $\epsilon=1/n$ cover $K$ by a minimal number of
$\epsilon$-balls, and consider the graph $\Gamma_n$, with vertex set
the set of balls, and edges connecting adjacent balls. Take the limit
$\Gamma$ of a convergent subsequence of $(\Gamma_n)_{n\ge0}$ (in the
local topology), with each $\Gamma_n$ based at the ball containing
$*$. Then $\Gamma$ is a graph of growth degree $d$.

We shall see in Subsection~\ref{ss:sgrcontr} examples of
Schreier graphs of polynomial growth, with associated compact spaces
of finite box dimension.

\subsection{Amenability}
\label{ss:amenab}
\begin{defi}
  Let $G$ act on a set $X$. The action is \emph{amenable} (in the
  sense of von Neumann~\cite{vneumann:masses}) if there exists a
  finitely additive measure $\mu$ on $X$, invariant under the action
  of $G$, with $\mu(X)=1$.
\end{defi}
We then say a group is amenable if its left- (or right-)
multiplication action on itself is amenable.

Amenability can be tested using the following criterion, due to F\o
lner for the regular action~\cite{folner:banach} (see
also~\cite{ceccherini-:amen} and the literature cited there):
\begin{theorem}
  Assume the group $G$ acts on a discrete set $X$. Then the action is
  amenable if and only for every if for every $\lambda>0$ and every
  $g\in G$ there exists a finite subset $F\subset X$ such that $|F\triangle
  gF|<\lambda|F|$, where $\triangle$ denotes symmetric difference and
  $|\cdot|$ cardinality.
\end{theorem}

\begin{corollary}
  Let $G$ act on a discrete space $X$, and assume the growth of $X$ is
  subexponential. Then the action is amenable.
\end{corollary}

Many other characterizations of amenability were discovered ---
see the reference~\cite{ceccherini-:amen}. The following are
equivalent:
\begin{enumerate}
\item The discrete $G$-space $X$ is non-amenable;
\item $X$ admits a \emph{paradoxical decomposition}, i.e.\ a partition
  $X=X_1\cup X_2$ with $X_1\cong X\cong X_2$, the $\cong$ sign
  indicating there is a piecewise-translational bijection between the
  spaces;
\item There exists a piecewise-translational map $X\to X$ with
  cardinality-$2$ fibers;
\item For any generating set $S$ of $G$, the
  simple random walk on the graph of the action of $G$ on $X$ has
  spectral radius strictly less than $1$.
\end{enumerate}

\subsection{Languages}
Most of the properties of growth of discrete spaces can be expressed
in terms of growth of languages; given a rooted locally finite graph
$(\Gamma,v)$, this is done by labeling the edges of $\Gamma$ by some
alphabet $X$, and identifying each vertex $w$ with some shortest path
from $v$ to $w$, usually the lexicographically minimal such path. In
case $\Gamma$ is a Cayley graph of a group $G$, whose edges are then
naturally labelled by generators, such a choice of paths is a
\emph{geodesic normal form}.

Let then $X$ be an alphabet. A \emph{language} is a subset $\lang$ of
$X^*$. Its \emph{growth} is the function
$\gamma_\lang:n\mapsto|\{w\in\lang:|w|\le n\}|$. Polynomial,
intermediate and exponential growth are defined as in
Subsection~\ref{ss:discretegrowth}. The language's \emph{growth
  series} is the formal power series
\[\Phi_\lang(t)=\sum_{w\in\lang}t^{|w|}=\sum_{n\ge0}\gamma_\lang(n)t^n.\]
As an excellent reference on more general power series
consult~\cite{salomaa-s:fps}.

In addition to their asymptotic behaviour, languages are classified by
their \emph{complexity} in the Chomsky hierarchy (see for
instance~\cite{hopcroft-u:automata} as a good reference).  The
simplest languages in that hierarchy are the \emph{regular languages},
and it is the smallest class of languages containing $\{x\}$ for all
$x\in X$, and closed under the operations of union, intersection,
concatenation, and iteration; this last one is defined as
$\lang^*=\{w^n:\,w\in\lang,n\ge0\}$, and actually intersection is not
necessary; if one replaces union by ``$+$'' and concatenation by
``$\cdot$'', one may write $\lang$ as a word over the alphabet $X\cup
\{(,),+,\cdot,^*\}$, called a \emph{regular expression}.

Let $N$ be a set disjoint from $X$; its elements are called
\emph{nonterminals}.  A \emph{grammar} is a collection of rules of the
form $v\to w$ for some $v,w\in(X\cup N)^*$; we write $v\models w$ if
$v=\alpha v_0\beta$, $w=\alpha w_0\beta$, and $v\to w$ for some
$\alpha\in X^*$ and $\beta,v_0,w_0\in(X\cup N)^*$. Given an initial
$n\in N$, the grammar \emph{produces} the language $\{w\in
X^*:n\models\dots\models w\}$. A language is regular if and only if it
is produced by a grammar with rules of the form $N\to X^*$ and $N\to
X^*N$.

Equivalently, a language $\lang$ is regular if and only if there
exists a finite state automaton $\autom$ (see
Subsection~\ref{ss:automata}), an ``initial'' state $q\in Q$ and a
set $F\subset Q$ of ``accepting'' states, such that $\lang$ is the set
of words in $X^*$ for which $\autom$, if started in state $q$,
ends in $F$. (The output of the automaton is discarded.)

A language is \emph{context-free} if it is produced by a grammar with
rules of the form $N\to(X\cup N)^*$; or equivalently if it is
recognized by a \emph{push-down automaton}, i.e., a finite-state
device that has access to the top of a linear stack.

A language is \emph{indexed} if it is recognized by a finite-state
device that has access to a tree-shaped stack; or equivalently if it
is produced by an \emph{indexed grammar}.

A language is \emph{context-sensitive} if it is produced by a grammar,
with no restrictions on the rules.

In addition to this hierarchy
``regular$\subset$context-free$\subset$indexed$\subset$context-sensitive'',
a language can be \emph{unambiguously} in one of these classes if
there exists such a grammar producing it, and such that or each
$w\in\lang$ the derivation $n\models\dots\models w$ is unique.
``Unambiguously regular'' is equivalent to ``regular'', but this does
not hold for the other classes.

It turns out that the formal power series captures essential
properties of $\lang$. We sum up the main facts:
\begin{proposition}\begin{enumerate}
  \item If the language $\lang$ is regular, then $\Phi_\lang$ is a
    rational function.
  \item If $\lang$ is unambiguously context-free, then $\Phi_\lang$ is
    an algebraic function.
  \item If $\lang$ is ambiguously context-free, then $\Phi_\lang$ may
    be transcendental, but its growth is either polynomial or
    exponential.
  \item If $\lang$ is indexed, then $\lang$ may have intermediate
    growth.
  \item If $\lang$ has intermediate growth, then $\Phi_\lang$ is
    transcendental.
  \end{enumerate}
\end{proposition}
\begin{proof}
  For the first point: If $\lang$ is regular, then it may be expressed
  unambiguously as a regular expression. In that expression, replace
  each variable by $t$ and $e^*$ by $(1-e)^{-1}$ for any
  expression $e$; the result is a rational function in $t$
  expressing the growth series of $\lang$.
  
  \noindent The second point, due to
  N.~Chomsky and M.-P.~Sch\"utzenberger~\cite{chomsky-s:alg}, is
  proved along similar lines.
  
  \noindent The third point follows, in its first part,
  from~\cite{flajolet:ambiguity} and, in its second part, from the
  independent work of~\cite{bridson-g:cflanguages}
  and~\cite{incitti:contextfree}.
  
  \noindent The fourth point follows from an example supplied
  in~\cite{grigorchuk-m:interm}.

  \noindent The fifth point follows from a result by
  P.~Fatou~\cite[page~368]{fatou:series}.
\end{proof}

\subsection{Groups}
\label{ss:grgroups}
In this section we sum up results on growth of groups, i.e.\ the
growth of $G$ seen as a metric space for the word metric. As a
reference, see~\cite[Chapters~VI--VIII]{harpe}
and~\cite{grigorchuk-h:groups}.

Fix a symmetric monoid-generating set $S$ for $G$. It induces a
\emph{word metric}
\[d(g,h)=\min\{n\in\N| g=s_1\dots s_nh,\text{ with }s_i\in S\},\]
on $G$, turning it into a (discrete) metric space. This metric is the
natural metric on the Cayley graph of $G$ introduced in
Subsection~\ref{ss:schdef}. If $S$ is finite, then $G$ is a locally
compact discrete space, and hence has a growth function
$\gamma_{G,S}(n)=|B(1,n)|$, and associated power series
$\Phi_{G,S}=\sum_{g\in G}t^{d(1, g)}$. Often we write $|g|=d(1,g)$, and
inversely $d(g,h)=|gh^{-1}|$. Note that we need $S$ to be symmetric
for $d$ to be a distance, but that last fact is not in itself
necessary.

An even stronger algebraic invariant of $(G,S)$ is the \emph{complete
  growth series}, introduced by F.~Liardet~\cite{liardet:phd}. It is
\[\tilde\Phi_{G,S}=\sum_{g\in G}gt^{|g|}\in\Z G[[t]],\]
i.e., a power series over the group algebra. A power series in $\Z
G[[t]]$ is rational if it lies in the closure of the polynomials ring
$\Z G[t]$ under the operations of addition, subtraction,
multiplication and \emph{quasi-inversion} $f^*=\sum_{n\ge0}f^n$ of a
series with no constant term.

There are many good reasons to consider properties of the growth
series of $G$; it is connected to the Hilbert-Poincar\'e series of the
classifying space $K(G,1)$, and in particular one often (but not
always) has $\Phi_G(1)=1/\chi(G)$, the latter being the Euler
characteristic~\cite{floyd-p:fuchsian,smythe:euler,MR97d:20031a}.
Other results are listed in the following
\begin{theorem}\begin{enumerate}
  \item If $G$ either has a finite-index abelian subgroup, or is
    hyperbolic in the sense of Gromov, or is a Coxeter group, then its    
    growth series is rational.
  \item There are soluble groups with algebraic, but not rational
    growth series.
  \item There are nilpotent groups $G$ that have rational growth for
    some generating sets and transcendental growth for others.
  \end{enumerate}
\end{theorem}
\begin{proof}
  The first point is due to M.~Benson~\cite{benson:zn} (``virtually
  abelian''), M.~Gromov~\cite[Chapitre~9]{ghys-h:gromov}
  (``hyperbolic''), and N.~Bourbaki~\cite[``exercice~26'']{bour:lie}
  (``Coxeter'').  Actually in all these cases the complete growth
  series is also rational.  F.~Liardet extended the virtually abelian
  case to the complete growth series, R.~Grigorchuk and
  T.~Smirnova-Nagnibeda computed explicitly the series for orientable
  surface groups~\cite{gri_tatiana}, and L.~Paris (unpublished)
  extended the result to Coxeter groups.

  \noindent The second point is due to W.~Parry~\cite{parry:wreath};
  explicitly, he computes the growth of some wreath products, and in
  particular of the wreath product $F_2\wr \Z/2\Z$ of a free group of
  rank $2$ with a cyclic group of order $2$.

  \noindent The last point is due to M.~Shapiro~\cite{shapiro:n2}
  (``rational'') and M.~Stoll~\cite{stoll:heisenberg}; for
  generalizations to many central extensions see~\cite{shapiro:psl2r}.
  
  The first author (unpublished) showed that although $\Phi$ is
  rational for some nilpotent groups, the complete growth series
  $\tilde\Phi$ is always transcendental. He also obtained examples of
  soluble groups for which $\Phi$ is rational and $\tilde\Phi$ is
  algebraic but non-rational.
\end{proof}

The growth function $\gamma_{G,S}$ depends heavily on the choice of
$S$. Write $\gamma\precsim\delta$ for two functions
$\gamma,\delta:\N\to\N$ if there exists a constant $C\in\N$ with
$\gamma(n)\le\delta(Cn)$ for all $n\in\N$, and $\gamma\sim\delta$ if
$\gamma\precsim\delta\precsim\gamma$. (Note that this equivalence
relation is coarser than the one introduced in
Subsection~\ref{ss:discretegrowth}.) Then the $\sim$-equivalence class
of $\gamma_{G,S}$ depends only on $G$.

\begin{theorem}\label{th:groupgrowth}\begin{enumerate}
  \item A group has polynomial growth if and only if it has a
    nilpotent finite-index subgroup; this implies that the growth
    degree is an integer.
  \item A soluble or linear group has either polynomial or exponential
    growth.
  \item There exist uncountable chains and antichains of groups of
    intermediate growth.
  \end{enumerate}
\end{theorem}
\begin{proof}
  The first point is due to Y.~Guivarc'h and independently
  H.~Bass~\cite{guivarch:poly1,bass:nilpotent} (``if'') and
  M.~Gromov~\cite{gro:gr} (``only if'').

  \noindent The second point is due to J.~Milnor and
  J.~Wolf~\cite{milnor:solvable,wolf:solvable} (``soluble'') and
  J.~Tits~\cite{tits:linear} (``linear'').
  
  \noindent The last point is due to the second
  author~\cite{grigorchuk:milnor_en,grigorchuk:kyoto}, who constructed the first
  example of a group $G$ of intermediate growth.
  See~\cite{bartholdi:upperbd,bartholdi:lowerbd} for the best-known
  estimates
  \[e^{n^{0.5157}}\precsim\gamma_G(n)\precsim e^{n^{0.7675}}.\]
\end{proof}

The following notions have also been introduced
(see~\cite{grigorchuk-h:groups,bozejko:uniform}): a group $G$ has
\emph{uniformly exponential growth} if there is a $b>1$ such that
$\gamma_{G,S}(n)\ge b^n$ for all $S$; the main point is that $b$ does
not depend on $S$.  Soluble~\cite{osin:solvable}, non-elementary
Gromov-hyperbolic~\cite{koubi:unifexpo},
one-relator~\cite{grigorchuk-h:1rel}, and most amalgamated products
and $HNN$ extensions~\cite{bucher-h:hnn} are known to have uniformly
exponential growth; no example is known of a group with exponential
growth, but not uniformly so.

The group $G$ is \emph{growth tight} if
$\sqrt[n]{\gamma_{G/N,S/N}(n)}<\sqrt[n]{\gamma_{G,S}(n)}$
for any infinite normal subgroup $N\triangleleft G$.  Hyperbolic
groups have uniformly exponential growth, and are growth
tight~\cite{arzhantseva-l:growthtight}.

It is unknown whether all groups of exponential growth have uniformly
exponential growth; however, there are examples of groups that do not
reach the infimum of their growth rates~\cite{sambusetti:hnn}.

Suppose $G$ is residually a $p$-group, and let $\{G_n\}$ be
its lower $p$-central series. It is known~\cite{bartholdi-g:lie} that
if the rank of $G_n/G_{n+1}$ increases exponentially in $n$,
then $G$ has uniformly exponential growth.

Growth tightness can also be defined for languages: given a language
$\lang$ and a word $w\in\lang$, one defines
$\lang_w=\{v\in\lang:\,v\text{ does not contain }w\}$; then $\lang$ is
\emph{growth tight} if for all $w\in\lang$ the growth rate of
$\lang_w$ is strictly less than that of $\lang$.
T.~Ceccherini-Silberstein and W.~Woess have shown that ``ergodic''
context-free languages (i.e., such that for any $n_1,n_2\in N$ we have
$n_1\models\dots\models un_2v$ for some $u,v\in(T\cup N)^*$) are
growth tight.

\subsection{$L$-presentations}
Let $S$ be a finite alphabet, and let $\Phi$ be a finite set of monoid
endomorphisms: $S^*\to S^*$. In effect, each $\phi\in\Phi$ is
determined by the values $\phi(s)$ for all $s\in S$. Let $I\subset
S^*$ be a finite set of initial words.  A \emph{D0L
  system}~\cite{MR82g:68053} is the closure of $I$ under iterated
application of $\Phi$. This formalism, named after A.~Lindenmayer, was
introduced to study growth of biological organisms.

The following definition was introduced in~\cite{bartholdi:lpres}:
\begin{defi}
  A group $G$ has a \emph{finite $L$-presentation} if it can be
  presented as $F_S/N$, where $F_S$ is the free group on $S$ and $N$
  is the normal closure of a D0L system.
\end{defi}

The class of groups with a finite $L$-presentation clearly contains
finitely presented groups, and enjoys various closure properties, in
particular that of being closed under wreath products by finite
groups.

The main motivation in studying $L$-presentations is the following
result, also taken from~\cite{bartholdi:lpres}:
\begin{theorem}
  Let $G$ be a contracting regular branch group. Then $G$ has a
  finite $L$-presentation, but is not finitely presented.
\end{theorem}

This is a generalization to branch groups of an earlier result due to
I.~Lysionok~\cite{lysionok:pres}:
\begin{theorem}
  The Grigorchuk group $\mathbf{G}$ admits the following presentation:
  \[\mathbf{G} = \langle a,c,d|\sigma^i(a^2),\sigma^i(ad)^4,
  \sigma^i(adacac)^4\;(i\ge0)\rangle,\] where
  $\sigma:\{a,c,d\}^*\to\{a,c,d\}^*$ is defined by
  $\sigma(a)=aca,\sigma(c)=cd,\sigma(d)=c$.
\end{theorem}
Note incidentally that the Schreier graphs $\G_n(\mathbf{G},S)$ of $\mathbf{G}$
described in Subsection~\ref{ss:exschreier} are isomorphic to line
segments of length $2^n-1$; the labeling on the edge is $\sigma^n(a)$.
This is the tip of a connection between the lower central series of a
branch group and its Schreier graph structure~\cite{bartholdi:lcs}.

\subsection{Semigroups and automata}
The automata from Subsection~\ref{ss:automata} are not necessarily
invertible; the most general setting in which growth questions can be
studied is that of \emph{semigroups}. This subsection discusses growth
of semigroups, and in particular semigroups generated by automata.

Let $T$ be a semigroup generated by a finite set $S$. Analogously to
the group situation, the \emph{growth function} of $T$ is the function
$\gamma(n)=|\{t\in T| t=s_1\dots s_n,\,s_i\in S\}|$.

On the other hand, let $\autom$ be a finite automaton, and let
$\autom^n$ be the $n$-fold composition $\autom*\dots*\mathcal A$ of
$\autom$ defined in Subsection~\ref{ss:automata}. The \emph{growth
  function} of $\autom$ is the function $\gamma(n)=\text{minimal
  number of states of }\autom^n$.

The following connection is clear:
\begin{proposition}
  The growth function of $\autom$ is the growth function of the
  semigroup generated by $\{\autom_q|q\in Q\}$.
\end{proposition}

Let $T$ be a cancellative semigroup (i.e.\ a semigroup satisfying the
axiom $(xz=yz\vee zx=zy)\Rightarrow x=y$), and let $G$ be its group
of (left) quotients. The second author studied in~\cite{MR92k:20116} the
connections between amenability, and growth, of $T$ and $G$. He showed
in~\cite{MR89f:20065} that a cancellative semigroup has polynomial
growth if and only if its group of left quotients $G$ is virtually
nilpotent --- and in that case that the growth degrees of $G$ and $T$
are the same. This generalizes Gromov's statement in
Theorem~\ref{th:groupgrowth}.

Essentially, semigroups can have any growth function at least
quadratic and at most exponential; and it was long known that Milnor's
question (``Do all finitely generated groups have either polynomial or
exponential growth?'') has a negative answer in the context of
semigroups --- see the example by V.~Belyaev, N.~Sesekin and
V.~Trofimov~\cite{belyaev-:interm}.

Later many other examples of semigroups of intermediate growth were
discovered~\cite{shneerson:interm}. Consult also the
book~\cite{2001g:20076}, and~\cite{lavrik:interm}, where semigroups of
intermediate growth are found even among $2\times2$-matrices.

\subsection{Growth and dimension of self-similar groups}
For subgroups $G$ of $\aut T(X)$, we have the following extra results
on growth. First, $\aut T(X)$ is a compact metric space, for which
$\{\stab_n\}_{n\in\N}$ is a basis of neighbourhoods: define a metric
on $\aut T(X)$ by
\[d(g,h)=\max\{|\aut T(X)/\stab_n|:g^{-1}h\in\stab_n\}.\]
Then the \emph{Hausdorff dimension} of $G$ is the Hausdorff dimension
of its closure $\overline G$, with the restricted metric, in $\aut
T(X)$. Since in the given metric all balls
are cosets of $\stab_n$ for some $n$, we have the simpler
definition~\cite{barnea-s:hausdorff}
\[\dim_H(G)=\lim_{n\to\infty}\frac{\log|(G\stab_n)/\stab_n|}{\log|\aut T(X)/\stab_n|}.\]

For instance, consider again the Grigorchuk group $\mathbf{G}$. It is
known that $\mathbf{G}\cap \stab_1$ has index $2$ in $\mathbf{G}$, and
that $\mathbf{G}\cap \stab_1$ embeds in $\mathbf{G}\times\mathbf{G}$
with index $8$. A simple calculation gives
$|\mathbf{G}/(\mathbf{G}\cap \stab_3)|=2^7$, so by induction
$\mathbf{G}/(\mathbf{G}\cap \stab_n)=2^{\frac582^n+2}$ for $n\ge3$. On
the other hand, $|\aut T(X)/\stab_n|=2^{2^n-1}$. Therefore
$\dim_H(\mathbf{G})=\frac58$. This generalizes as follows:
\begin{theorem}[\cite{bartholdi:imgexamples}]
  Let $G$ be a regular branch group. Then its Hausdorff dimension
  $\dim_H(G)$ is a rational number in $(0,1)$.
\end{theorem}

All known groups of intermediate growth are subgroups of $\aut T(X)$
for some finite $X$, and are described by the following result:
\begin{theorem}[\cite{bartholdi:interm}]
  Let $G<\aut T(X)$ be generated by an automaton $X$. Assume the
  states of $T$ are partitioned as $\{1\}\cup T_1\cup T_2$, where $1$
  is an inactive state that maps to itself; all states in $T_1$ map to
  $1$; and all states in $T_2$ are inactive, and all their maps except
  one are to $T_1\cup\{1\}$. Then $G$ has subexponential growth.
\end{theorem}

\subsection{Schreier graphs of contracting actions and their growth}
\label{ss:sgrcontr}
Let $(G, \xo)$ be a self-similar action of a finitely generated group.
We assume that the orbits of the action are infinite. This is the case
for any level-transitive action of a group.

A graph $\G$ has polynomial growth if and only if the number
\[\alpha=\limsup_{r\to\infty}\frac{\log|B(v, n)|}{\log n}\]
is finite. The number $\alpha$ is then called the \emph{degree} of the
growth.

\begin{proposition}
  \label{pr:polgr}
  The growth of every orbit Schreier graph of a contracting action
  $(G, \xo)$ is polynomial.
  
  The growth degree is not greater than $-\frac{\log |X|}{\log\rho}$,
  where $\rho$ is the contraction coefficient of the action.
\end{proposition}

\paragraph{Examples.}
\begin{enumerate}
\item The orbit Schreier graph of the Grigorchuk group have linear
  growth (i.e., the growth degree is equal to $1$). This follows
  directly from their description.  The contraction coefficient of the
  Grigorchuk group is equal to $1/2$.
\item The orbit Schreier graph of the iterated monodromy group
  $\img(z^2-1)$ has polynomial growth of degree $2$. For example, the
  Schreier graph on the orbit of the point $111\ldots=1^\infty$ is
  described as follows:
  
  There is a $b$-labeled loop at $1^\infty$.  Between $0^{2n}1^\infty$
  and $0^{2n+1}1^\infty$ there is a $2^{n+1}$-gon labeled by
  $a$-edges, and between $0^{2n-1}1^\infty$ and $0^{2n}1^\infty$ there
  is a $2^n$-gon labeled by $b$-edges. The vertices on these polygons
  are labeled by strings in $(00|10)^n(0|1)1^\infty$ and
  $(00|01)^n1^\infty$ respectively.
  
  Then at each of these new vertices on the polygons finite graphs are
  attached; if the vertex is labeled $0^k1w1^\infty$ for some
  $w\in\{0,1\}^*$, then the attached graph has $2^k$ vertices labeled
  by all words in $(0|1)^k1w^\infty$.
  
  It therefore follows that the ball of radius $2^n-2$ at $1^\infty$
  contains only vertices with labels in $(0|1)^{2n-2}1^\infty$, and
  contains all vertices with labels in $(0|1)^{2n-4}1^\infty$. It
  follows that the cardinality of $B(1^\infty,2^n)$ is approximately
  $2^{2n}$.

  This graph is self-similar under the graph contraction $f=(f_V,f_E)$
  given by $f_V=\shift: s_1s_2\dots\mapsto s_2\dots$, and
  \[f_E:\begin{cases}
    (0w,a)\mapsto (w,b),\\
    (0w,b)\mapsto (w,a),\\
    (1w,a)\mapsto\eth,\\
    (1w,b)\mapsto\eth.
  \end{cases}\]
  It contracts distances by a factor of $\sqrt2$, while collapsing $2$
  points to $1$;
\end{enumerate}

\section{Limit spaces of self-similar group actions}
In this section we return from self-similar (semi)group actions to
self-similar topological spaces, showing that a naturally-defined
self-similar topological space is associated with every contracting
self-similar action. This space can be defined in different ways: as a
quotient of the Cantor set by an \emph{asymptotic equivalence
  relation} (Definition~\ref{defi:real_eq}), as a limit of finite
Schreier graphs (Theorem~\ref{th:limscrit}) or as the boundary of a
naturally-defined hyperbolic graph (Theorem~\ref{th:hypsp}).

\subsection{The limit space $\lims{G}$}
\label{ss:limsp}
Let us fix a self-similar contracting action of a group $G$ on the
space $\xo$.

Denote by $\xmo$ the space of all sequences infinite to the left:
\[\xmo=\{\ldots x_3x_2x_1\}\]
equipped with the product topology.

\begin{defi}
  \label{defi:real_eq} Two elements $\ldots x_3x_2x_1,
  \ldots y_3y_2y_1\in\xmo$ are said to be \emph{asymptotically
    equivalent} with respect to the action of the group $G$, if there
  exist a finite set $K\subset G$ and a sequence $g_k\in K, k\in\N$
  such that
  \[(x_kx_{k-1}\ldots x_2x_1)^{g_k}= y_ky_{k-1}\ldots y_2y_1\]
  for every $k\in\N$.
\end{defi}

It follows directly from the definition that the asymptotic
equivalence is an equivalence relation.

\begin{proposition}
  \label{pr:regequiv}
  Let $\nuke$ be the nucleus of the action. Then two sequences $\ldots
  x_2x_1, \ldots y_2y_1\in\xmo$ are asymptotically equivalent if and
  only if there exists a sequence $h_n\in \nuke, n\geq 0$ such that
  \begin{equation}\label{eq:aseq2}
    x_n^{h_n}=y_n, \quad\text{and}\quad h_n|_{x_n}=h_{n-1}
  \end{equation}
  for all $n\geq 1$.
\end{proposition}

Proposition~\ref{pr:regequiv} can be reformulated in the following
terms:
\begin{proposition}
  \label{pr:regequiv2}
  Let $\Gamma$ be the Moore diagram of the nucleus $\nuke$ of the
  action of the group $G$. Two sequences $\ldots x_2x_1, \ldots
  y_2y_1\in\xmo$ are asymptotically equivalent if and only if $\Gamma$
  has a path $\ldots e_2e_1$ such that every edge $e_i$ is labeled by
  the pair $(x_i, y_i)$.
\end{proposition}

\begin{defi}
  \label{defi:real}
  The \emph{limit space} of the self-similar action (written
  $\lims{G}$) is the quotient of the topological space $\xmo$ by the
  asymptotic equivalence relation.
\end{defi}

It follows from the definition of the asymptotic equivalence relation
that $\lims{G}$ is invariant under the shift map $\shift:\ldots
x_3x_2x_1\mapsto\ldots x_4x_3x_2$, and therefore the shift
$\shift:\xmo\to\xmo$ induces a surjective continuous map
$\si:\lims{G}\to\lims{G}$ on the limit space. Every point
$\xi\in\lims{G}$ has at most $|X|$ preimages under the map $\si$.

\begin{defi}
  The dynamical system $(\lims{G}, \si)$ is called the \emph{limit
    dynamical system} associated to the self-similar action.
\end{defi}

\begin{examp}
  In the case of the adding machine action of $\Z$, one sees from the
  diagram of the nucleus given in Figure~\ref{fig:nuc} that two
  sequences in $\xmo$ are asymptotically equivalent if and only if
  they are equal or are of the form
  \[
  \ldots 0001x_mx_{m-1}\ldots x_1,\text{ and } \ldots 1110x_mx_{m-1}\ldots x_1.
  \]
  This is the usual identification of dyadic expansions of reals in
  $[0,1]$
  \[
  0.x_1x_2\ldots x_m0111\ldots=0.x_1x_2\ldots x_m1000\ldots.
  \]
  
  Consequently, the limit space $\lims\Z$ is the circle that we obtain
  after identifying the endpoints of the unit segment (since the
  asymptotic equivalence relation glues the sequences $\ldots 000$ and
  $\ldots 111$). The map $\si$, induced on the circle from the shift
  on the space $\xmo$, is the two-fold covering map
  $\si(x)=2x\;(\mathop{\mathrm{mod}} 1)$.
\end{examp}

Proposition~\ref{pr:regequiv} implies now the following properties of
the limit spaces (see~\cite{nek:lim}).
\begin{theorem}
  \label{th:contrsp}
  Let $\nuke$ be the nucleus of a contracting self-similar action of a
  group $G$. The asymptotic equivalence is a closed equivalence
  relation on $\xmo$ and every point is equivalent to at most
  $|\nuke|$ points.

  The limit space $\lims{G}$ is metrizable and has topological
  dimension $\leq |\nuke|-1$.
  
  If the group is $G$ finitely generated then the limit space
  $\lims{G}$ is connected. 
 \end{theorem}

\subsection{Self-similarity of the space $\lims{G}$}
\label{ss:markov}
Here we construct a self-similarity structure on the space $\lims{G}$
which will be defined by a Markov partition of the dynamical system
$(\lims{G}, \si)$

\begin{defi}
  For a given finite word $v\in\xs$, define the \emph{tile} $\til_v$
  to be the image of the set $\xmo v=\{\ldots x_2x_1v\}$ under the
  canonical map $\xmo\to\lims{G}$.
\end{defi}

It follows from the definitions that $\til_{\emp}=\lims{G}$ and that
\[\si(\til_v)=\til_{v'},\]
where $v'$ is the word obtained from the word $v$ by deletion of
its last letter. We also have
\[\til_v=\bigcup_{x\in X}\til_{xv}.\]

Consequently, for every fixed $n$ the set of the tiles $\{\til_v: v\in
X^n\}$ is a Markov partition of the dynamical system $(\lims{G},
\si)$. Therefore, the tiles together with the restrictions of the maps
inverse to the shift $\si$ define a self-similarity structure on the
space $\lims{G}$.

Every tile is a compact set and any point of the limit space belongs
to not more than $|\nuke|$ tiles, where $\nuke$ is the nucleus of the
action.

\begin{proposition}
  \label{pr:tilgr}
  The tiles $\til_v, \til_u$, for $u, v\in X^n$, intersect if and only
  if there exists $h\in\nuke$ such that $v^h=u$.
\end{proposition}

\begin{defi}
  \label{defi:adjgr}
  Denote by $\mathrm{J}_n(G)$ the simplicial graph whose vertices are
  the tiles $\til_v$ for $v\in X^n$, with two vertices connected by an
  edge if and only if the respective tiles have a nonempty
  intersection.
\end{defi}

Then Proposition~\ref{pr:tilgr} can be formulated in the following
way:
\begin{corollary}
  \label{cor:tilgr}
  The map $v\to\til_v$ is an isomorphism between the simplicial
  Schreier graph $\G(\left\langle\nuke\right\rangle, \nuke, X^n)$ and
  the graph $\mathrm{J}_n(G)$.
\end{corollary}

If the action is recurrent, then the nucleus $\nuke$ generates the
group $G$, and thus the graphs $\mathrm{J}_n(G)$ are the Schreier graphs
of the group $G$.

The following theorem shows that the limit space $\lims{G}$ is a limit
of the graphs $\mathrm{J}_n(G)$:
\begin{theorem}
  \label{th:limscrit}
  A compact Hausdorff space $\X$ is homeomorphic to the
  limit space $\lims{G}$ if and only if there exists a collection
  $\mathfrak{T}=\{T_v: v\in\xs\}$ of closed subsets of $\X$
  such that the following conditions hold:
  \begin{enumerate}
  \item $T_\emp=\X$ and $T_v=\cup_{x\in X}T_{xv}$ for every
    $v\in\xs$;
  \item For every word $\ldots x_3x_2x_1\in\xmo$ the set
$\cap_{n=1}^\infty T_{x_nx_{n-1}\ldots x_1}$ contains only
    one point;
 \item The intersection $T_v\cap T_u$ for $u, v\in X^n$ is non-empty
    if and only if there exists an element $s$ of the nucleus such
    that $v^s=u$.
  \end{enumerate}
\end{theorem}
In particular, it follows that if for some sequence of numbers $R_n>0$
the metric spaces $(\mathrm{J}_n(G), d(u, v)/R_n)$ converge in the
Gromov-Hausdorff metric to a metric space $\mathfrak X$, then
$\mathfrak X$ is homeomorphic to the limit space $\lims{G}$.

If $\X$ is a metric space then condition~(2) can be replaced by the condition:
\[\lim_{n\to\infty}\max_{v\in X^n}\mathop{\mathrm{diam}}(T_v)=0.\]

\subsection{Examples of limit spaces}
\label{ss:exlimit}
\paragraph{The Grigorchuk group.} It follows from the Moore diagram of
the nucleus of the Grigorchuk group (see Figure~\ref{fig:grig}) that
two sequences are asymptotically equivalent if and only if they are
either equal or have the form
\[\xi = \ldots 111101w\quad \zeta = \ldots 111100w,\]
where $w\in\xs$.

Let us define a homeomorphism $F$ of the space $\xmo$:
\[F(\ldots x_3x_2x_1)=\ldots y_3y_2y_1,\]
with $y_i=(1+x_1)+(1+x_2)+\cdots+(1+x_i)\mod 2$.

Then two elements $\xi, \zeta\in\xmo$ are asymptotically equivalent
with respect to the action of the Grigorchuk group if and only if the
$F(\xi)$ and $F(\zeta)$ are asymptotically equivalent with respect to
the adding machine action, except in the case of the points $\ldots
1111$ and $\ldots 11110$. These two points are not equivalent, while
their $F$-images $\ldots 0000$ and $\ldots 1111$ are equivalent.
Therefore the limit space of the Grigorchuk group is homeomorphic to
the real segment $[0,1]$.

The shift $\si$ on the space $\lims{G}$ is the ``tent'' map
$\si:x\mapsto 1-|2x-1|$, folding the interval in two.

The fact that the limit space of the Grigorchuk group is the segment
also follows from Theorem~\ref{th:limscrit} and the description of the
Schreier graphs of its action on the sets $X^n$.

\begin{figure}[ht]
  \begin{center}
    \psfrag{21}[rb]{$(2, 2)$}
    \psfrag{12}[lb]{$(1, 1)$}
    \psfrag{0}[rb]{$(0, 0)$}
    \psfrag{s}[ct]{$\begin{array}{c} (0, 1)\\ (1, 2)\\ (2, 0)\end{array}$}
    \psfrag{e}[lc]{$\begin{array}{c}(0, 0)\\ (1, 1)\\ (2, 2)\end{array}$}
    \includegraphics[height=50mm]{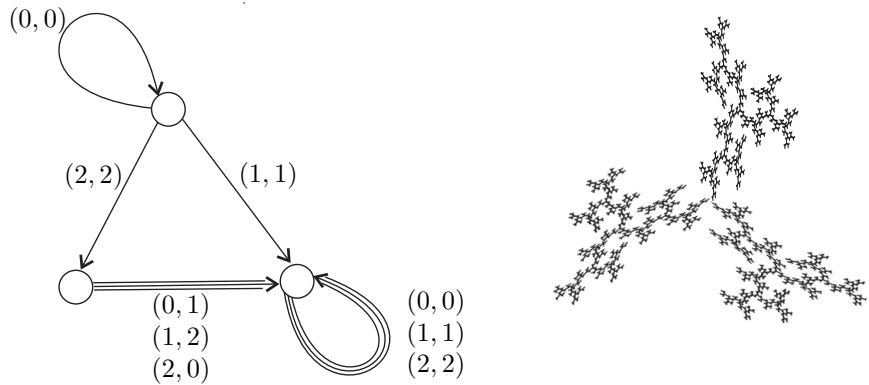}
  \end{center}
  \caption{The Fabrykowski-Gupta group and its limit space}
  \label{fig:bart}
\end{figure}

\paragraph{The Fabrykowski-Gupta group.} The limit space of this
group, defined by J.~Fabrykowski and N.~D.~Gupta~\cite{gufabr}, is the
dendrite fractal shown on Figure~\ref{fig:bart}. The picture is drawn
using the description of Schreier graphs of the group
(see~\cite{bgr:spec}) and Theorem~\ref{th:limscrit}. This dendrite
fractal is also homeomorphic to the fractal described on Figure~8.2 of
the book~\cite{falcon:geom}.

\paragraph{The Sierpi\'nski gasket.} Let the alphabet $X$ be
$\{0,1,2\}$. Define three transformations $b_i$, $i\in X$ of the space
$\xo$ as follows: put $(iw)^{b_i}=i\left(w^{b_i}\right)$ and
$(jw)^{b_i}=kw$ for all $i, j, k$ such that $\{i, j, k\}=\{0,1,2\}$.
In the standard notation:
\[b_0=(b_0,1,1)\sigma_{12},\quad b_1=(1,b_1,1)\sigma_{02},\quad
b_2=(1,1,b_2)\sigma_{01}.\]

One can show, using Theorem~\ref{th:limscrit}, that the limit space of
the group $\langle b_1, b_2, b\rangle$ is the Sierpinski gasket.

\paragraph{The Julia sets of rational mappings.} The \emph{Julia set}
of a rational function $f\in\C(z)$ can be defined as the closure of the
union of its repelling cycles (see~\cite{milnor,lyubich:top}). If $f$ is polynomial, its Julia set is
the boundary of the attraction basin of $\infty$.

The proof of the next theorem will appear in~\cite{nek:img}.
\begin{theorem}
  \label{th:limjul}
  The iterated monodromy group of a sub-hyperbolic rational function
  (respectively to any natural action on a regular tree) is
  contracting; and its limit space is homeomorphic to the Julia set of
  the rational function.
\end{theorem}

\paragraph{Examples.}
\begin{enumerate}
\item The Julia set of the polynomial $z^2$ is the circle. Its
  iterated monodromy group is the adding machine action, which is
  contracting and with the limit set homeomorphic to the circle.
\item The Julia set of the polynomial $z^2-2$ is the segment $[-2,2]$.
  The iterated monodromy group is the dihedral group
  $\mathbb{D}_\infty$ generated by the automaton in
  Figure~\ref{fig:dihedr}. The limit space of this action is
  homeomorphic to the segment, and this is proved in the same way as
  for the Grigorchuk group (in fact the asymptotic equivalence
  relations in these two cases are the same).
\item The Julia set of the polynomial $z^2-1$ is shown in
  Figure~\ref{fig:jul1}. Compare its shape with the pictures of the
  Schreier graphs of the iterated monodromy group of this
  polynomial given in Figure~\ref{fig:g2shr}.
  \begin{figure}[ht]
    \begin{center}
      \includegraphics[width=100mm]{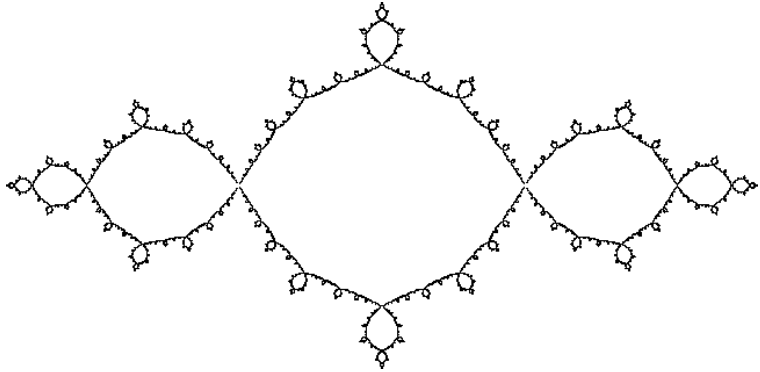}
    \end{center}
    \caption{The Julia set of the polynomial $z^2-1$.}
    \label{fig:jul1}
  \end{figure}
\item The Julia set of the polynomial $z^2+i$ is the dendrite shown in
  Figure~\ref{fig:i}.  Its tree-like structure parallels the fact that
  the Schreier graphs of its iterated monodromy group are
  all trees. This holds in fact for all Misiurewicz polynomials, i.e.,
  quadratic polynomials $z^2+c$ for which the critical point $0$ is
  strictly pre-periodic. See the paper~\cite{kameyama:julia}, where
  the self-similarity of Julia sets of Misiurewicz polynomials is
  studied.
  \begin{figure}[ht]
    \begin{center}
      \includegraphics[width=60mm]{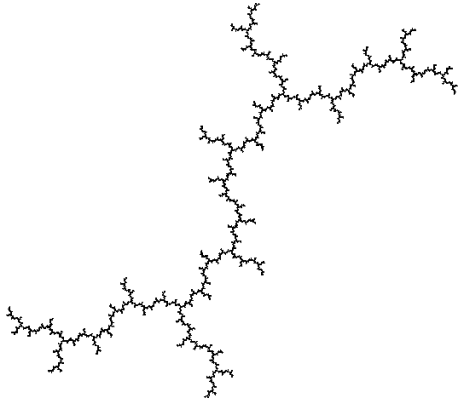}
    \end{center}
    \caption{The Julia set of the polynomial $z^2+i$.}
    \label{fig:i}
  \end{figure}
\item The Julia sets of the polynomials $z^2+c$, where
  $c^3+2c^2+c+1=0$, are shown in Figure~\ref{fig:jul2}. The left one,
  called the ``airplane'', corresponds to the real root $c$, while the
  right one, called ``Douady's rabbit'', corresponds to the complex
  root with positive imaginary part. Their corresponding groups have
  almost the same recursion --- see Subsection~\ref{ss:exmono}.
  \begin{figure}[ht]
    \begin{center}
      \includegraphics[width=60mm]{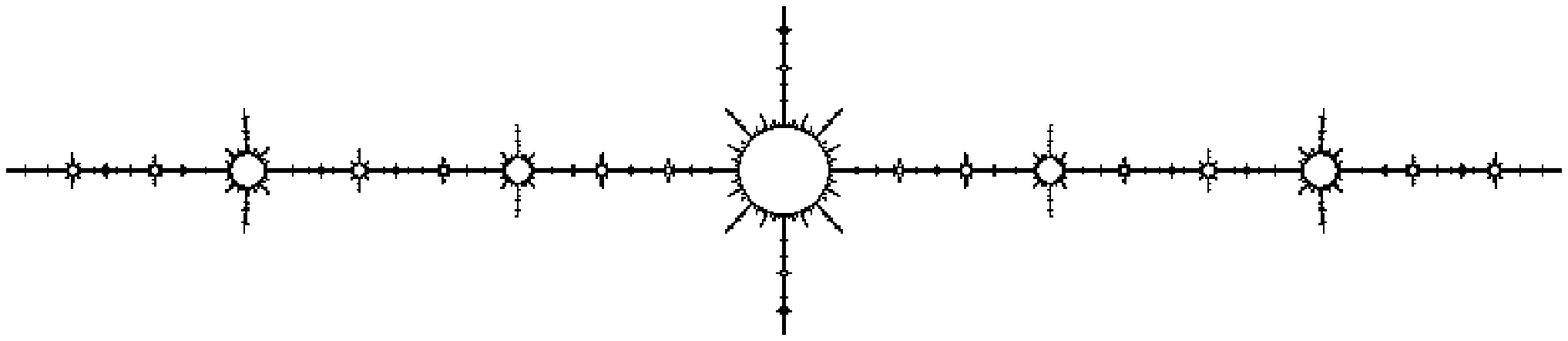}\qquad
      \includegraphics[width=60mm]{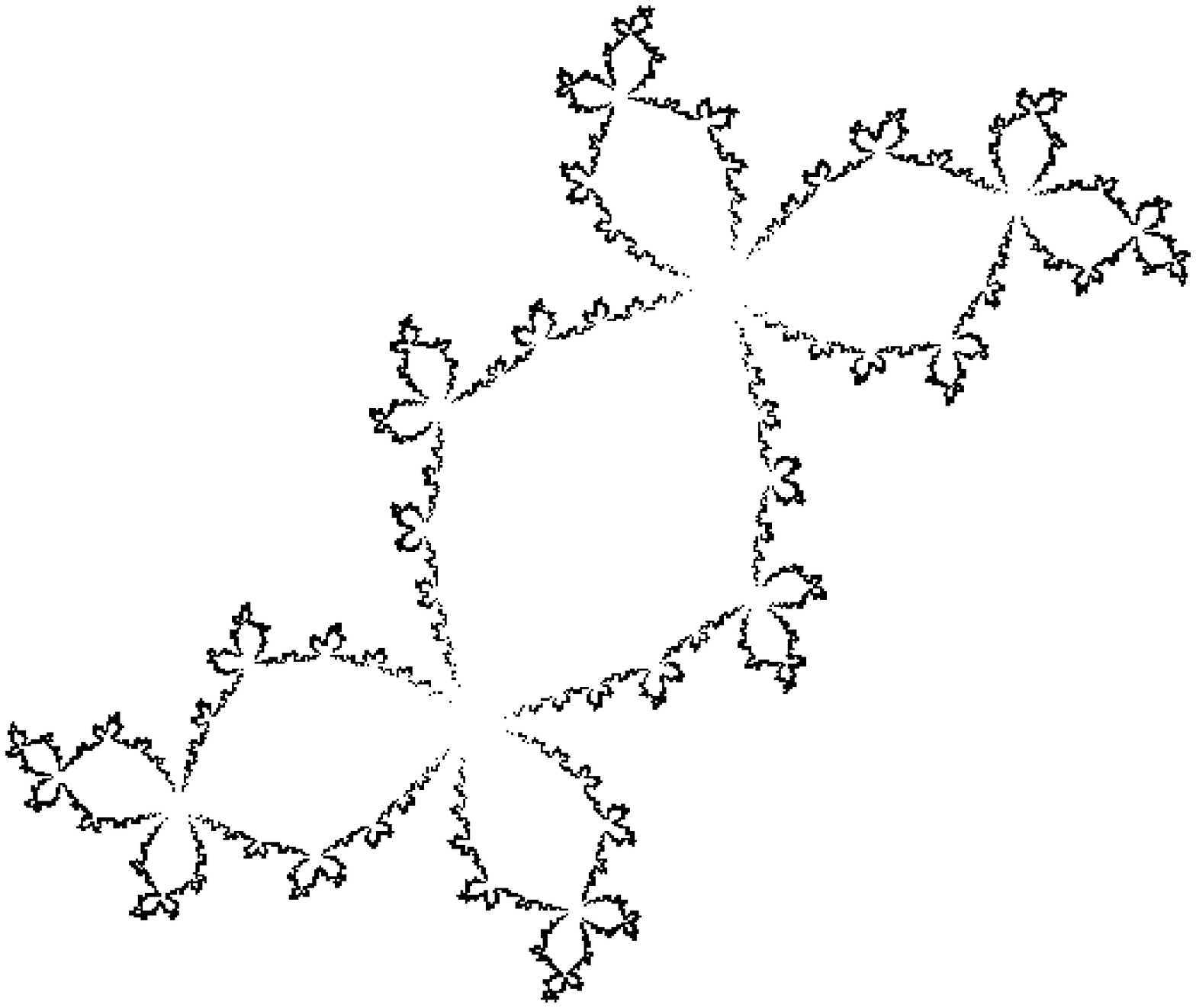}
    \end{center}
    \caption{The Julia sets of the polynomials $z^2+c$, with
      $(c^2+c)^2+c=0$.}
    \label{fig:jul2} 
  \end{figure}
\end{enumerate}

\subsection{The solenoid $\solen{G}$}
\label{ss:solen}
Let us fix a self-similar contracting action of a group $G$ over the alphabet $X$.
Denote by $\xz$ the space of all two-sided infinite sequences over the
alphabet $X$ with the product topology. The elements of this space
have the form
\[\xi=\ldots x_{-3}x_{-2}x_{-1}\;.\;x_0x_1x_2\ldots,\]
with $x_i\in X$, and where the dot marks the place between the
$(-1)$-st and $0$th coordinates. The sequence $x_0x_1x_2\ldots$ is
called the \emph{integer part} of the sequence $\xi$ and is written
$[\xi]$.

The map
\[
\shift:\ldots x_{-3}x_{-2}x_{-1}\;.\;x_0x_1x_2\ldots\mapsto\ldots
x_{-4}x_{-3}x_{-2}\;.\;x_{-1}x_0x_1\ldots
\]
is called the \emph{shift}. It is a homeomorphism of the space
$\xz$.

We say that two sequences $\ldots x_{-2}x_{-1}\;.\;x_0x_1\ldots$ and
$\ldots y_{-2}y_{-1}\;.\;y_0y_1\ldots\in\xz$ are \emph{asymptotically
  equivalent} if there exists a sequence $\{g_k\}_{k=1}^{\infty}$
taking a finite number of different values in $G$, such that
\[
\left(x_{-k}x_{-k+1}x_{-k+2}\ldots\right)^{g_k}=y_{-k}y_{-k+1}y_{-k+2}\ldots.
\]

\begin{defi}
  The \emph{(limit) solenoid} of the contracting action of the group $G$ is the
quotient of the space $\xz$ by the asymptotic equivalence relation.
We denote it by $\solen{G}$.
\end{defi}

Proposition~\ref{pr:regequiv2} gives a more convenient description of
the asymptotic equivalence relation, which is also true for the space
$\xz$ (the only difference is that the path must be infinite in both
directions).

The next proposition shows that the solenoid is uniquely defined by
the dynamical system $(\lims{G}, \si)$:
\begin{proposition}
  \label{pr:invlimtor}
  The space $\solen{G}$ is homeomorphic to the inverse limit of the
  topological spaces
  \[
  \lims{G}\stackrel{\si}{\longleftarrow}
  \lims{G}\stackrel{\si}{\longleftarrow}\cdots.
  \]
\end{proposition}

\begin{defi}
  \label{defi:tiles}
  For a given $w\in\xo$, the \emph{tile} $\til_w$ is the image of the
  set $\{\xi\in\xz: [\xi]=w\}$ under the canonical quotient map
  $\xz\to\solen{G}$.
  
  If $O\subset\xo$ is a $G$-orbit, then the corresponding \emph{leaf}
  is the image in $\solen{G}$ of the set of all $\xi\in\xz$ with
  $[\xi]\in O$.
\end{defi}

It follows from the definition of the asymptotic equivalence that the
integer parts of asymptotically equivalent elements of $\xz$ belong
to one $G$-orbit. Thus the solenoid is a disjoint union of leaves.

The shift $\shift:\xz\to\xz$ induces a homeomorphism $\si$ on the
space $\solen{G}$, which will be also called the \emph{shift}.  We
have
\[\si(\til_w)=\bigcup_{x\in X}\til_{xw},\]
so that every tile $\til_w$ is homeomorphic to a union of $|X|$ tiles.

The \emph{type} of the tile $\til_w$ is the set of elements of the
nucleus $\nuke$ which fix the word $w$.

\begin{proposition}
  \label{prop:tiltyp}
  If the tiles $\til_{w_1}$ and $\til_{w_2}$ have the same type, then
  the map $w_2:u\;.\;w_1\mapsto u\;.\;w_2$, for all $u\in\xmo$,
  induces a homeomorphism $\til_{w_1}\to\til_{w_2}$.
\end{proposition}

Consequently, there exist not more than $2^{|\nuke|}$ different tiles
up to homeomorphism.

Therefore, every contracting action defines an iterated function system
on the sets $\til_w$, with maps $\mathsf{T_x}:\til_w\mapsto\til_{xw}$.
In this system we identify the tiles of same type. Then we get a
finite number of sets. The type of the tile $\til_{xw}$ is uniquely
defined by the type of the tile $\til_w$ and the letter $x$, since
$xw$ is fixed under the action of $h_0\in\nuke$ if and only if
$x^{h_0}=x$ and $h_0|_x=h_1$ for an element $h_1\in\nuke$ fixing $w$.
We call this system the \emph{tiling} iterated function system.

In fact, the term ``tile'' may be misleading, since in general the
tiles $\til_w$ may overlap. See the paper~\cite{reptile} where this
problem is discussed in the abelian case. The following theorem gives
a criterion determining when the tiles intersect only on their
boundaries.
\begin{theorem}
  \label{th:openset}
  If for every element $g$ of the nucleus there exists a finite word
  $v\in\xs$ such that $g|_v=1$, then any two different tiles $\til_w$
  have disjoint interiors.
  
  On the other hand, if for some element $g$ of the nucleus all the
  restrictions $g|_v$, $v\in\xs$ are non-trivial, then there exists a
  tile $\til_w$ covered by other tiles.
\end{theorem}

We say that a contracting action satisfies the \emph{open set
  condition} if every element of its nucleus has a trivial
restriction. It is easy to see that this is equivalent to the
condition that every element of the group has a trivial restriction.

Therefore, in the case of contracting actions satisfying the open set
condition, we get tilings of the leaves in the usual sense. In the
next subsection we show how the self-affine tilings of Euclidean space
appear as tilings associated with self-similar actions of free abelian
groups.

\subsection{Self-affine tilings of Euclidean space}
\label{s:eucrep}
Let us fix a recurrent action $(\Z^n, \xo)$ over the alphabet
$X=\{0,1,\ldots,d-1\}$. Let $R=\{r_0=0,r_1, \ldots, r_{d-1}\}$ be a
digit set for the action. We keep the notation of
Subsection~\ref{ss:abel}, and in particular write the group
additively.

Suppose that the action is finite-state; then by
Theorem~\ref{th:finab} the associated virtual endomorphism yields
a linear map $\phi:\Q^n\to\Q^n$, which is a contraction.

Then for all sequences $w=i_1i_2\ldots\in\xo$ the series
\[F(w)=\sum_{k=1}^\infty\phi^k(r_{i_k})\]
is convergent in $\R^n$. Let
$\til(\phi, R)=\{F(w): w\in\xo\}$ be the set of their sums. We will
call $\til(\phi, R)$ the \emph{set of fractions}, or the \emph{tile}.

In the classical situation of the dyadic numeration system, given by
$\phi(n)=n/2$ and $R=\{0,1\}$, the set of fractions $\til(\phi, R)$
is the interval $[0,1]$.

The set $\til(\phi, R)$ is an attractor of the following affine
iterated function system (see~\cite{hutchinson}) on $\R^n$:
\[\big\{p_i(r)=\phi(r_i+r)\text{ for all }r\in\R^n\big\}, i\in X\]
which means that it is the unique fixed point of the transformation
\[P(C)=\bigcup_{i=0}^{d-1}p_i(C)\]
defined on the space of all non-empty compact subsets of $\R^n$.
Moreover, for any nonempty compact set $C\subset\R^n$, the sequence
$P^n(C)$ converges in this space to $\til(\phi, R)$ with respect to
the Hausdorff metric.  This can be used in practice to obtain approximations of these sets.

Sometimes, the set of fractions is just a rectangle; this is the case
for instance with $\phi=\left(
  \begin{smallmatrix}
    0 & 1 \\
    1/2 & 0
  \end{smallmatrix}
\right)$ and $R=\{(0, 0), (1, 0)\}$ (its set of fractions is then the
square $[0,1]\times [0,1]$), but often the set of fractions has an
interesting fractal appearance.

\begin{figure}[ht]
  \begin{center}
    \includegraphics[height=55mm]{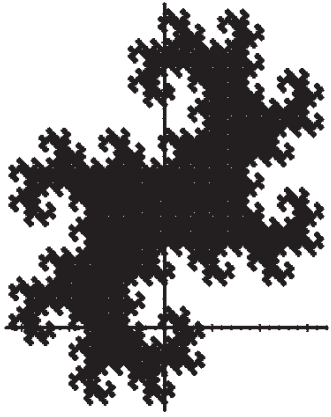}
  \end{center}
  \caption{The set of fractions bounded by the dragon curve.}
  \label{fig:dragon}
\end{figure}

One of the most famous examples is the region bounded by the ``dragon
curve'', corresponding to the case $\phi=\left(
  \begin{smallmatrix}
    1/2 & -1/2 \\
    1/2 & 1/2
  \end{smallmatrix}
\right)$ and $R=\{(0, 0), (1, 0)\}$. An approximation of this set is
shown in Figure~\ref{fig:dragon}.
The associated numeration system can be interpreted as a numeration
system of the Gaussian integers $\Z[i]$ in base $(1+i)$, with the set of
digits $\{0,1\}$. See its discussion in~\cite{knuth}. See
also~\cite{gilbert:three,giordano:cstar}, where similar numeration
systems of the complex numbers and their related set of fractions are
considered.

\begin{figure}[ht]
  \begin{center}
    \includegraphics[width=120mm]{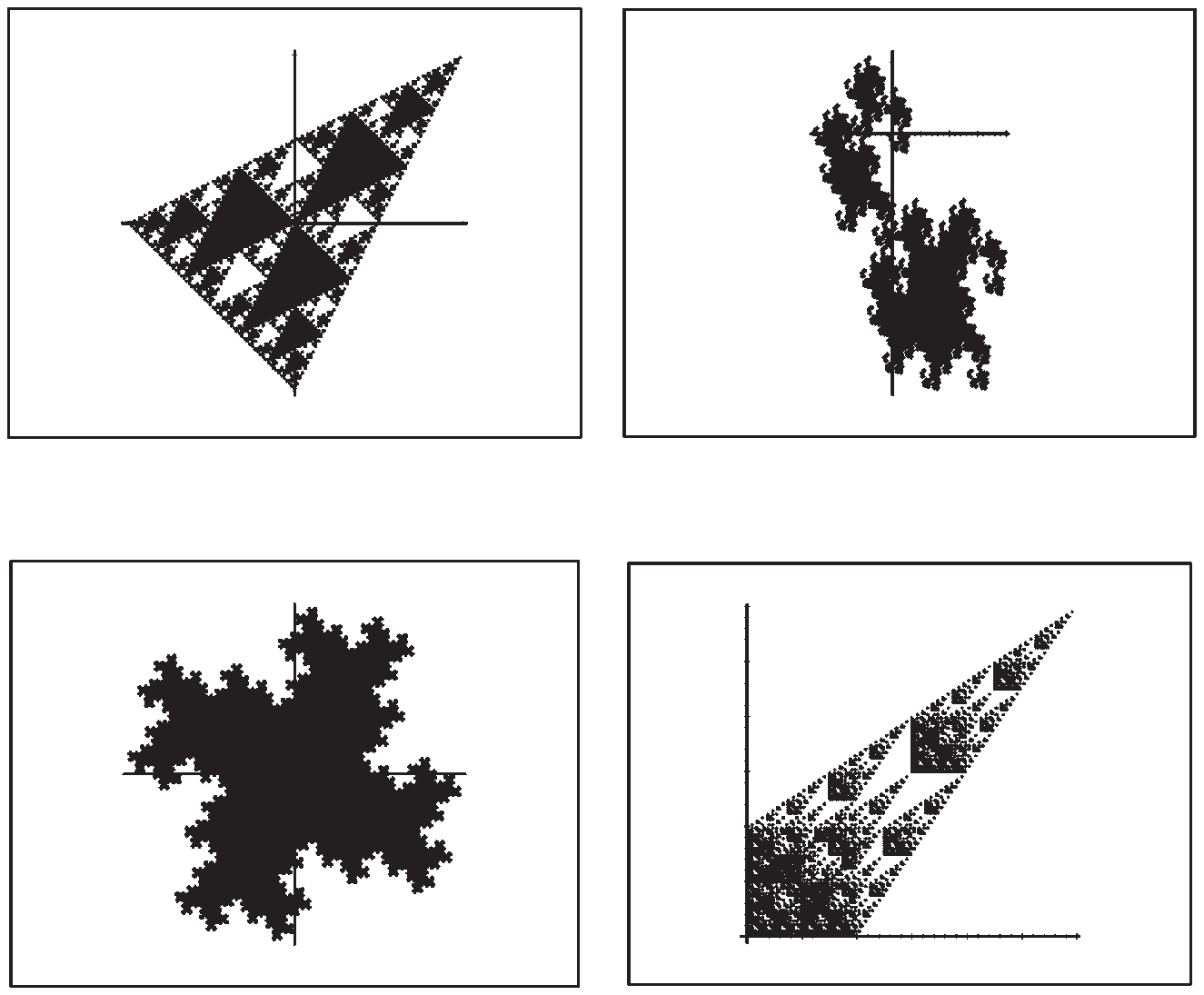}
  \end{center}
  \caption{Tiles of abelian groups}
  \label{fig:examps}
\end{figure}

Other examples of sets of fractions are shown on
Figure~\ref{fig:examps}.

See~\cite{reptile,vince:digtile,pragat:numsyst} for examples of
fraction sets and their properties. These sets have applications to
wavelet theory, computer image processing, toral dynamical systems and
other fields. See also the book~\cite{cuntz_rep} for relations with
the representation theory of the Cuntz $C^*$-algebra.

\begin{theorem}
  \label{pr:frabel}
  Two sequences $(\ldots x_2x_1)$ and $(\ldots y_2y_1)\in\xmo$ are
  asymptotically equivalent if and only if
  \begin{equation*}
    \label{eq:aseqab2}
    \sum_{k=1}^\infty\phi^k(r_{x_k})-\sum_{k=1}^\infty\phi^k(r_{y_k})\in\Z^n.
  \end{equation*}
  
  Two sequences $\xi=(\ldots x_{-2}x_{-1}\;.\;x_0x_1\ldots)$ and
  $\zeta=(\ldots y_{-2}y_{-1}\;.\;y_0y_1\ldots)\in \xz$ are
  asymptotically equivalent if and only if
  \begin{equation*}
    \label{eq:aseqab1}
    \sum_{k=1}^\infty\phi^k(r_{x_{-k}})-\sum_{k=1}^\infty\phi^k(r_{y_{-k}})=
    \sum_{k=0}^\infty\phi^{-k}(r_{y_k})-\sum_{k=0}^\infty\phi^{-k}(r_{x_k}),
  \end{equation*}
  where the left-hand side part is calculated in $\R^n$, while the right one is
  calculated in the closure $\hat\Z^n$ of $\Z^n$; both differences
  must belong to $\Z^n$.
\end{theorem}

Let $\mathcal{L}$ be a leaf of the solenoid $\solen{G}$. Then it
decomposes into the union of its tiles, and thus can be equipped with
the direct limit topology coming from this decomposition.  More
explicitly, a set $A\subseteq\mathcal{L}$ is open in the direct limit
topology if and only if for any finite union of tiles $B$ the set
$A\cap B$ is open in the relative topology of $B$.

\begin{corollary}
\label{cor:abeltil}
  Let $(\Z^n, \xo)$ be a self-similar recurrent finite-state action.
  Then
  \begin{enumerate}
  \item the limit space $\lims{\Z^n}$ is homeomorphic to the torus
    $\mathbb{T}^n=\R^n/\Z^n$;
  \item for every leaf $\mathcal{L}$ (with its direct limit topology)
    of the solenoid $\solen{\Z^n}$ there exists a homeomorphism
    $\Phi:\mathcal{L}\to\R^n$ such that for every tile $\til_w$ of
    $\mathcal{L}$ we have $\Phi(\til_w)=\til(\phi,R)+r(w)$ for some
    $r(w)\in\Z^n$.
  \end{enumerate}
\end{corollary}
Essentially, $r(w)$ is the base-$\phi$ evaluation of $w$.

It follows from the description we obtained of the limit space
$\lims{\Z^n}$ that the shift $\si$ on it coincides with the map on the
torus $(\R/\Z)^n$ given by the linear transformation $\phi^{-1}$.
This map is obviously a $d$-to-$1$ covering. The tiles, just as in the
general case, define a Markov partition for this toral dynamical
system.

\begin{figure}[ht]
  \begin{center}
    \includegraphics{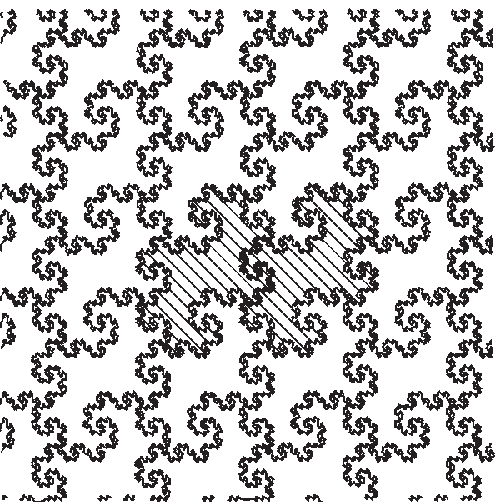}
  \end{center}
  \caption{Plane tiling by dragon curves}
  \label{fig:dragtile}
\end{figure}

Corollary~\ref{cor:abeltil} shows that the tiled leaves of contracting
recurrent self-similar actions of abelian groups are the classical
digit tilings of Euclidean space.  For example, a part of the tiling
by ``dragons'' is shown on Figure~\ref{fig:dragtile}.  The union of
the two marked central tiles is similar to the original tile.

\subsection{Limit spaces of the inverse semigroups}
\label{ss:lsis}
\begin{defi}
  A self-similar inverse semigroup $G$ acting on a topological Markov
  chain $\f\subseteq\xo$ is \emph{contracting} if there exists a
  finite set $\nuke\subset G$ such that for every $g\in G$ and for
  every word $w\in\f$ there exists a finite prefix $v$ of $w$ such
  that
  \[
  \mathsf{T}_v g=h_1 \mathsf{T}_{u_1}+h_2\mathsf{T}_{u_2}+\cdots+h_k\mathsf{T}_{u_k},
  \]
  where $h_i\in\nuke$ and $u_i$ are words of length $|v|$, and
  $\mathsf{T}_u$ denotes the partial permutation
  $\mathsf{T}_u:x\mapsto ux$.
\end{defi}

It is easy to prove that, if $G$ is a group, this definition agrees
with Definition~\ref{defi:contr} of a contracting self-similar
action of a group on a space $\xo$.

In the case of inverse semigroup actions we can not define the limit
space $\lims{G}$, since we do not have a canonical action of the
semigroup on the set of finite words.  However, we can define the
limit solenoid $\solen{G}$ using the action on infinite sequences.

Let $G$ be an inverse semigroup acting on a shift of finite type
$\f\subseteq\xo$. Define $\f^\Z\subseteq\xz$ to be the two sided shift
space defined by the same set of admissible words as $\f$.  More
prosaically, $\f$ is the set of all biinfinite sequences $\ldots
x_{-2}x_{-1}. x_0x_1x_2\ldots\in\xz$ such that for every $n\in\Z$ the
sequence $x_nx_{n+1}x_{n+2}\ldots$ belongs to $\f$.

We define the asymptotic equivalence relation on $\f^\Z$ for the
action of the semigroup $G$ exactly in the same way as it is defined
on the space $\xz$ for group actions. The \emph{limit solenoid}
$\solen{G}$ is then the quotient of the topological space $\f^\Z$ by
the asymptotic equivalence relation.

We also define the tiles of the solenoid in the same way as it is done
for groups (Definition~\ref{defi:tiles}).

By the \emph{type} of a tile $\til_w$ in the case of an inverse
semigroup we mean a pair $(D, F)$ of subsets of the nucleus. The set
$D$ is the set of those elements $g$ of the nucleus for which $w$ is
contained in the domain $\dom g$. The set $F$ is, as in the case of
the group actions, the set of the elements fixing the word $w$.

Proposition~\ref{prop:tiltyp} remains true for contracting actions of
inverse semigroups.  Therefore, for self-similar contracting inverse
semigroups the tiling iterated function systems are also well defined.

\paragraph{The Fibonacci transformations.} The semigroup generated by
the Fibonacci transformations is contracting with contracting
coefficient $\tau^{-1}$, where $\tau=\frac{1+\sqrt{5}}{2}$.

The corresponding iterated function system on the tiles is the
\emph{Fibonacci iterated function system} described among the Examples
of Subsection~\ref{ss:dirlim}.

\paragraph{Penrose tilings.} The semigroup related to the Penrose
tilings is also contracting, with contraction coefficient $\tau^{-1}$.
The tiling iterated function system on the tiles corresponding to this
semigroup is exactly the Penrose iterated functions system.

\section{Hyperbolic spaces and groups}
\label{s:hyp}
\subsection{Definitions}
\begin{defi}
  A metric space $(\X, d)$ is $\delta$-hyperbolic (in the
  sense of M.~Gromov) if for every $x_0,x, y, z\in\X$ the
  inequality
  \[
  \langle x\cdot y\rangle_{x_0}\geq\min\left\{\langle x\cdot
    z\rangle_{x_0}, \langle y\cdot z\rangle_{x_0} \right\}-\delta
  \]
  holds, where
  \[
  \langle x\cdot y\rangle_{x_0}=\frac{1}{2}\left(d(x_0,x)+d(x_0,y)-d(x, y)\right)
  \]
  denotes the \emph{Gromov product} of the points $x$ and $y$ with
  respect to the base point $x_0$.
\end{defi}

Examples of hyperbolic metric spaces are all bounded spaces (with
$\delta$ equal to the diameter of the space), trees (which are
0-hyperbolic) and the usual hyperbolic space $\mathbb{H}^n$, which is
hyperbolic with $\delta=\log 3$.

\begin{defi}
  A finitely generated group is \emph{hyperbolic} if it is hyperbolic
  as a word metric space.
\end{defi}

The definition is independent of the choice of the generating set
with respect to which the word metric is defined. For the proof of
this fact, and for the proof of other properties of hyperbolic groups,
see~\cite{gro:hyperb,cdp:hyp,curnpap:symb,ghys-h:gromov}.

Here is a short summary of examples and properties of hyperbolic
groups:
\begin{enumerate}
\item Every finite group is hyperbolic.  
\item Every finitely generated free group is hyperbolic.
\item If $G_1$ is a subgroup of finite index of the group $G$, then
  $G$ is hyperbolic if and only if $G_1$ is hyperbolic.
\item A free product of two hyperbolic groups is hyperbolic.
\item The fundamental group of a compact Riemannian space of negative
  curvature is hyperbolic.
\item A hyperbolic group is finitely presented.
\item The word problem in a hyperbolic group is solvable in linear
  time.
\item A subgroup of a hyperbolic group either contains the free
  subgroup $F_2$, or is a finite extension of a cyclic group.
\item Hyperbolic groups have a rational growth series, and are either
  virtually cyclic or have uniformly exponential growth (see
  Subsection~\ref{ss:grgroups}).
\end{enumerate}

\subsection{The boundary of a hyperbolic space}
\label{ss:boundary}

Let $(\X, d)$ be a hyperbolic space. We say that a sequence
$(x_n)_{n\geq 1}$ of points of the space \emph{converges to infinity}
if for a fixed $x_0\in\X$
\[
\langle x_n\cdot x_m\rangle_{x_0}\to+\infty\qquad\text{ when $n, m\to+\infty$}.
\]

It is easy to prove that the definition does not depend on the choice
of the point $x_0$.

Two sequences $(x_n), (y_n)$ converging to infinity are said to be
\emph{equivalent} if $\langle x_n\cdot y_n\rangle_{x_0}\to+\infty$,
when $n\to+\infty$.

This definition also does not depend on $x_0$. The quotient of the set
of sequences converging to infinity by this equivalence relation is
called the \emph{boundary} of the hyperbolic space $\X$, and
is denoted $\partial\X$.

If a sequence $\left(x_n\right)$ converges to infinity, then its
\emph{limit} is the equivalence class $a\in\partial\X$ to
which the sequence $\left(x_n\right)$ belongs, and we say that
$\left(x_n\right)$ \emph{converges to $a$}.

If $a, b\in\partial\X$ are two points of the boundary, then
their \emph{Gromov product} is defined as
\[
\langle a\cdot b\rangle_{x_0}=\sup_{\left(x_n\right)\in a,
\left(y_m\right)\in b} \liminf_{m, n\to\infty}\langle x_n\cdot
y_m\rangle_{x_0}.
\]

For every $r>0$ define
\[
V_r=\{(a, b)\in\partial\X\times\partial\X:
\langle a\cdot b\rangle\geq r\}_{x_0}.
\]
Then the set $\{V_r: r\geq 0\}$ is a fundamental neighborhood basis of
a uniform structure on $\partial\X$ (see~\cite{bour:top}
and~\cite{ghys-h:gromov} for the necessary definitions and proofs). We
topologize the boundary $\partial\X$ by this uniform
structure.

Another way to define the topology on the boundary is to introduce the
\emph{visual metric} on it.

Namely, let $\X$ be a geodesic hyperbolic metric space with
a base-point $x_0$. Recall that a metric space $\X$ is said
to be \emph{geodesic} if any two points $x, y\in\X$ can be
connected by a path defined by an isometric embedding of the real
segment $[0,d(x, y)]$ into the space.

Let $a\in(0,1)$ be a number close to $1$. Then for every path
$\gamma:[0,t]\to\X$ we define its $a$-length $l_a(\gamma)$
as the integral
\[l_a(\gamma)=\int_{x\in\gamma[0,t]}a^{-d(x, x_0)}dx.\]
Then we define the $a$-distance between two points $x,
y\in\X$ as the infimum of $a$-lengths of the continuous
paths connecting $x$ and $y$. There exists $a_0\in(0,1)$ such that for
all $a\in(a_0,1)$ the completion of the $a$-metric on $\X$
is $\X\cup\partial\X$. The restriction of the
extended metric on $\partial\X$ is called \emph{visual
  metric} on the boundary.

The boundary of a hyperbolic group has a rich self-similar structure
(see Section~\ref{s:fpds} and the book~\cite{curnpap:symb}). Some of
the classical fractals (for instance, the Sierpinski carpet) can also
be realized as the boundary of a hyperbolic group
(see~\cite{curnpap:symb,kapklein:bound} and their bibliography).

\paragraph{Questions.} What topological spaces can be realized as the
boundaries of hyperbolic groups?  How can one compute the Hausdorff
dimension of the boundary of a hyperbolic group respectively to the
visual metric?

\begin{defi}
\label{def:limitset}
  Suppose a group $G$ acts by isometries on a hyperbolic space
  $\X$. Then its \emph{limit set} is the set of all the
  limits in $\partial\X$ of sequences of the form
  $x_k=x_0^{g_k}, k\geq 1$, where $g_k\in G$ and $x_0\in\X$
  is a fixed point.
\end{defi}

It is easy to prove that the limit set does not depend on the point
$x_0$. An important case is the limit set of a Kleinian group acting on $\mathbb{H}^n$. It is a subset
of $S^{n-1}=\partial \mathbb{H}^n$.

\paragraph{Example.} Let $\Gamma$ be the group generated by the four
inversions $\gamma_i$ defining the self-similarity structure of the
Apollonian net (see ``The group associated with the Apollonian
gasket'' in Subsection~\ref{ss:gofss}).  Since the group of conformal automorphisms of the Riemannian
sphere is isomorphic to the group of isometries of the hyperbolic space $\mathbb{H}^3$, the
action of $\Gamma$ on the sphere is extended in the standard way  to an action
by isometries on the space $\mathbb{H}^3$ (see~\cite{elstrgrun:groups}), where the sphere is
identified with the boundary $\partial\mathbb{H}^3$. In the extended
action the generators $\gamma_i$ are reflections.  For instance, in
the Kleinian model of the space $\mathbb{H}^3$, the generators
$\gamma_i$ can be defined by the matrices
\[
\left(\begin{smallmatrix}-1 & 2 & 2 & 2\\ 0 & 1 & 0 & 0 \\ 0 & 0 & 1 & 0 \\ 0 & 0 & 0 & 1\end{smallmatrix}\right),\quad
\left(\begin{smallmatrix} 1 & 0 & 0 & 0\\ 2 & -1 & 2 & 2 \\ 0 & 0 & 1 & 0 \\ 0 & 0 & 0 & 1\end{smallmatrix}\right),\quad
\left(\begin{smallmatrix} 1 & 0 & 0 & 0\\ 0 & 1 & 0 & 0 \\ 2 & 2 & -1 & 2 \\ 0 & 0 & 0 & 1\end{smallmatrix}\right),\quad
\left(\begin{smallmatrix} 1 & 0 & 0 & 0\\ 0 & 1 & 0 & 0 \\ 0 & 0 & 1 & 0 \\ 2 & 2 & 2 & -1\end{smallmatrix}\right).
\]

It is easy to see that the limit set of the obtained group is the
Apollonian net $\mathcal{P}$.

Let $\Gamma_1<\Gamma$ be the subgroup of orientation-preserving
transformations. It is of index $2$ in $\Gamma$ and is a free group of
rank 3 generated by the elements $\gamma_1\gamma_2, \gamma_1\gamma_3$
and $\gamma_1\gamma_4$. It is a subgroup of the Moebius group and its
limit set coincides with the limit set of the group $\Gamma$, i.e.,
with the Apollonian net.

\subsection{The limit space of a self-similar action as a hyperbolic boundary}
\begin{defi}
  \label{defi:ssimcomp}
  Let $(G, \xo)$ be a self-similar action of a finitely generated
  group. For a given finite generating system $S$ of the group $G$ we
  define the \emph{self-similarity complex} $\Ss(G, S)$ as the
  $1$-complex with set of vertices $\xs$, and with two vertices $v_1,
  v_2\in\xs$ belonging to a common edge if and only if either
  $v_i=xv_j$ for some $x\in X$ (the edges of the \emph{first type}) or
  $v_i^s=v_j$ for some $s\in S$ (the edges of the \emph{second type});
  here $\{i, j\}=\{1,2\}$.
\end{defi}

The set of vertices of the self-similarity complex splits into the
levels $X^n, n\in\N$. Every edge of the first type connects two
vertices from neighboring levels, while every edge of the second type
connects the vertices of the same level. The set of edges of the
second type spans $\G(G,S,\xs)$, the disjoint union of all the finite
Schreier graphs $\G_n(G, S)$ of the group $G$.

\begin{figure}[ht]
  \begin{center}
    \includegraphics{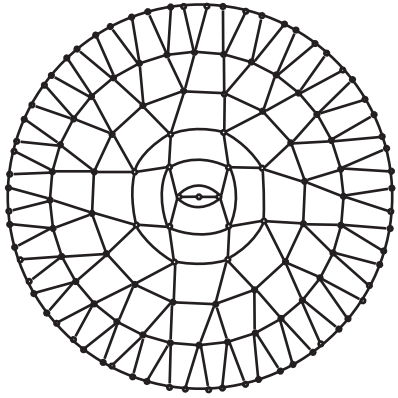}
  \end{center}
  \caption{The self-similarity complex of the adding machine}
  \label{fig:admcom}
\end{figure}

A part of the self-similarity complex of the adding machine is shown
in Figure~\ref{fig:admcom}.

If all the restrictions of the elements of the generating set $S$ also
belong to $S$ (this is the case, for instance, when $S$ is the nucleus
$\nuke$ of $G$), then the self-similarity complex $\Ss(G, S)$ is an
\emph{augmented tree} in the sense of V.~Kaimanovich (see~\cite{kaim:augtree}).

\begin{theorem}
  \label{th:hypsp}
  If the action of a finitely generated group $G$ is contracting, then
  the self-similarity complex $\Ss(G, S)$ is a Gromov-hyperbolic
  space.
  
  The limit space $\lims{G}$ is then homeomorphic to the hyperbolic
  boundary $\partial\Ss(G, S)$ of the self-similarity complex $\Ss(G,
  S)$. Moreover, there exists a homeomorphism
  $\lims{G}\to\partial\Ss(G, S)$, which makes the diagram
  \begin{equation*}\label{eq:cdb}
    \xymatrix{
      {\xmo} \ar[d]_{\pi} \ar[dr]^{\ell}\\
      {\lims{G}} \ar[r] & {\partial\Ss(G, S)}}
  \end{equation*}
  commutative. Here $\pi$ is the canonical projection and $\ell$
  carries every sequence $\ldots x_2x_1\in\xmo$ to its limit
  \[\lim_{n\to\infty}x_nx_{n-1}\ldots x_1\in\partial\Ss(G, S).\]
\end{theorem}


\section{Finitely presented dynamical systems and semi-Markovian spaces}
\label{s:fpds}
Let $\left(\{F_v\}_{v\in V}, \{\phi_e\}_{e\in E}\right)$ be a
self-similarity structure on a compact space $F=\cup_{v\in V} F_v$.

Suppose that the self-similarity structure is such that for every
infinite path $\ldots e_2e_1$ in the structural graph of the iterated
function system, the intersection
\[
\bigcap_{k\geq 1}\phi_{e_1}\left(\phi_{e_2}\left(\ldots \phi_{e_k}\left(F_{\be(e_k)}\right)\right) \ldots\right),
\]
contains exactly one point of the set $F$.

Let $\f\subset E^\omega$ be the set of all sequences $e_1e_2\ldots$
for which $\ldots e_2e_1$ is a path in the structural graph of the
self-similarity structure.

\begin{proposition}
  The map $\Pi:\f\to F$ defined by the rule
  \[
 \left\{\Pi(e_1e_2\ldots)\right\}= \bigcap_{k\geq
    1}\phi_{e_1}\left(\phi_{e_2}\left(\ldots
      \phi_{e_k}\left(F_{\be(e_k)}\right)\right) \ldots\right)
  \]
  is continuous.
\end{proposition}

The space $E^\omega\times E^\omega$ is naturally identified with the
one-sided shift space $(E\times E)^\omega$ via the homeomorphism
\[(e_1e_2\ldots, f_1f_2\ldots)\mapsto (e_1, f_1)(e_2, f_2)\ldots.\]

\begin{defi}
  The self-similarity structure $\left(\{F_v\}_{v\in V},
    \{\phi_e\}_{e\in E}\right)$ on the space $F$ is \emph{finitely
    presented} if the set
  \[
  \left\{(w_1, w_2): \Pi(w_1)=\Pi(w_2)\right\}\subseteq E^\omega\times E^\omega
  \]
  is a subshift of finite type in $(E\times E)^\omega$.
  
  A topological space $F$ which has a finitely presented
  self-similarity structure is called \emph{Markovian}.
\end{defi}

If the maps $\phi_e^{-1}$ are restrictions of a map $f:F\to F$, and
the self-similarity structure is finitely presented, then the
dynamical system $(F, f)$ is also called \emph{finitely presented}
(see~\cite{curnpap:symb}).

The notion of a finitely presented dynamical system and a
(semi-)Markovian space was formulated for the first time by M.~Gromov
(see~\cite{gro:hyperb}).

We recall here the definitions of a shift of finite type and of a
finitely presented dynamical system in the general case of a semigroup
action, following~\cite{curnpap:symb}. The one-sided shift is then the
case of the semigroup $\N$; the bilateral shift is the case of the
group $\Z$.

Let $X$ be an alphabet and let $G$ be a semigroup. The \emph{shift
  space} is the Cartesian power $X^G$, i.e., the set of all functions
$G\to X$. We put on the shift space the direct product topology. The
semigroup $G$ acts on the elements of the shift space $X^G$ by the
rule
\[\xi^g(h)=\xi(gh).\]

A subset $C\subseteq X^G$ is called a \emph{cylinder} if there exists
a finite set $F\subseteq G$ and a finite set $A$ of functions $F\to X$
such that $\xi$ is an element of $C$ if and only if the restriction
$\xi|_F$ belongs to $A$. Any cylinder is a closed and open subset of
the space $X^G$.

In the same way as for the unilateral shift, we can identify the
direct product $X^G\times X^G$ with the full shift space $(X\times
X)^G$ over the alphabet $X\times X$. Namely, the pair $(\xi, \zeta)\in
X^G\times X^G$ is identified with the function $h\mapsto (\xi(h),
\zeta(h))$.

\begin{defi}
  A subset $\f\subseteq X^G$ is a \emph{subshift of finite type} if
  there exists a cylinder $C\subset X^G$ such that
  \[\f=\bigcap_{g\in G} C^{g^{-1}},\]
  where $C^{g^{-1}}=\{\xi: \xi^g\in C\}$.
\end{defi}

It is easy to see that every subshift of finite type is a closed
$G$-invariant subset of $X^G$.

\begin{defi}
  A dynamical system $(\X, G)$ (i.e., a topological space
  $\X$ with a continuous action of a semigroup $G$) is a
  \emph{system of finite type} if there exists a finite alphabet $X$,
  a subshift of finite type $\f\subseteq X^G$ and a continuous,
  surjective and $G$-equivariant map $\pi:\f\to\X$.
  
  A dynamical system $(\X, G)$ is \emph{finitely presented}
  if in addition the set $\{(\xi, \zeta)\in \f\times\f:
  \pi(\xi)=\pi(\zeta)\}$ is a subshift of finite type in $(X\times
  X)^G$.
\end{defi}

\begin{defi}
  A dynamical system $(\X, G)$ is \emph{expansive} if there
  exists an open set $U\subset\X\times\X$ such
  that
  \[\Delta=\bigcap_{g\in G}U^{g^{-1}},\]
  where $\Delta=\{(x, x) | x\in\X\}$ is the diagonal in $\X\times\X$.
\end{defi}

For the proof of the following theorem see~\cite{curnpap:symb} (see
also~\cite{fried}, where the case $G=\Z$ is considered):
\begin{theorem}
  Let $(\X, G)$ be a dynamical system.  Let $\f\subseteq
  X^G$ be a subshift of finite type and $\pi:\f\to\X$ a
  continuous, surjective and $G$-equivariant map. Then the space
  $\{(\xi, \zeta)\in \f\times\f: \pi(\xi)=\pi(\zeta)\}$ is a subshift
  of finite type if and only if the system $(\X, G)$ is
  expansive.
\end{theorem}

\begin{corollary}
  A dynamical system $(\X, G)$ is finitely presented if and
  only if it is both expansive and of finite type.
\end{corollary}

\paragraph{Examples.} \begin{enumerate}
\item Any subshift of finite type is a finitely presented dynamical
  system.
\item Let $G$ be a hyperbolic group. Then the dynamical system
  $(\partial G, G)$ is finitely presented. See the
  book~\cite{curnpap:symb}, where different finite presentations of
  this dynamical system are given.
\end{enumerate}

\begin{defi}
  A Hausdorff topological space $\X$ is
  \emph{semi-Markovian} if there exists a finite alphabet $X$, a
  cylinder $C\subseteq\xo$, a subshift of finite type $\f\subseteq\xo$
  and a continuous surjection $\Pi:C\cap\f\to\X$ such that
  the set $\{(w_1, w_2): \Pi(w_1)=\Pi(w_2)\}$ is an intersection of a
  cylinder and a subshift of finite type in $(X\times X)^\omega$.
\end{defi}

\paragraph{Examples.} \begin{enumerate}
\item Obviously the shifts of finite type and the Cantor space $\xo$
  are semi-Markovian spaces.
\item The real segment $[0,1]$ is a semi-Markovian space. The
  semi-Markovian presentation is defined, for instance, by the dyadic
  expansion of the reals.
\item In general, any finite simplicial complex is semi-Markovian.
\item The boundary of a torsion-free hyperbolic group is semi-Markovian.
\end{enumerate}
For more details on these examples see~\cite{curnpap:symb}.


\section{Spectra of Schreier graphs and Hecke type operators}
\label{s:spectrum}
Let $\G$ be a graph, and consider the Hilbert space
$\hilb=\ell^2(\G,\deg)$, the complex space of finite-norm functions on
$V(\G)$ determined by the scalar product
\[\scalar{f}{g}=\sum_{v\in V(\G)}f(v)\overline{g(v)}\deg(v).\]
The adjacency of $\G$ defines an operator $M$ on $\hilb$, called
its \emph{Markov} or \emph{transition} operator, by
\[M(f)(v) = \frac1{\deg(v)}\sum_{\substack{e\in E(\G)\\ \alpha(e)=v}}f(\omega(e)).\]
The \emph{spectrum} of $\G$ is defined as the spectrum of $M$. If $\G$
has bounded degree, then $M$ is a bounded operator, and hence the
spectrum of $\G$ is compact.

The spectrum of $\G$ is intimately connected to random walk properties
of $\G$: the Green function is essentially the resolvent of $M$; the
graph $\G$ is amenable if and only if the spectral radius of $M$ is
$1$, as was shown by Kesten~\cite{kesten:amen}.

\begin{defi}
  Let $G$ be a group with fixed generating set $S$, and let $\pi$ be a
  unitary representation of $G$ on a Hilbert space $\hilb$. The
  associated \emph{Hecke type operator} is
  \[H_\pi = \frac1{|S|}\sum_{s\in S}\pi(s)\in\mathcal B(\hilb).\]
\end{defi}
In the case when $S$ is symmetric (meaning $S=S^{-1}$), $H$ is a
self-adjoint operator. The \emph{spectrum} of $\pi$ is the spectrum of
its Hecke type operator $H_\pi$. Note that it depends also on $S$,
although this is not apparent in the notation.

In our situation, we have a group $G$ acting on the boundary $\xo$ of
a $d$-regular tree $T(X)$, and hence $G$ acts unitarily on the Hilbert
space $\hilb=L^2(\xo,\mu)$, where $\mu$ is the uniform measure.  By
restriction to the subsets $x\xo$ for $x\in X$, we get an isomorphism
$\phi:\hilb\cong\hilb^X$ given by $\phi(f)(x)=f\mathsf{T}_x:w\mapsto
f(xw)$. We use this isomorphism to obtain finite approximations to the
spectrum of the representation on $\hilb$.

This isomorphism also makes possible to interprete the self-similarity
of the group actions in terms of operator algebras. We do this in
Section~\ref{s:cstar}.

Consider the constant functions $\hilb_0=\C\subset\hilb$, and define
subspaces $\hilb_{n+1}=\phi^{-1}(\hilb_n^X)$ for all $n\ge0$. All
$\hilb_n$ are closed $G$-invariant subspaces of $\hilb$, and
$\dim\hilb_n=d^n$. We will compute the spectrum of
$\pi_n=\pi|_{\hilb_n}$; since $\pi_n$ contains $\pi_{n-1}$ as a
subrepresentation, we write $\pi_n\ominus\pi_{n-1}$ an orthogonal
complement. The convergence of the spectrum of $\pi_n$ towards the
spectrum of $\pi$ is given by the following result
from~\cite{bgr:spec}:
\begin{theorem}
  $\pi$ is a reducible representation of infinite dimension whose
  irreducible components are precisely those of the
  $\pi_n\ominus\pi_{n-1}$ (and thus are all finite-dimensional).
  Moreover
  \[\spec(\pi)=\overline{\bigcup_{n\ge0}\spec(\pi_n)}.\]
  $H_\pi$ has a pure-point spectrum, and its spectral radius
  $r(H_\pi)=s\in\R$ is an eigenvalue.
\end{theorem}

Consider also the following family of representations: let $H$ be a
subgroup of $G$, and let $G$ act on $\ell^2(G/H)$ by left translations.
The resulting unitary representation $\rho_{G/H}$ is called
\emph{quasi-regular}; its Hecke-type operator is the Markov operator
of the Schreier graph of $\G(G,S,G/H)$. If $H=1$, this is the
\emph{regular} representation, and more generally if $H$ is normal
$\rho_{G/H}$ is the regular representation of the quotient.

Let $P_n$ be the stabilizer in $\xs$ of a sequence $w_n$ of length
$n$.  Then $\rho_{G/P_n}$ is naturally isomorphic to $\pi_n$, since
the cosets of $P_n$ are in $1-1$ correspondence with the orbit $X^n$
of $w$. Consider also $P$, the stabilizer of an infinite sequence
$w\in\xo$, and the quasi-regular representation $\rho_{G/P}$.
\begin{theorem}
  The spectrum of $\rho_{G/P}$ is contained in the spectrum of $\pi$,
  and if either $G/P$ or $P$ is amenable, these spectra coincide. If
  $P$ is amenable, then they are contained in the spectrum of
  $\rho_G$.

  The spectral radius of $\rho_{G/P}$ is never an eigenvalue.
\end{theorem}
Note then that, under the amenability condition, $\rho_{G/P}$ and
$\pi$ are two distinct representations with same spectrum.

Each $\pi_n(g)$ is a $d^n\times d^n$-permutation matrix, giving $g$'s
action on $\hilb_n$. By decomposition, each $\pi_n(g)$ is a $d\times
d$-matrix, whose entries are either $0$ or $\pi_{n-1}(g')$ for some
$g'\in G$. In particular, if $G$ is contracting, there is a finite
generating set $K$ such that for each $k\in K$, the matrix $\pi_n(k)$
has entries of the form $0$ or $\pi_{n-1}(k')$ for some $k'\in K$.

The computation of the spectrum of $\pi_n$ can be seen in a more
general context, that of \emph{Frobenius determinants}. Given a
finite-dimensional representation $\rho$ of a group $G$, let
$\{\lambda_g\}$ be a set of formal variables indexed by $G$, and
define the Frobenius determinant as
\[\Phi_\rho=\det\Big(\sum_{g\in G}\lambda_g\rho(g)\Big).\]
It is known that the factorization of $\Phi_\rho$ in prime components
exactly parallels the decomposition of $\rho$ in irreducible
components.

Even though $G$ is not assumed to be finite, it is clear that
$\Phi_\rho$ depends only on $G/\ker\rho$. In our case, $\rho=\pi_n$ is
a permutational representation, so $G/\ker\pi_n$ is a finite group
acting on $X^n$.

The spectrum of $\rho$ is obtained by substituting $\lambda_s=1/|S|$
for all $s\in S$, preserving $\lambda_1$, and setting $\lambda_g=0$
for all other variables, in $\Phi_\rho$; and then solving the
resulting polynomial in $\lambda_1$.  The following somewhat
miraculous facts make it sometimes possible to compute this particular
value of $\Phi_\rho$:
\begin{itemize}
\item If $G$ is sufficiently contracting, it may be possible to
  consider only a finite subset $V=\{\lambda_g\}_{g\in F}$ of
  variables, and to express $\Phi_{\pi_n}(V)$ in the form
  $R(\Phi_{\pi_{n-1}}(S(V)))$ for a rational function $R\in\C(z)$ and
  $S\in\C(V)$;
\item Even though $\Phi_\rho$ may not factor in many low-degree terms,
  it may happen that after imposing some extra conditions among the
  $\lambda_g$, that are both weak enough so that they allow
  $\lambda_s=1$ for $s\in S$ and $\lambda_1=-|S|\lambda$ to be chosen,
  and strong enough so that a recursion still holds between
  $\Phi_{\pi_n}$ and $\Phi_{\pi_{n-1}}$, the determinant $\Phi_\rho$
  does factor nicely.
\end{itemize}
These two properties hold for the Grigorchuk group, the Gupta-Sidki
group, the Fabrykowski-Gupta group and the Bartholdi-Grigorchuk
group~\cite{bgr:spec}, and for the lamplighter
group~\cite{gr_zu:lamp}, allowing a complete computation of the
spectrum of $\pi$.

\begin{figure}[ht]
  \begin{center}
    \includegraphics{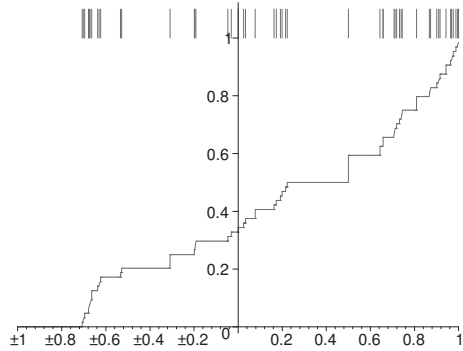}
  \end{center}
  \caption{The spectrum and (counting) spectral measure of $\img{z^2-1}$}
  \label{fig:g2}
\end{figure}

We make this process explicit for the group $\img{z^2-2}$. It is the group, generated by the
transformations $a=(b, 1)\sigma, b=(a, 1)$ (see Subsection~\ref{ss:exmono}). Define the
polynomial
\[\Phi_k(\lambda; \lambda_1, \lambda_2)=
\det\Big(\lambda+\lambda_1\big(\pi_k(a)+\pi_k(a^{-1})\big)+
\lambda_2\big(\pi_k(b)+\pi_k(b^{-1})\big)\Big).
\]
Write $\Lambda=\lambda+2\lambda_2$. Then we have
\begin{align*}
  \Phi_{k+1}(\lambda;\lambda_1,\lambda_2) &=\det\begin{pmatrix}
    \lambda+\lambda_2\big(\pi_k(a)+\pi_k
    (a^{-1})\big)     & \lambda_1(1+\pi_k(b^{-1})) \\
    \lambda_1(1+\pi_k(b)) & \Lambda\end{pmatrix}\\
  &= \Phi_k(\Lambda\lambda-2\lambda_1^2;\Lambda\lambda_2, -\lambda_1^2).
\end{align*}
Consider the curve $\Lambda$ in $\R^3$, and the polynomial
transformation
\[
P:(\lambda;\lambda_1,\lambda_2)\mapsto(\Lambda\lambda-2\lambda_1^2;\Lambda\lambda_2,
-\lambda_1^2)
\]
on $\R^{r+1}$; consider the intersection of $P^{-k}(\Lambda)$ with
$\{\lambda_1=\lambda_2=-1/4\}$. This is the spectrum of $\pi_k$. Therefore, the spectrum of
$\pi$ is the intersection of $\overline{\bigcup_{k\ge 0}P^{-k}(\Lambda)}$ with
$\{\lambda_1=\lambda_2=-1/4\}$.  This spectrum seems to be a Cantor set of null measure, although
we have no complete proof of this fact; see Figure~\ref{fig:g2}. The problem is
that $\Phi_k(\lambda;-1/4,-1/4)$ does not have any nice factorization properties
(it has an irreducible factor of degree $2^{k-2}+1$)

\subsection{The Fabrykowski-Gupta group}
\label{subs:fgspectrum}
The spectrum of the quasi-regular representation of this group,
introduced in Subsection~\ref{ss:exlimit}, was computed by the first
two named authors in~\cite{bgr:spec}. In standard notation, it is the group
\[G = \langle a=(123), \quad s=(a,1,s)\rangle.\]
We only quote the statements, whose proof appear in the
above-mentioned paper.  Denote by $a_n$ and $s_n$ the matrices of the
action on $\{0,1,2\}^n$. We have
\begin{gather*}
  a_0 = s_0 = (1),\\
  a_n = \begin{pmatrix}0&1&0\\0&0&1\\1&0&0\end{pmatrix},\qquad s_n = \begin{pmatrix}a_{n-1}&0&0\\0&1&0\\0&0&s_{n-1}\end{pmatrix}.\\
  \intertext{Let us define operators}
  A_n = a_n + a_n^{-1},\qquad S_n = s_n + s_n^{-1},\qquad
  Q_n(\lambda,\mu) = S_n + \lambda A_n - \mu;\\
  \intertext{then the combinatorial Laplacian of $G$ on $X^n$ is}
  \Delta_n = a_n + a_n^{-1} + s_n + s_n^{-1} =
  \begin{pmatrix}A_{n-1}&1&1\\1&2&1\\1&1&S_{n-1}\end{pmatrix}.\\
  \intertext{Define the polynomials}
  \alpha=2-\mu+\lambda,\qquad \beta=2-\mu-\lambda,\\
  \gamma=\mu^2-\lambda^2-\mu-2,\qquad\delta=\mu^2-\lambda^2-2\mu-\lambda.
\end{gather*}

\begin{lemma}
  \begin{align*}
    \det Q_0(\lambda,\mu) &= 2+2\lambda-\mu=\alpha+\lambda,\\
    \det Q_1(\lambda,\mu) &= (2+2\lambda-\mu)(2-\lambda-\mu)^2=(\alpha+\lambda)\beta^2,\\
    \det Q_n(\lambda,\mu) &= (\alpha\beta\gamma^2)^{3^{n-2}}\det Q_{n-1}\left(\frac{\lambda^2\beta}{\alpha\gamma},\mu+\frac{2\lambda^2\delta}{\alpha\gamma}\right)\qquad(n\ge2).
  \end{align*}
\end{lemma}

Consider now the quadratic forms
\[H_\theta = \mu^2-\lambda\mu-2\lambda^2-2-\mu+\theta\lambda,\]
and the function $F:[-4,5]\to[-4,5]$ given by $F(\theta) =
4-2\theta-\theta^2$.  Set $X_2 = \{-1\}$, and iteratively define $X_n
= F^{-1}(X_{n-1})$ for all $n\ge3$.  Note that $|X_n| = 2^{n-2}$.
\begin{lemma}
  We have for all $n\ge 2$ the factorization
  \[\det Q_n(\lambda,\mu)=(2+2\lambda-\mu)(2-\lambda-\mu)^{3^{n-1}+1}\prod_{\substack{2\le m\le n\\ \theta\in X_m}}H_\theta^{3^{n-m}+1}.\]
\end{lemma}

Thus, according to the previous proposition, the spectrum of $Q_n$ is
a collection of two lines and $2^{n-1}-1$ hyperbol\ae\ $H_\theta$ with
$\theta\in X_2\sqcup X_3\sqcup\dots\sqcup X_n$. The spectrum of
$\Delta_n$ is obtained by solving $\det Q_n(1,\mu)=0$, as given in
Figure~\ref{fig:fg} and the following theorem:
\begin{figure}[ht]
  \begin{center}
    \psfrag{mu}{$\mu$}
    \psfrag{lambda}{$\lambda$}
    \psfrag{\2612}{$-2$}
    \psfrag{\2611}{$-1$}
    \psfrag{0}{$0$}
    \psfrag{1}{$1$}
    \psfrag{2}{$2$}
    \psfrag{3}{$3$}
    \psfrag{4}{$4$}
    \includegraphics[width=100mm]{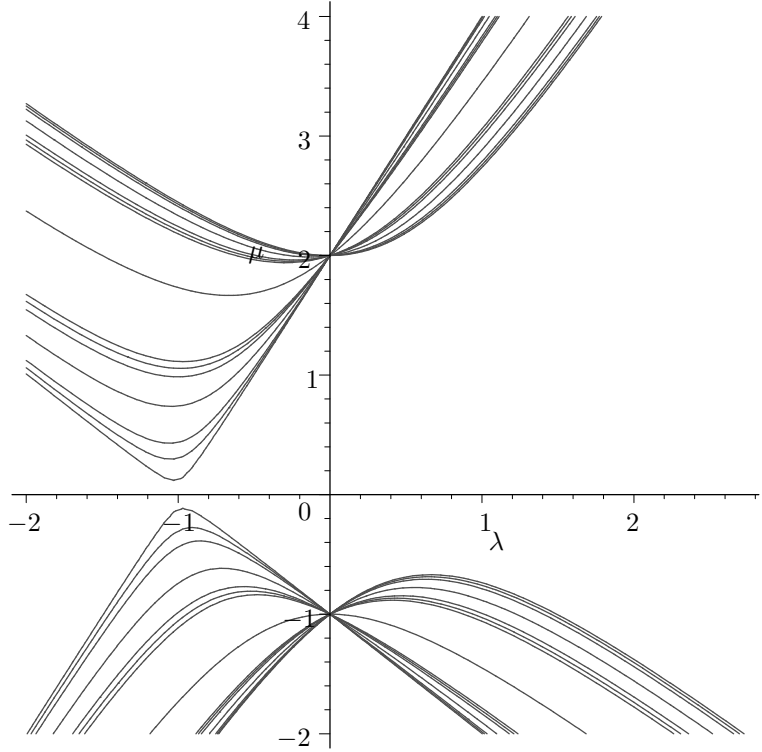}
  \end{center}
  \caption{The zeroes of $\det Q_n(\lambda,\mu)$, giving the spectrum of the
    Fabrykowski-Gupta group $G$}
  \label{fig:fg}
\end{figure}

\begin{theorem}
  Let $\pi_\pm:[-4,5]\to[-2,4]$ be defined by
  $\pi_\pm(\theta)=1\pm\sqrt{5-\theta}$. Then
  \begin{gather*}
    \spec{\Delta_0} = \{4\};\\
    \spec{\Delta_1} = \{1,4\};\\
    \spec{\Delta_n} = \{1,4\} \cup \bigcup_{2\le m\le n}\pi_\pm(X_m)\qquad(n\ge 2).
  \end{gather*}
  The spectrum of $\pi$ for the Fabrykowski-Gupta group $G$ is the
  closure of the set (here $\mu$ indicates the dimension of the
  eigenspace)
  \[
  \left\{\begin{array}{rl}
      4 & (\mu=0)\\
      1 & (\mu=\frac13)\\
      1\pm\sqrt6 & (\mu=\frac19)\\
      1\pm\sqrt{6\pm\sqrt6} & (\mu=\frac1{3^3})\\
      1\pm\sqrt{6\pm\sqrt{6\pm\sqrt6}} & (\mu=\frac1{3^4})\\
      \dots\end{array}\right\}.
  \]
  It is the union of a Cantor set of null Lebesgue measure that is
  symmetrical about $1$, and a countable collection of isolated points
  supporting the empiric spectral measure, which has the values
  indicated as $\mu$.
\end{theorem}

\section{Self-similarity and $C^*$-algebras}
\label{s:cstar}
\subsection{The self-similarity bimodule}
The self-similarity of the sets, group actions and inverse semigroups
can be treated from a common point of view through the notion of a
\emph{$C^*$-bimodule} or \emph{correspondence}. For the notions of
$C^*$-bimodules, correspondences and their applications
see~\cite[Appendix A in Chapter~II and Appendix B in
Chapter~V]{connes:noncomg} and the papers~\cite{conjon:t,subfactors}.

\begin{defi}
  Let $A$ be a $*$-algebra. An \emph{$A$-bimodule} $\Phi$ is a right
  $A$-module with an $A$-valued sesquilinear inner product and a
  $*$-homomorphism $\phi:A\to\mathrm{End}(\Phi)$.
\end{defi}

By an $A$-valued sesquilinear inner product on $\Phi$ we mean a
function $\scalar{\cdot}{\cdot}$ from $\Phi\times\Phi$ to $A$ such that
\begin{enumerate}
\item
$\scalar{v}{v_1+v_2}=\scalar{v}{v_1}+\scalar{v}{v_2}$,
\item
$\scalar{v_1}{v_2a}=\scalar{v_1}{v_2}a$,
\item
$\scalar{v_1}{v_2}=\left(\scalar{v_2}{v_1}\right)^*$,
\item
$\scalar{v}{v}$ is a positive element of $A$ and is equal to zero if and only if $v=0$,
\end{enumerate}
where $v, v_1, v_2\in\Phi$ and $a\in A$ are arbitrary.

For $a\in A$ and $v\in\Phi$ we will write $av$ instead of
$\phi(a)(v)$, so that the map $\phi$ defines the left multiplication
of the bimodule $\Phi$.

Let $\left(\{F_v\}_{v\in V}, \{\phi_e\}_{e\in E}\right)$ be a
graph-directed iterated function system.  Let $A=\bigoplus_{v\in
  V}C(F_v)$ be the direct sum of the $C^*$-algebras of continuous
$\C$-valued functions on the spaces $F_v$. The algebra $A$ can be also
defined as the algebra of continuous functions on the disjoint union
$\tilde F$ of the spaces $F_v$.

Let $\Phi_R$ be the right $A$-module which is the direct sum
$\bigoplus_{e\in E} e\otimes C(F_{\alpha(e)})$ (here $e\otimes\cdot$
is only an index), with the $A$-valued inner product
\[
\scalar{\sum_{e\in E}e\otimes f_e}{\sum_{e\in E}e\otimes h_e}=\sum_{e\in E} \overline{f_e} h_e,
\]
and the right multiplication defined by the rule
\[
\left(\sum_{e\in E}e\otimes f_e\right)\cdot h=\sum_{e\in E}e\otimes (f_e\cdot h|_{F_{\be(e)}}).
\]

We can also define the left multiplication by the formula
\[
h\cdot \left(\sum_{e\in E}e\otimes f_e\right)=\sum_{e\in E}e\otimes \left(h\circ\phi_e\cdot f_e\right),
\]
where $h\circ\phi_e$ is the composition $
h\circ\phi_e(x)=h|_{F_{\en(e)}}\left(\phi_e(x)\right)$.

In this way we get a well-defined \emph{self-similarity bimodule}
$\Phi$ over the algebra $A$ associated to the iterated function
system.

The self-similarity of the group and of inverse semigroup actions can
be also encoded in terms of bimodules.

Let us fix some self-similar action of a group $G$ on the space $\xo$
for the alphabet $X=\{\mathit{1}, \mathit{2}, \ldots, \mathit{d}\}$.
Denote by $\Phi_R$ the free right module over the group algebra $\C G$
with the free basis equal to the alphabet $X$.

The algebra $\C G$ is a $*$-algebra with the standard involution
$(\alpha\cdot g)^*=\overline{\alpha}g^{-1}$ (where $\alpha\in\C$,
$g\in G$). On the module $\Phi_R$ we define a $\C G$-valued
sesquilinear form by the equality
\[
\scalar{\sum_{x\in X}x\cdot a_x}{\sum_{y\in X} y\cdot b_y}=
\sum_{x\in X} a_x^*b_x,
\]
where $a_x, b_y\in\C G$.

Now the definition of self-similarity gives us a structure of a left
$\C G$-module on $\Phi_R$. Namely, for any element $g\in G$ and any
vector $x\in X$ from the basis of $\Phi_R$ we define
\begin{equation}
  \label{eq:stcl2}
  g\cdot x=y\cdot h,
\end{equation}
where $h$ and $y$ are as in Equation~\eqref{eq:ssim} of
Definition~\ref{defi:ssact2}. This multiplication is extended by
linearity to the whole module $\Phi_R$. In this way we get a map from
$G$ to $\mathrm{End}(\Phi_R)$.  This map extends to a morphism of
algebras $\phi:\C G\to\mathrm{End}(\Phi_R)$. The algebra
$\mathrm{End}(\Phi_R)$ is isomorphic to the algebra $M_d\left(\C
  G\right)$ of $d\times d$-matrices over $\C G$.

The morphism $\phi:\C G\to M_d\left(\C G\right)$ is called the
\emph{linear recursion} of the self-similar action.  The obtained
$\C G$-bimodule $\Phi$ is called the \emph{self-similarity} bimodule of
the action.

For instance, in the case of the adding machine action, the
self-similarity bimodule $\Phi$ is $2$-dimensional as a right
$\C\langle a\rangle$-module and  left multiplication by the
generator $a$ is the endomorphism of the right module defined by the
matrix
\[
\phi(a)=\left(
  \begin{array}{cc}
    0 & a \\
    1 & 0
  \end{array}\right).
\]

The linear recursions of self-similar actions were used by the first
two named authors and A.~\.Zuk to compute the spectra of random walks on the
Cayley graphs and on the Schreier graphs of some self-similar groups
(see~\cite{gr_zu:lamp,bgr:spec}). See Section~\ref{s:spectrum} for
details.

Suppose that $\rho$ is a unitary representation of the group $G$ on a
Hilbert space $\hilb$. Then we have a representation
$\rho_1=\Phi\otimes_{\C G}\rho$ of the group $G$ on the space
\[\Phi\otimes_{\C G}\hilb=\sum_{x\in X} x\otimes \hilb,\]
where each $x\otimes \hilb$ is a copy of $\hilb$ with the natural
isometry $T_x:\hilb\to x\otimes \hilb:v\mapsto x\otimes v$.  Then the
representation $\rho_1$ acts on the space $\Phi\otimes_{\C G}\hilb$ by
the formula $\rho_1(g)(x\otimes v)=y\otimes(\rho(h)(v)),$ where $g\in
G$, $v\in \hilb$, and $h\in G$, $y\in X$ are such that $g\cdot
x=y\cdot h$.

A unitary representation $\rho$ of the group $G$ is said to be
\emph{self-similar} if $\rho$ and $\Phi\otimes\rho$ are unitary
equivalent. Therefore, the representation $\rho$ is self-similar if
and only if there exit a decomposition of $\hilb$ into a direct sum
$\hilb=\hilb_{\mathit{1}}\oplus \hilb_{\mathit{2}}\oplus\cdots\oplus
\hilb_{\mathit{d}}$ and isometries $S_x:\hilb\to \hilb$ with range
$\hilb_x$, such that
\begin{equation}
  \rho(g)S_x=S_y\rho(h)
  \label{eq:phiinv}
\end{equation}
whenever $g\cdot x=y\cdot h$.

\paragraph{Examples.} \begin{enumerate}
\item Since the group $G$ acts on the space $\xo$ by
  measure-preserving transformations, we get a natural unitary
  representation $\rho$ of the group $G$ on the space
  $\hilb=L^2(\xo)$. Since the set $\xo$ is the disjoint union
  $\cup_{x\in X} x\xo$, the space $\hilb$ is the direct sum of the
  spaces $\hilb_x$ of functions with support in $x\xo$. We have a
  natural isometry $S_x:\hilb\to \hilb_x\subset \hilb$ defined by the
  rule:
  \[
  S_x(f)(w)=\left\{\begin{array}{cl}
      0     & \text{if $w\notin x\xo$},\\
      \sqrt{|X|}\cdot f(w') & \text{if $w=xw'$}.
    \end{array}
  \right.
  \]
  It is checked directly that condition~\eqref{eq:phiinv} holds, so
  the natural representation of $G$ on the space $L^2(\xo)$ is
  self-similar.
\item A set $\mathcal{M}\subseteq \xo$ is \emph{self-similar} if
  $\mathcal{M}=\cup_{x\in X}x\mathcal{M}$ is a disjoint union.
  
  Let $\mathcal{M}\subset\xo$ be a countable self-similar
  $G$-invariant set.  Set $\hilb=\ell^2(\mathcal{M})$. Then the space
  $\hilb$ is also a direct sum $\hilb_{\mathit{1}}\oplus
  \hilb_{\mathit{2}}\oplus\cdots\oplus \hilb_{\mathit{d}}$, where
  $\hilb_x$ is the subspace spanned by the set $x\mathcal{M}$. The
  isometry $S_x:\hilb\to \hilb_x$ is the linear extension of the map
  $w\mapsto xw$.  Since the set $\mathcal{M}$ is $G$-invariant, we
  have a natural permutational representation of the group $G$ on the
  space $\hilb$. It is also easy to check that this representation
  satisfies condition~\eqref{eq:phiinv}, thus is self-similar.
\end{enumerate}

Suppose that the representation $\rho$ of the group $G$ is
self-similar and let $S_x$ be the isometries for
which~\eqref{eq:phiinv} holds.
  
Let $C^*_\rho(G)$ be the completion of the group algebra $\C G$ with
respect to the operator norm induced by the representation $\rho$. Let
$\Phi_\rho$ be the right $C^*_\rho(G)$-module with free basis $X$.
Then formula~\eqref{eq:stcl2} gives a well-defined
$C^*_\rho(G)$-bimodule structure on $\Phi_\rho$.

Therefore, if the representation $\rho$ is self-similar, then the
$\C G$-bimodule $\Phi$ can be extended to a $C^*_\rho(G)$-bimodule.

In general, we adopt the following definition:
\begin{defi}
  A completion $A$ of the algebra $\C G$ respectively to some
  $C^*$-norm is said to be \emph{self-similar} if the self-similarity
  bimodule $\Phi$ extends to an $A$-bimodule.
\end{defi}

We need the following auxiliary notion:
\begin{defi}
  Let $G$ be a countable group acting by homeomorphisms on the set
  $\xo$. A point $w\in\xo$ is \emph{generic} with respect to $g\in G$
  if either $w^g\neq w$ or there exists a neighborhood $U\ni w$
  consisting of the points fixed under the action of $g$.
  
  A point $w\in\xo$ is \emph{$G$-generic} if it is generic with
  respect to every element of $G$.
\end{defi}

It is easy to see that for any countable group $G$ acting by
homeomorphisms on the set $\xo$, the set of all $G$-generic points is
residual, i.e., is an intersection of a countable set of open dense
subsets. In particular, it is non-empty.

For any point $w\in\xo$ we denote by $G(w)$ the $G$-orbit of $w$. Let
$\ell^2(G(w))$ be the Hilbert space of all square-summable functions
$G(w)\to\C$. We have the permutation representation $\pi_w$ of $G$ on
$\ell^2(G(w))$. Let $\|\cdot\|_w$ be the operator norm on $\C G$
defined by the representation $\pi_w$.

\begin{proposition}
  \label{pr:gen_norm}
  Let $w_1, w_2\in\xo$ and suppose that $w_1$ is $G$-generic.  Then
  for every $a\in\C G$
  \[\|a\|_{w_1}\leq \|a\|_{w_2}.\]
\end{proposition}

By $A_w$ we denote the completion of the algebra $\C G$ with respect to
the operator norm defined by the permutation representation of $G$ on
the orbit of the point $w\in\xo$. As an immediate corollary of
Proposition~\ref{pr:gen_norm}, we get

\begin{theorem}
  \label{th:gen_rep}
  Let $w\in\xo$ be a $G$-generic point. Then for every $u\in\xo$, the
  algebra $A_w$ is a quotient of the algebra $A_u$. If $u$ is also
  generic, then the algebras $A_w$ and $A_u$ are isomorphic.
\end{theorem}

Let us denote the algebra $A_w$, where $w$ is a generic point, by
$\redu$.

\begin{theorem}
  \label{th:gen_rec}
  The algebra $\redu$ is self-similar. If $A$ is another self-similar
  completion of $\C G$ then the identity map $G\to G$ extends to a
  surjective homomorphism of $C^*$-algebras $A\to\redu$.
\end{theorem}

\subsection{Cuntz-Pimsner algebras}
In~\cite{pimsner} M.~Pimsner associates to every bimodule $\Phi$ an
algebra $\OPhi$, which generalizes the Cuntz algebra $\On$
(see~\cite{cuntz}), and the Cuntz-Krieger algebra $\mathcal{O}_A$
from~\cite{cuntzkrie}. In the case of the self-similarity bimodule the
Cuntz-Pimsner algebra $\OPhi$ can be defined in the following way.

\begin{defi}
  \label{defi:cp}
  Let $\Phi$ be the self-similarity bimodule of the action of a group
  $G$ over the alphabet $X$.  The \emph{Cuntz-Pimsner} algebra $\OPhi$
  is the universal $C^*$-algebra generated by the algebra $\redu$ and
  $|X|$ operators $\{S_x:x\in X\}$ satisfying the relations
  \begin{gather}
    S_x^*S_x=1, \quad S_x^*S_y=0\text{ if }x\neq y,\label{eq:cp1}\\
    \sum_{x\in X}S_xS_x^*=1,\label{eq:cp2}\\
    aS_x=\sum_{y\in X}S_ya_{x, y}\label{eq:cp3},
  \end{gather}
  where $a_{x, y}\in\redu$ are such that $a\cdot x=\sum_{y\in X}
  y\cdot a_{x, y}$, so $a_{x, y}=\scalar{y}{a\cdot x}$.
\end{defi}

We will use the multi-index notation, so that for $v=x_1x_2\ldots
x_n\in X^*$ the operator $S_v$ is equal to $S_{x_1}S_{x_2}\cdots
S_{x_n}$ and for the empty word $\emp$ the operator $S_\emp$ is equal
to 1. Then, since $g\cdot v=u\cdot h$, we have $gS_v=S_uh$ in the
algebra $\OPhi$.

The \emph{Cuntz algebra} $\On$ is the universal algebra generated by
$d$ isometries $S_1, S_2, \ldots, S_d$ such that $\sum_{i=1}^d
S_iS_i^*=1$.  This algebra is simple (see~\cite{cuntz,davidson}) and
thus any isometries satisfying such a relation generate an algebra
isomorphic to $\On$. In particular, since the operators $\{S_x: x\in
X\}$ satisfy relations~\eqref{eq:cp1},~\eqref{eq:cp2}, they generate a
subalgebra of $\OPhi$, isomorphic to $\On$ for $d=|X|$.

If $\rho:\OPhi\to\mathcal{B}(\hilb)$ is a representation of the
algebra $\OPhi$, then its restriction onto the subalgebra generated by
$\redu$ is a self-similar representation of the algebra $\redu$,
by~\eqref{eq:cp3}. Conversely, if $\rho$ is a self-similar
representation, then by~\eqref{eq:phiinv} we get a representation of
the algebra $\OPhi$.

\begin{theorem}
  \label{th:ogn}
  The algebra $\OPhi$ is simple. Moreover, for any non-zero
  $x\in\OPhi$ there exist $p, q\in\OPhi$ such that $pxq=1$.
\end{theorem}

The Cuntz-Pimsner algebras of self-similarity bimodules of iterated
function systems were studied by V.~Deaconu (in the more general
situation of a \emph{continuous graph}) in~\cite{deac:solen}.  See
also the work~\cite{putnam}, where the \emph{Ruelle algebras} of a
\emph{Smale space} are studied. The connection between the Ruelle
algebras and the Cuntz-Pimsner algebras of topological graphs is
described in~\cite{ar}.

\subsection{The algebra $\UPhi$}
For $k\in\N$ denote by $\mathcal{F}_k$ the linear span in $\OPhi$ over
$\C$ of the products $S_vaS_u^*$, with $u, v\in X^k$ and $a\in\redu$.

The space $\mathcal{F}_k$ is an algebra isomorphic to the algebra
$M_{d^k}(\redu)$ of $d^k\times d^k$-matrices over the algebra $\redu$.

By the linear recursion we have inclusions
$\mathcal{F}_k\subset\mathcal{F}_{k+1}$. Denote by $\UPhi$ the closure
of the union $\cup_{k\geq 1}\mathcal{F}_{k}$.

The algebra $\UPhi$ is isomorphic to the direct limit of the matrix
algebras $M_{d^n}(\redu)$ with respect to the embeddings defined by
the linear recursion $\phi$.

\begin{theorem}
  \label{th:cross}
  Suppose the self-similar action of $G$ is recurrent and that the
  point $w_0\in\xo$ and all its shifts $\shift^n(w_0)$ are
  $G$-generic.
  
  Let $\pi_1$ be the permutation representation of $G$ on the orbit of
  $w_0$ and let $\pi_2$ be the natural representation on
  $\ell^2(G(w_0))$ of the algebra $C(\xo)$ of continuous functions on
  $\xo$. Then the algebra $\UPhi$ is isomorphic to the $C^*$-algebra
  generated by $\pi_1(G)\cup\pi_2(C(\xo))\subset
  \mathcal{B}\left(\ell^2\left(G\left(w_0\right)\right)\right)$.
\end{theorem}

\paragraph{Example.} In the case of the adding machine action the algebra
$\redu$ is isomorphic to the algebra $C(\mathbb{T})$ with the linear
recursion $C(\mathbb{T})\to M_2(C(\mathbb{T}))$ coming from the double
self-covering of the circle.  Thus the algebra $\UPhi$ in this case is
the Bunce-Deddence algebra. Then Theorem~\ref{th:cross} in this case
is the well known fact that the Bunce-Deddence algebra is the
cross-product algebra of the odometer action on the Cantor space $\xo$
(see~\cite{davidson}).


\bibliographystyle{amsalpha}
\bibliography{mymath}
\end{document}